\documentclass[12pt]{amsart}
\begin{document}
\theoremstyle{plain}
\newtheorem{thm}{Theorem}[section]
\newtheorem*{thm*}{Theorem}
\newtheorem{prop}[thm]{Proposition}
\newtheorem*{prop*}{Proposition}
\newtheorem{lemma}[thm]{Lemma}
\newtheorem{cor}[thm]{Corollary}
\newtheorem*{conj*}{Conjecture}
\newtheorem*{cor*}{Corollary}
\newtheorem{defn}[thm]{Definition}
\theoremstyle{definition}
\newtheorem*{defn*}{Definition}
\newtheorem{rems}[thm]{Remarks}
\newtheorem*{irem*}{Important Remark}
\newtheorem*{rems*}{Remarks}
\newtheorem*{proof*}{Proof}
\newtheorem*{not*}{Notation}
\newcommand{\npartial}{\slash\!\!\!\partial}
\newcommand{\Heis}{\operatorname{Heis}}
\newcommand{\Solv}{\operatorname{Solv}}
\newcommand{\Spin}{\operatorname{Spin}}
\newcommand{\SO}{\operatorname{SO}}
\newcommand{\ind}{\operatorname{ind}}
\newcommand{\Index}{\operatorname{index}}
\newcommand{\ch}{\operatorname{ch}}
\newcommand{\rank}{\operatorname{rank}}

\newcommand{\abs}[1]{\lvert#1\rvert}

\newcommand{\blankbox}[2]{%
  \parbox{\columnwidth}{\centering
    \setlength{\fboxsep}{0pt}%
    \fbox{\raisebox{0pt}[#2]{\hspace{#1}}}%
  }%
}
\newcommand{\supp}[1]{\operatorname{#1}}
\newcommand{\norm}[1]{\parallel\, #1\, \parallel}
\newcommand{\ip}[2]{\langle #1,#2\rangle}
\setlength{\parskip}{.3cm}
\newcommand{\nc}{\newcommand}
\nc{\nt}{\newtheorem}
\nc{\gf}[2]{\genfrac{}{}{0pt}{}{#1}{#2}}
\nc{\mb}[1]{{\mbox{$ #1 $}}}
\nc{\real}{{\mathbb R}}
\nc{\comp}{{\mathbb C}}
\nc{\ints}{{\mathbb Z}}
\nc{\Ltoo}{\mb{L^2({\mathbf H})}}
\nc{\rtoo}{\mb{{\mathbf R}^2}}
\nc{\slr}{{\mathbf {SL}}(2,\real)}
\nc{\slz}{{\mathbf {SL}}(2,\ints)}
\nc{\su}{{\mathbf {SU}}(1,1)}
\nc{\so}{{\mathbf {SO}}}
\nc{\hyp}{{\mathbb H}}
\nc{\disc}{{\mathbf D}}
\nc{\torus}{{\mathbb T}}
\newcommand{\tk}{\widetilde{K}}
\newcommand{\boe}{{\bf e}}\newcommand{\bt}{{\bf t}}
\newcommand{\vth}{\vartheta}
\newcommand{\CGh}{\widetilde{\CG}}
\newcommand{\db}{\overline{\partial}}
\newcommand{\tE}{\widetilde{E}}
\newcommand{\tr}{\mbox{tr}}
\newcommand{\ta}{\widetilde{\alpha}}
\newcommand{\tb}{\widetilde{\beta}}
\newcommand{\txi}{\widetilde{\xi}}
\newcommand{\hV}{\hat{V}}
\newcommand{\IC}{\mathbf{C}}
\newcommand{\IZ}{\mathbf{Z}}
\newcommand{\IP}{\mathbf{P}}
\newcommand{\IR}{\mathbf{R}}
\newcommand{\IH}{\mathbf{H}}
\newcommand{\IG}{\mathbf{G}}
\newcommand{\CC}{{\mathcal C}}
\newcommand{\CD}{{\mathcal D}}
\newcommand{\CS}{{\mathcal S}}
\newcommand{\CG}{{\mathcal G}}
\newcommand{\CL}{{\mathcal L}}
\newcommand{\CO}{{\mathcal O}}
\nc{\ca}{{\mathcal A}}
\nc{\cag}{{{\mathcal A}^\Gamma}}
\nc{\cg}{{\mathcal G}}
\nc{\chh}{{\mathcal H}}
\nc{\ck}{{\mathcal B}}
\nc{\cl}{{\mathcal L}}
\nc{\cm}{{\mathcal M}}
\nc{\cn}{{\mathcal N}}
\nc{\cs}{{\mathcal S}}
\nc{\cz}{{\mathcal Z}}
\nc{\sind}{\sigma{\rm -ind}}


\title{}

\titlepage

\begin{center}{\bf SPECTRAL FLOW IN FREDHOLM
MODULES, ETA INVARIANTS AND THE JLO COCYCLE}
\footnote{
spectral flow, $\theta$-summable Fredholm module, eta invariant, index.}
\\
\vspace{.75 in}
{\bf by}\\
\vspace{.75 in}
{\bf Alan Carey}\\Mathematical Sciences Institute\\Australian National
University\\
Canberra A.C.T.\\AUSTRALIA\\
\vspace{.25 in}
{\bf and}\\
\vspace{.25 in}
{\bf John Phillips}\\Department of Mathematics and Statistics\\
University of Victoria\\Victoria, B.C. V8W 3P4\\CANADA\\
\vspace{1.5 in}
Supported by grants from ARC (Australia) and
NSERC (Canada)
\vspace{.5 in}

\end{center}

\begin{abstract}
We give a comprehensive account of an analytic approach to spectral flow 
along paths of self-adjoint Breuer-Fredholm operators
in a type $I_{\infty}$ or $II_\infty$ von Neumann algebra ${\mathcal N}$.
The framework is that of
{\it odd unbounded} $\theta$-{\it summable} {\it Breuer-Fredholm modules} for a
unital Banach $*$-algebra, $\mathcal A$. 
In the type $II_{\infty}$ case spectral flow is real-valued,
has no topological definition as an intersection number
and our formulae encompass all that is known.
We borrow Ezra Getzler's idea (suggested by I. M. Singer) of considering 
spectral flow (and eta invariants) as the integral of a closed one-form on
an affine space. Applications in both the type I and type II cases
include a general formula for the relative index
of two projections,
representing truncated eta functions as integrals of one forms and
expressing spectral flow in terms of the JLO cocycle
to give the pairing of the $K$-homology and $K$-theory 
of $\mathcal A$.

\end{abstract}

\maketitle

\section{\bf INTRODUCTION}

The notion of spectral flow has been
a useful analytic tool in geometry ever since its invention by
Atiyah and Lusztig \cite{APS1,APS3, BW}.
 Motivated by observations of I.M. Singer \cite{Si}
on eta invariants which suggest
that spectral flow should be expressible as the integral
of a one-form, there has been a succession of contributions in
\cite{DHK},\cite{H},\cite{Kam},\cite{G},\cite{Ph, Ph1} and \cite{P1,P2}.
A key step in synthesising these developments was taken in
\cite {CP1} where we exploited previous work by one of us \cite{Ph,Ph1}
and ideas of \cite{G} to produce spectral flow formulae as
integrals of one-forms on affine spaces arising from finitely
summable Fredholm modules. Significantly, we also established spectral
flow formulae
along paths of self-adjoint
Breuer-Fredholm operators in a type $II_\infty$ von Neumann algebra,
a development that was first hinted at in \cite{APS3}.
This has special relevance to recent developments in the study
of $L^2$ spectral invariants for manifolds whose fundamental group
has a non-type $I$ regular representation, see \cite{Ma}
for a review of these ideas.
The casual reader of say \cite{G} might be puzzled as to the
reasons for our lengthy treatment. They are threefold,
 first of all we establish much more general
formulae than are described in \cite{G}. Second, in the type $II_\infty$
case there is no topological definition of spectral flow
as an intersection number
(partly due to the fact that 
Breuer-Fredholm operators may have zero in their continuous spectrum).
Third, the use of the Duhamel Principle in \cite{G} requires strong
domain assumptions (which are not stated)
for the unbounded operators that appear there. 
Inserting these assumptions then excludes
interesting examples and introduces ugly 
 unnecessary complications. We adopt new methods
which entail analytic
subtleties in both the type $I_\infty$ and type $II_\infty$ case
but lead to analytic results of wider interest
and applicability and considerably extend the existing literature 
even in the classical type I case (for example, our
formula for the relative index of two projections in Section 3).

The present paper also goes much further in
explaining these earlier developments
in terms of a single general formula for spectral flow along paths
of self-adjoint Breuer-Fredholm operators in $\theta$-summable
Breuer-Fredholm modules. Moreover our proofs are
the same for
both type $I_\infty$ and type $II_\infty$ von Neumann algebras.
Our analytic approach has already found application in
the type $I_\infty$ case in \cite{BF} and
 in proving index theorems for generalised Toeplitz operators \cite{CPSu}.
 In this paper our main application is on the relationship 
with the JLO cocycle (and hence the pairing of the $K$-homology of 
$\mathcal A$ and the odd $K$-theory of $\mathcal A$).

To explain our results we need some notation and definitions.
We denote by $\mathcal N$ a type $I_\infty$ or type $II_\infty$ von Neumann
algebra acting on a Hilbert space $H$ with faithful
normal semifinite trace $\tau$.
The ideal of ``compact operators'' in $\mathcal N$
is denoted $\mathcal K_{\cn}$ (see \cite {B1, B2}
for this and also the notion of Fredholm operator in the type $II_\infty$
setting).
We fix an unbounded self-adjoint operator $D_0$
on $H$ affiliated with
$\mathcal N$ and assume that $A$ is a self-adjoint element in $\mathcal N$
(the set of all such being denoted ${\mathcal N}_{sa}$).
We consider paths of the form $D(t)=D_0+A(t)$ where $t$ is a real parameter
and $A(t)\in {\cn}_{sa}$ for each $t$.

\begin{defn*}
 We say
that $(\mathcal N,D_0)$ is an {\bf unbounded $\theta$-summable
Breuer-Fredholm module}
for a Banach $*$-algebra $\mathcal A$ if $\mathcal A$ is represented in
$\mathcal N$ and if $e^{-tD^2_0}$ is trace-class for all $t>0$
and $[D_0,a]$ is bounded for all $a$ in a dense $*$-subalgebra of $\mathcal A$.
\end{defn*}

We let
$${\mathcal M}_0= \{D=D_0+A\ |\ A\in {\mathcal N}_{sa}\}$$
Clearly ${\mathcal M}_0$  is an affine space modelled on  ${\mathcal N}_{sa}$.
Introduce the \lq gauge group' $\mathcal G$ defined by
$${\mathcal G} = \{U\in{\mathcal N}\ |\ U\text{ is unitary},
 [D_0,U]\text{ is bounded}\}.$$

Let $\gamma=\{D_t=D_0+A(t), a\leq t\leq b\}$ be a
piecewise $C^1$ path in  ${\mathcal M}_0$ with
$D_a$ and $D_b$ invertible.
The spectral flow formula
of \cite{G} when ${\mathcal N}=B(H)$  is
\begin{eqnarray}
\text{sf}(D_a,D_b)&=& -\int_{\gamma}\alpha_\epsilon +
\frac{1}{2}\eta_\epsilon(D_b)-\frac{1}{2}\eta_\epsilon(D_a)\nonumber\\
&=&\sqrt{\frac{\epsilon}{\pi}}\int_a^b\tau\left(\frac{d}{dt}(D_t)
e^{-\epsilon D_t^2}\right)dt+
\frac{1}{2}\eta_\epsilon(D_b)-\frac{1}{2}\eta_\epsilon(D_a).\nonumber
\end{eqnarray}
where $\eta_\epsilon(D)$
are approximate eta invariant correction terms (we define these later).
The difference vanishes if the endpoints are in the same gauge
group orbit.

Our aim is to establish not only this formula
 but a much more general one for the spectral flow of
bounded self-adjoint Breuer-Fredholm operators.
We will use this to establish
formulae for spectral flow along paths in ${\mathcal M}_0$.
To describe the bounded case we need some further
notation and definitions.

\begin{defn*}
A bounded, {\bf odd Breuer-Fredholm module}
for a unital Banach
$*$-algebra $\mathcal A$ represented in $\mathcal N$
is a pair $({\mathcal N}, F)$ with  $F$ a self-adjoint operator
in $\mathcal N$ such that $F^2=1$ and $[F,a] \in \mathcal K_{\cn}$ for all
$a\in \mathcal A$.
\end{defn*}

In this introduction let us restrict to the type $I_\infty$ case and introduce
the two-sided ideal of operators $Li_0(H)$ consisting of those
compact operators $T$
whose $n^{th}$ singular value is $o((\log n)^{-1})$ as $n\to \infty$.
Then we say $(H,F)$ is $\theta$-summable if $[F,a]\in Li_0^{1/2}(H)$
for a dense set of $a\in {\mathcal A}$. The choice $F=sign(D_0)$ relates unbounded to
bounded Fredholm modules.
If we replace the ideal  $Li_0(H)$ by the ideal $Li(H)$
  consisting of compact operators $T$ whose $n$th
singular value is $O((\log n)^{-1})$,
then this defines
a weakly $\theta$-summable Fredholm module. (Connes' most recent definition
of $\theta$-summable is what we call weakly $\theta$-summable \cite{Co4}.)
We prove similar spectral flow formulae in each of these cases.

$\spadesuit$ Our most general formula in the bounded case deals with a pair 
of self-adjoint Fredholm operators
$\{F_j,\;\;\ j=1,2\}$, joined by
a piecewise $C^1$ path $\{F_t\}$, $t\in [1,2]$
in a certain affine subspace (specified in terms of $Li_0(H)$) of the
space of all self-adjoint Fredholm operators.
The spectral flow along
such a path
 is given by
$$sf(F_1,F_2)=\frac{1}{C}
\int_1^2\tau\left(\frac{d}{dt}(F_t)|1-F_t^2|^{-r}
e^{-|1-F_t^2|^{-\sigma}}\right)dt
+\gamma(F_2)-\gamma(F_1) $$
where the $\gamma(F_j)$ are eta invariant type correction terms and $C$ is a
normalization constant depending on the parameters $r\geq 0$ and $\sigma\geq1$.
To see one important place where such complicated formulae arise, one takes 
the Getzler expression with $\epsilon=1$: 

$$\frac{1}{\sqrt{\pi}}\int_a^b\tau\left(\frac{d}{dt}(D_t)
e^{-D_t^2}\right)dt$$

\noindent and does the change of variable $F_t=D_t(1+D_t^2)^{-1/2}$. Then,
$(1+D_t^2)^{-1}=(1-F_t^2)=|1-F_t^2|$, and if one is 
careless and just differentiates formally (not worrying about the order of the
factors), one obtains the expression:

$$\frac{e}{\sqrt{\pi}}\int_a^b\tau\left(\frac{d}{dt}(F_t)|1-F_t^2|^{-3/2}
e^{|1-F_t^2|^{-1}}\right)dt.$$

While the actual details are much more complicated, this is the heuristic 
essence of our reduction of the unbounded case to the bounded case: see
Propositions 6.5 and 6.6.
$\spadesuit$

With the extra flexibility afforded by the parameters $r\geq 0$ and
$\sigma\geq 1$, we also show that this same formula, given a 
type $II_\infty$ analogue of the ideal $Li(H)$, holds for spectral flow 
along a path in an affine subspace of bounded
self-adjoint Breuer-Fredholm operators associated with a $\theta$-summable
Breuer-Fredholm module.

Our approach is very general; formulae
studied elsewhere follow from it
(with the proviso that side conditions are needed
in some cases).

The plan of the paper is to relegate many technical
functional analytic issues to appendices. This is not to say these
results are not in themselves of interest, rather that they
could detract from the flow of the main arguments.

We begin in Section 2 by laying out all of our definitions and assumptions
and where appropriate indicating how they relate to the existing literature.

Section 3 contains a 
formula for the essential codimension (or relative index)
of two projections (Theorem 3.1). Such
results have a long history and we believe our formula is the most
general possible. The relevance of this to the notion of spectral flow in
Fredholm modules is the following. Given a (bounded) Fredholm module
$(H,F_0)$ for a Banach $*$-algebra $\mathcal A$ and a unitary $u\in \mathcal A$,
then $P=\frac{1}{2}(F_0+1)$ and $Q=uPu^*$ are two projections with the property
that the operator $QP: P(H) \to Q(H)$ is a Fredholm operator whose
(relative) index equals the spectral flow of the straight line path from
$F_0$ to $uF_0u^*$. Depending on the summability flavour of the Fredholm module
we are able to obtain explicit integral formulae for this index in Section 4
based on these results of Section 3.

Thus, Section 4 contains the first of our spectral flow formulae (eg.,
Theorems 4.1 and 4.2).
We single these out because they are elegant, their proofs are relatively
short and ultimately the proofs of all of our later formulae are based
on them.

In Section 5 we show that the integral formulae of Section 4 make sense in the
more general context of pre-Fredholm modules where $1-F_0^2$ is not $0$
but merely compact of some summability flavour (these arise naturally from
transforming {\bf unbounded} Fredholm modules). We show
that the integrals involved are actually integrals of exact one-forms on
appropriate affine spaces of the form: $F_0+{\mathcal I_{F_0}}$ where
$\mathcal I_{F_0}$ is a certain subspace of compact operators depending on the
the summability flavour of the module. The exactness of these one-forms
is a difficult analytic problem involving Cauchy integrals along unbounded
contours. This problem is a key difficult step in our approach,
but in the text we focus on the ideas, relegating most of the technical issues
to appendix C. 
 For arbitrary piecewise $C^1$ paths in the space
$F_0+{\mathcal I_{F_0}}$ from
say $F_a$ to $F_b$ we must tack on one extra path at each end so that our new
path runs between the corresponding symmetries
$\tilde F_a=sign(F_a)$ to $\tilde F_b=sign(F_b)$.
We then invoke the invariance of the integral of an exact one-form and the
results of Section 4 to obtain our formulae (eg., Theorems 5.7 and 5.9).
Thus the correction terms arise
naturally as the integrals of the one-forms from an operator $F$ in
$F_0+{\mathcal I_{F_0}}$ to its associated symmetry $\tilde F$.

In Section 6 we show that the transformation from unbounded
$\theta$-summable modules $(H,D)$ to bounded modules $(H,F)$ via the map
$F_D=D(1+D^2)^{-1}$ carries $C^1$ paths to $C^1$ paths
provided we consider $\theta_q$-summable modules in the bounded case (where
$0<q<1$). The integral formulae are shown to transform as expected, thanks to
the commutativity property of the trace (Proposition 6.6).

In Section 7 we then use this transformation to obtain our formulae in the
unbounded case from those in the bounded case (Section 5). The direct
translations of these formulae involve the parameter $q$: we obtain
our final versions by taking the limit as $q\to 1$ (Theorem 7.8 and
Corollary 7.10).

In Section 8 we shed further light on the idea that
the eta invariant is the integral of a one-form. We
see that the correction terms in our spectral flow formula correspond to
integrating our one-form along a particular path. These correction terms are
seen to be truncated
eta invariants by using path independence of the integral 
of our exact one form and 
integrating along a different path
joining the same endpoints (with a discontinuity at one end).
In the type $I_\infty$ setting
considered by \cite{G} where he assumed that the endpoints of his path,
$D_a$ and $D_b$ are invertible, we show that our correction terms are
identical with his. In the general case (Theorem 8.9
and Corollary 8.10) we show how to modify the eta terms
to remove the invertibility assumption.

In Section 9 we use the Laplace Transform and the results of
Section 7 to give a ``best possible'' finitely summable unbounded version
of the spectral flow formula (Theorem 9.3). This should be compared with
Theorem 2.16 of \cite{CP1}. In a later paper \cite{CPSu}, this result is 
crucial in proving Connes' Dixmier-trace formula for the odd index in the
general setting of $(1,\infty)$ Breuer-Fredholm modules.

In Section 10 we explain how our formulae lead to the JLO cocycle 
by generalising the main theorem of \cite{G}. Our formula for type II spectral
flow in terms of the JLO cocycle may be used as the starting point
for a proof of the local index formula of Connes and Moscovici 
\cite{CoMo}in the setting 
of Breuer-Fredholm modules. This is a lengthy matter however and
we leave it to another place.

We reiterate that our aim has been to create a usable
theory of spectral flow in 
a type $II_\infty$ von Neumann algebra. In a separate paper \cite{CPSu}
we show that our formulae imply
index theorems such as those of \cite{CDSS}
bringing them into the fold of noncommutative geometry.

{\bf ACKNOWLEDGEMENTS} The authors would like to thank their many colleagues
for interesting and useful discussions on the topics of this paper.
In particular Chris Bose, Varghese Mathai, Adam Rennie 
and Fyodor Sukochev have
dispelled our ignorance of many aspects of the background material.

\section {\bf DEFINITIONS AND NOTATIONS}

Throughout this paper, $\mathcal N$ will denote a semifinite
von Neumann algebra (with separable predual) and $\tau$ will denote
a faithful, normal semifinite trace on
$\mathcal N.$ The norm-closed two-sided ideal in $\mathcal N$
generated by the elements of finite trace, will be denoted by
$\mathcal K_{\cn}.$ We will be concerned with certain normed ideals
$\mathcal I$
contained in $\mathcal K_{\cn}$ and which are best defined in terms
of $generalized\; singular\; values$ a notion due to Fack and Kosaki
(and others), see \cite{FK}.

\subsection{\bf Singular Values}

\begin{defn}
If $S\in \mathcal {N}$ and $t>0$, then the {\bf t-th (generalized)
singular value of S} is given by
$$\mu_t(S)=inf\{||SE||\:|\:E\:is\:a\:projection\:in\:\mathcal{N}\:with\:
\tau(1-E)\leq t\}.$$
Although we refer to \cite{FK} for the properties of these singular
values, we note that we are restricting ourselves only to bounded
operators (in $\mathcal N$). Hence,  we have that $0\leq \mu_t(S)
\leq ||S||$ for all $t>0$ and that for $S\in \mathcal {K_N}$ we
have $\mu_t(S)\to 0.$ Thus, for us it is reasonable to set $\mu_0=||S||.$
\end{defn}

\subsection{\bf Operator Ideals}

$\spadesuit$ For most of the paper we would like to concern ourselves with the ideals,
$Li$ and $Li_0$ defined below. However, when we transform from the
setting of unbounded modules to bounded modules we lose a little
control and are forced to consider powers of these ideals, $Li^q$
and $Li_0^q$ for $0<q\leq 1.$ In the end, we are able to rid
ourselves of these irritating exponents in our formulae by taking a limit
as $q\to 1.$ We observe that $Li_0\subset Li\subset Li_0^q\subset Li^q$ for all 
$q$ with $0<q<1$. $\spadesuit$

\begin{defn}
We define $Li=\{T\in \mathcal{N}\:|\:\mu_t(T)=O(1/\log{t})\}.$ The norm
on this ideal is $$||T||_{Li}=\sup_{x>0}\left\{\frac{\int_0^x \mu_t(T)dt}
{\int_0^x (\log{(t+e)})^{-1}dt} \right\}.$$ We note that $||T||_{Li} \geq ||T||.$

The ideal $Li_0$ is the closed subspace of $Li$ in the norm $||.||_{Li}$
of those operators $T\in \mathcal N$ satisfying $\mu_t(T)=o(1/\log{t}).$

For $0<q\leq 1,$ we consider also the powers of these ideals
$Li^q\supset Li_0^q$ with the norm
$$||T||_{Li^q}=\left(||\;|T|^{1/q}\;||_{Li}\right)^q.$$
For more on these ideals we refer the reader to appendix A.
\end{defn}

We also consider the ideals of finitely summable operators, $L^p,$
which we define for $p\geq 1$ by $L^p=\{T\in \mathcal{N}\:|\:
\tau(|T|^p)<\infty\}.$ Since the trace of a positive operator
is expressible in terms of singular values \cite{FK}, we could also
define these ideals in terms of singular values via:
$$L^p=\{T\in \mathcal{N}\:| \int_0^\infty \mu_t(T)^pdt < \infty\}.$$
The norm on $L^p$ is $$||T||_p=max\{||T||\:,\: (\tau(|T|^p))^{1/p}\}.$$
We observe that if $\mathcal N$ is type $I$ then we can omit the operator
norm on the right hand side given the usual normalization of $\tau.$
In the type $II$ case, our ideal $L^p$ is strictly contained in the space
of $p$-summable measurable operators \cite{FK}.

\subsection{\bf Breuer-Fredholm Modules}

\begin{defn}
An {\bf odd pre-Breuer-Fredholm module} for a unital Banach $*$-algebra
$\mathcal A$  is a pair $({\mathcal N},F_0)$ where $\mathcal A$ is
(continuously) represented in $\mathcal N$ and $F_0$ is a self-adjoint
Breuer-Fredholm operator in $\mathcal N$ satisfying:

$1.\: 1-F_0^2 \in {\mathcal {K_N}},\:and$

$2.\: [F_0,a] \in {\mathcal {K_N}}\: for\: a \in {\mathcal A}.$

\noindent If $1-F_0^2=0$ we drop the prefix "pre-".

\noindent If, in addition, our module satisfies:

$1.'\: 1-F_0^2\in Li_0$ (respectively, $Li$; $Li_0^q$ for $0<q\leq 1$) and

$2.'\: [F_0,a]\in Li_0^{1/2}$ (respectively, $Li^{1/2}$; $Li_0^{q/2}$)
for a dense set of $a\in \mathcal A$,

\noindent then we call $({\mathcal N},F_0)$ {\bf $\theta$-summable}
(respectively, {\bf weakly $\theta$-summable}; {\bf $\theta_q$-summable}).
By the Remark $\spadesuit\cdots\spadesuit$ above, $\theta$-summable 
implies weakly
$\theta$-summable implies $\theta_q$-summable for $0<q<1$.
We note that {\bf $\theta_1$-summable} = {\bf $\theta$-summable.}
\end{defn}
\noindent We warn the reader that what we call  weakly $\theta$-summable
for a Fredholm module is what Connes calls $\theta$-summable
in \cite{Co4} chapter IV. We do this to be consistent with Connes'
original definitions in the unbounded case discussed below
\cite{Co2}, \cite{Co3}.We now define the closely
related notion of unbounded Breuer-Fredholm modules and note that
in the above definition we do not bother with the extra adjective
{\it bounded}.

\begin{defn}
An {\bf odd unbounded Breuer-Fredholm module} for a unital Banach $*$-algebra
$\mathcal A$ is a pair $({\mathcal N},D_0)$ where $\mathcal A$ is
(continuously)
represented in $\mathcal N$ and $D_0$ is an unbounded self-adjoint operator
affiliated with $\mathcal N$ satisfying:

$1.\:(1+D_0^2)^{-1}\in {\mathcal {K_N}},\:and$

$2.\:[D_0,a]\in {\mathcal N}$ {\it for a dense set of} $a\in \mathcal A.$

\noindent If, in addition, our module satisfies:

$1.'\:(1+D_0^2)^{-1}\in Li_0$ (respectively, $Li$; $Li_0^q$ for $0<q\leq 1$),

\noindent then we call $({\mathcal N},D_0)$ {\bf $\theta$-summable}
(respectively, {\bf weakly $\theta$-summable}; {\bf $\theta_q$-summable}).
Again, $\theta$-summable $\Rightarrow$ weakly
$\theta$-summable $\Rightarrow$ $\theta_q$-summable for $0<q<1$.
We also note that {\bf $\theta_1$-summable} = {\bf $\theta$-summable}.
\end{defn}
\noindent By Corollary B.6 of Appendix B, we observe that for unbounded
Fredholm modules
our definition of $\theta$-summable agrees with Connes' original
definition, \cite{Co2}, \cite{Co3}, while our definition of weakly
$\theta$-summable coincides with his later definition of
$\theta$-summable, \cite{Co4}.

\subsection{\bf The Transformation $D\mapsto D(1+D^2)^{-1/2}$}

\begin{rems*}
In general, if $({\mathcal N},D_0)$ is an odd unbounded Breuer-Fredholm
module for some $\ca$ and
$F_0=D_0(1+D_0^2)^{-1/2}$, then
$({\mathcal N},F_0)$ is an
odd pre-Breuer-Fredholm module. Since
$1-F_0^2=(1+D_0^2)^{-1}$ the conditions
labelled 1. in the definitions coincide. The commutator conditions are
more subtle, but are handled by the strong-operator convergent integral:
$$F_0=\frac{1}{\pi}\int_0^\infty\lambda^{-1/2}D_0(1+D_0^2+\lambda)^{-1}
d\lambda,$$ \cite{BJ}, \cite{CP1}.

Using this integral formula, we showed in \cite{CP1} that an odd unbounded
$p$-summable Breuer-Fredholm module yields an odd $(p+\epsilon)$-summable
$pre$-Breuer-Fredholm module for each $\epsilon>0.$ By using completely
different techniques, F.A. Sukochev \cite{Suk}, and then \cite{CPS}
were able to eliminate the $\epsilon$ ($except$ in the case $p=1$
and $\mathcal N$ is type $II_\infty$). Unfortunately, these new techniques
do not allow us to handle certain smoothness  properties of the
transformation $D\mapsto D(1+D^2)^{-1/2}$ which are crucial in
obtaining integral formulae for spectral flow.

$\spadesuit$ In the present paper, we revisit the integral formula (specifically,
lemma 2.7 of \cite{CP1}) to show that if the unbounded module 
$({\mathcal N},D_0)$ is
$\theta$-summable (even weakly $\theta$-summable) then the {\bf bounded} module
$({\mathcal N},F_0)$ is $\theta_q$-summable for all $q$, $0<q<1.$
If we let $1/q=1+\epsilon$ then this is the same $\epsilon$-problem
we encountered in the finitely-summable situation (an artifact of
lemma 2.7 of \cite{CP1}). In a recent preprint \cite{Suk2} , F. A. Sukochev
has extended his techniques to the $\theta$-summable case to show
that if $({\mathcal N},D_0)$ is $\theta$-summable, then so is
$({\mathcal N},F_0)$.
However, we still need these integral techniques to handle the smoothness
of the transformation $D\mapsto D(1+D^2)^{-1/2}$, and so we are {\bf forced}
to consider $\theta_q$-summable modules for $0<q<1.$ $\spadesuit$

As in section 1 of \cite{CP1} we can obtain a genuine Breuer-Fredholm
module $({\mathcal N},\tilde{F}_0)$ from a pre-Breuer-Fredholm module,
$({\mathcal N},F_0)$ by letting $\tilde{F}_0=sign(F_0)$ where
$$sign(x)=\left\{\begin{array}{ll}
            +1 & \mbox{if $x\geq 0$} \\
            -1 & \mbox{if $x<0.$}
            \end{array}
          \right.$$
We observe that $({\mathcal N},\tilde{F}_0)$ has the same summability
flavour as $({\mathcal N},F_0)$ since
$$(\tilde{F}_0-F_0)=(1-F_0^2)(\tilde{F}_0+F_0)^{-1}.$$
We warn the reader that our sign function is never $0$ and so $sign(F_0)$
is a always a self-adjoint unitary: this differs from Connes' convention
\cite{Co4}.
\end{rems*}

\subsection{\bf Spectral Flow}

\begin{defn}
If $\{F_t\}$ is a continuous path of self-adjoint Breuer-Fredholm operators
in $\mathcal N$, then the definition of the {\bf spectral flow} of the
path, $sf(\{F_t\})$ is based on the following sequence of observations
in \cite{Ph1}:

\noindent 1. The map $t\mapsto sign(F_t)$ is usually discontinuous as is the
projection-valued mapping $t\mapsto P_t=\frac{1}{2}(sign(F_t)+1).$

\noindent 2. However, if $\pi:{\mathcal N}\to {\mathcal N}/{\mathcal K_{\cn}}$
is the canonical mapping, then $t\mapsto \pi(P_t)$ is continuous.

\noindent 3. If $P$ and $Q$ are projections in $\mathcal N$ and
$||\pi(P)-\pi(Q)||<1$
then $$PQ:rng(Q)\to rng(P)$$ is a Breuer-Fredholm operator and so
$ind(PQ)\in \real$ is well-defined.

\noindent 4. If we partition the parameter interval of $\{F_t\}$ so
that the $\pi(P_t)$ do not vary much in norm on each subinterval
of the partition then $$sf(\{F_t\}):=\sum_{i=1}^n ind(P_{t_{i-1}}P_{t_i})$$
is a well-defined and (path-) homotopy-invariant number which
agrees with the usual notion of spectral flow in the type $I_\infty$
case \cite{Ph},
and also agrees with all previous special definitions of type $II_\infty$
spectral flow, for example \cite{P1,P2}.

In particular, if the path, $\{F_t\}$ for $t\in [a,b]$ lies entirely in
$F_0+{\mathcal K_{\cn}}$, then $\pi(P_t)$ is constant by \cite{Ph1} and so
$$sf(\{F_t\})=ind(P_bP_a).$$
Thus, since the spectral flow of a path in $F_0+{\mathcal K_{\cn}}$
depends only on the endpoints, we will often denote it by $sf(F_a,F_b).$
All of the paths we consider in sections 3, 4, and 5 are of this kind.
\end{defn}

\subsection{\bf Spaces Of Breuer-Fredholm Operators}

\begin{rems*}
If $({\mathcal N},F_0)$ is a pre-Breuer-Fredholm module for a
Banach $*$-algebra ${\mathcal A},$ where $1-F_0^2$ and $|[F_0,a]|^2$
are in the invariant operator
ideal ${\mathcal I}$ (see Appendix A), then the operators
$F=uF_0u^*$ (for a dense set of
unitaries in ${\mathcal U}({\mathcal N})$) are self-adjoint Breuer-Fredholm
operators in ${\mathcal N}$ which also satisfy:

$1-F^2\in {\mathcal I},\:and$

$F-F_0\in {\mathcal I}^{1/2}.$

Moreover, any operator $F_t$ in the straight line path from
$F_0$ to $F=uF_0u^*$ also satisfies these conditions. Thus, we are led to
consider the affine space of ``allowable perturbations'' of $F_0$:
$${\mathcal I}_{F_0}:=F_0 + \{X\in {\mathcal I^{1/2}_{sa}}\:|1-(F_0+X)^2\in
{\mathcal I}\}
=F_0 + \{X\in {\mathcal I^{1/2}_{sa}}\:| F_0X+XF_0\in {\mathcal I}\}.$$
\end{rems*}
We observe that
if we let $P=\frac{1}{2}(sign(F_0)+1),$
then relative to the decomposition of $H$ determined by $P$,
$${\mathcal I}_{F_0} = \left(\begin{array}{cc}
                   {\mathcal I_{sa}} & {\mathcal I^{1/2}_{sa}} \\
                   {\mathcal I^{1/2}_{sa}} & {\mathcal I_{sa}}
                   \end{array} \right).$$

In appendix B
we show that for the operator ideals we are studying,
this space is a real Banach space in a natural norm and
$X\mapsto 1-(F_0+X)^2: {\mathcal I}_{F_0}\to {\mathcal I}$
is continuous. It is in these spaces that we study the spectral flow of
paths of self-adjoint Breuer-Fredholm operators.

In the case of unbounded modules, our space of ``allowable perturbations''
of $D_0$ will always be: $${\cm}_0=D_0+{\mathcal N}_{sa},$$
an affine space modelled on the real Banach space ${\mathcal N}_{sa}.$
Provided that we are careful in our choice of ${\mathcal I},$
(the $\epsilon$ problem) we show that the transformation
$D\mapsto F=D(1+D^2)^{-1/2}$ carries $D_0+{\mathcal N}_{sa}$
to $F_0+{\mathcal I}_{F_0}.$ and that this transformation is suitably
smooth. Thus we study the spectral flow of paths of ``unbounded
Breuer-Fredholm operators'', $\{D_t\}$ by considering the transformed
paths $\{F_t\}$ of genuine Breuer-Fredholm operators in the bounded setting.

\subsection{\bf One-forms}

\begin{rems*}
We will consider the affine spaces defined in the previous paragraphs
as real Banach manifolds, $M.$ For example, if $M=F_0+{\mathcal I}_{F_0},$
then for any $F$ in $M$ the tangent space at $F$ is
$T_F(M)={\mathcal I}_{F_0}.$

If $f$ is any continuous real-valued function on $\real$ then provided
$f:{\mathcal I}\to L^1$ is continuous, we can define a
one-form, $\alpha$ via:
$$\alpha(X)=\frac{1}{C}\tau(Xf(1-F^2))$$
where $$F\in M=F_0+{\mathcal I}_{F_0},$$
$$X\in T_F(M)={\mathcal I}_{F_0},\:and$$
$$C=\int_{-1}^{1}f(1-t^2)dt.$$

The integral of this one-form along a path $\{F_t\}$ for $t\in [0,1]$ in $M$
is given by $$\frac{1}{C}\int_{0}^{1}\tau(\frac{d}{dt}(F_t)f(1-F_t^2))dt.$$
To see that these integrals are independent of path, we show that such
one-forms (for suitable $f$) are closed: that is, their exterior derivatives
vanish identically. A Poincar\'{e} Lemma completes the proof of
path independence. 

We use the invariant definition of {\bf exterior differentiation} \cite{Sp}.
For $F$ in $M,$ we have $X,Y$ in $T_F(M)={\mathcal I}_{F_0}$ realized
as tangent vectors at $F$ by differentiating the curves $F+sX$ and $F+sY$
at $s=0.$ That is, we also consider $X$ and $Y$ as the canonical vector
fields on $M$ (or flows on $M$) given by flowing in the $X$ direction
or $Y$ direction. Then, by definition:
$$d\alpha(X,Y)=X\!\cdot\!(\alpha(Y))-Y\!\cdot\!(\alpha(X))-\alpha([X,Y]).$$
Since $X$ and $Y$ commute {\it as flows} the last term is $0$ and so
drops from the calculation. For the straight line paths above, this means 
showing that:

$$0=d\alpha(X,Y)=\frac{1}{C}\left[\tau(Xf[1-(F_0+Y)^2])-\tau(Yf[1-(F_0+X)^2])
\right].$$

$\spadesuit$ Thus, to prove that our integral formulae in the bounded 
setting yield
spectral flow (of the path) we are reduced to showing that: 1) our integrals
are independent of path (by the procedure indicated above), and 2) that for
certain special paths where we can actually calculate the integrals
we get the desired answer. It is this second calculation that we do in the
following section. $\spadesuit$
\end{rems*}

\section{\bf RELATIVE INDEX OF TWO PROJECTIONS}

The {\it essential codimension} (or {\it relative index}) of two projections
is a fundamental
tool. Formulae for this index have a long history: see for example,
\cite{ASS}, \cite{BW}, \cite{P2}, \cite{Ph1} and the references
contained therein. The next result
subsumes all previous ones to our knowledge.

\begin{thm}
Let $f:[-1,1]\rightarrow \bf{R}$ be a continuous odd function
with $f(1)\neq 0$. Let $P$ and $Q$
be projections with $P-Q\in\mathcal K_{\cn}$
and $f(P-Q)$ trace class.
Then ind$(QP)=\frac{1}{f(1)}\tau[f(P-Q)]$
where ind$(QP)$ is the index of $QP$ as an operator from
$PH$ to $QH$.
\end{thm}

The proof depends on a preliminary result.

\begin{prop}
Let $P$ and $Q$ be projections on $H$,
then the subspaces  $ranP\cap ker Q$ and
$ker P\cap ran Q$ are mutually orthogonal,
closed and invariant under both $P$ and $Q$.
Let $H_1$ be the orthogonal complement of their direct sum.
Then $H_1$ is invariant under both $P$ and $Q$ so that
$P_1=P|_{H_1}$ and $Q_1= Q_{|_{H_1}}$ are projections in
$B(H_1)$ with $P_1-Q_1=(P-Q)_{|_{H_1}}$, $Q_1-P_1=(Q-P)_{|_{H_1}}$.
Then, there exists a self-adjoint unitary $U$ in $\{P,Q,1\}^{\prime\prime}$
which is 1 on $H_1^\perp$ and is such that $U(P_1-Q_1)U^*=Q_1-P_1$.
\end{prop}

\begin{proof}[\bf Proof]
We work on $H_1$ and observe that on this space, $ranP_1\cap kerQ_1$
and $kerP_1\cap ranQ_1$ are both $\{0\}$. Let $B=1-(P_1+Q_1)$ and let
$B=U|B|$ be the polar decomposition of $B.$ Then $B$ anticommutes with
$(P_1-Q_1)$ and so $B^2$ commutes with $(P_1-Q_1)$, and hence any continuous
function of $B^2$ commutes with $(P_1-Q_1).$ In particular,
$|B|$ commutes with $(P_1-Q_1).$
One easily calculates that: $$U(P_1-Q_1)|B| = (Q_1-P_1)U|B|.$$
That is, $U(P_1-Q_1)$ agrees with $(Q_1-P_1)U$ on $ran|B|$. So, since
$B = B^* = |B|U,$ it suffices
to see that the self-adjoint operator $B$ has dense range (on $H_1$).
This is equivalent to $kerB=\{0\}$ (on $H_1$),
which is easily seen to be equivalent to the conditions:
$$kerP_1\cap ranQ_1 = \{0\} = ranP_1\cap kerQ_1.$$
Finally, we extend $U$ to be $1$ on $H_1^\perp.$
\end{proof}

\begin{rems*} If $P,Q$ are as in the proposition with
$QP$ regarded as mapping $PH$ to $QH$ then its kernel
is $ranP\cap ker Q$ and its cokernel is the kernel of
$PQ$ on $QH$ or, $ranQ\cap ker P$. In particular
if $P-Q$ is in ${\mathcal K}_{\mathcal N}$ then $QP:PH\rightarrow QH$
is a Breuer-Fredholm operator and
$${ind}(QP)={dim(ran}P\cap ker Q)-{dim(ran}Q\cap ker P).$$
$$=\tau([ranP\cap kerQ])-\tau([kerP\cap ranQ]).$$
\end{rems*}

We are now able to complete the proof of the theorem.

\begin{proof}[\bf Proof]
Using the notation of the previous proposition to define $H_1$ write
$$H=H_0\oplus  H_1=(ranP\cap kerQ)\oplus(kerP\cap ranQ)\oplus  H_1.$$
 Then
$$P=\left(\begin{array}{cc}
1 & 0 \cr 0 & 0
\end{array}\right)\oplus P_1,$$
and
$$Q=\left(\begin{array}{cc}
0 & 0\cr 0 & 1
\end{array}\right)\oplus Q_1.$$
So
$$f(P-Q)=\left(\begin{array}{cc}
f(1) & 0 \cr
0 & -f(1)
\end{array}\right)\oplus f(P_1-Q_1).$$
Now by the previous proposition and the fact that $f$ is odd:
$$Uf(P_1-Q_1)U^* = -f(P_1-Q_1).$$
By assumption this operator is trace class so we get
$\tau[f(P_1-Q_1)]=0$.
That is,
$$\tau[f(P-Q)]=f(1)\{\tau([ranP\cap ker Q])
-\tau([kerP\cap ranQ])\} =f(1)ind(QP)$$
by the remark.
\end{proof}

\begin{rems*}
At several places in this paper we consider the smooth functions,
$f:\real^{+} \to \real^{+}$ of the form $f(x)=x^{-r}e^{-x^{-\sigma}}$
where $r\geq 0$ and $\sigma \geq 1.$ These functions are (of course)
defined to be $0$ at $x=0.$
\end{rems*}

\begin{cor}
If $P$ and $Q$ are projections and $P-Q$ is $n$-summable for $n\geq 1$
 (not necessarily an
integer) then
$$ind(QP)=\tau[(P-Q)|P-Q|^{n-1}].$$
\end{cor}

The case $n=2k+1$ an odd integer is in \cite{ASS} and \cite{Ph1}.

\begin{cor}
If $P$ and $Q$ are projections and $P-Q$ is $\theta$-summable
(i.e., $(P-Q)^2\in Li_0$) then,
$$ind(QP)=e \tau[(P-Q)e^{-(P-Q)^{-2}}].$$
If $P-Q$ is weakly $\theta$-summable
(i.e., $(P-Q)^2\in Li$) then for $\epsilon>0$,
$$ind(QP)=e \tau[(P-Q)e^{-|P-Q|^{-(2+\epsilon)}}].$$
\end{cor}

\begin{proof}[\bf Proof]
The first statement follows from the theorem and Corollary B.5.
The second statement follows from Lemma B.4 since
the function $x\mapsto e^{-tx^{-1}}$ dominates $x\mapsto e^{-x^{-(1+\delta)}}$
as $x\to \infty$ for any fixed $t>0$ and $\delta>0.$
\end{proof}

\begin{cor}
If $P$ and $Q$ are projections and $P-Q$ is $\theta_q$-summable
(i.e., $0<q\leq 1$ and $(P-Q)^2\in Li^q_0$). Then,
for $1+\epsilon = 1/q$ we have
$$ind(QP)=e \tau[(P-Q)e^{-|P-Q|^{-2(1+\epsilon)}}]$$
\end{cor}

\begin{proof}[\bf Proof]
This follows from Corollary B.5.
\end{proof}

\section{\bf THE SPECTRAL FLOW FORMULA, SIMPLEST CASE}

Theorem 3.1 enables us to
establish a result which serves as a prototype
for all of the formulae in subsequent sections.

\begin{thm}
Let $P$ and $Q$ be infinite and co-infinite projections
in the semi-finite factor $\mathcal N$ and suppose $(P-Q)^2\in Li_0^q$
for some $0<q\leq 1$. Then $F_0=2P-1$ and $F_1=2Q-1$
are self-adjoint Breuer-Fredholm operators as is the
path $F_t= F_0+t(F_1-F_0)$.
Let $r\geq 0$
then
$$sf(\{F_t\})=\frac{1}{C_{r,q}}
\int_0^1 \tau\left(\frac{d}{dt}(F_{t})(1-F_t^2)^{-r} e^{-(1-F_t^2)^{-1/q}}
\right) dt$$
where
$$C_{r,q} =\int_{-1}^1 (1-u^2)^{-r} e^{-(1-u^2)^{-1/q}}du.$$
\end{thm}

\begin{proof}[\bf Proof]
We have
\begin{eqnarray} \frac{d}{dt}(F_{t})
&=&F_1-F_0\;\;=\;\;2(Q-P)\;\; {\text and}\nonumber\\
1-F_{t}^2&=&t(1-t)(F_1-F_0)^2\;\;=\;\;4t(1-t)(Q-P)^2\nonumber
\end{eqnarray}
and so by assumption
$|1-F_t^2|^{1/q}\in Li_0$. This means  $e^{-(1-F_t^2)^{-1/q}}$
is trace class for $t\in [0,1)].$
Thus,
\begin{eqnarray}
& &\int_0^1\tau\left(\frac{d}{dt}(F_t)(1-F_t^2)^{-r} e^{-(1-F_t^2)^{-1/q}}
\right) dt\nonumber\\
&=&\int_0^1\tau\left[2(Q-P)[4(t-t^2)(Q-P)^2]^{-r}e^{-[4(t-t^2)(Q-P)^2]^{-1/q}}
\right]dt.\nonumber
\end{eqnarray}
\noindent Now for each $t\in(0,1)$ define:
$$f_t(x)=2x[4(t-t^2)x^2]^{-r}e^{-[4(t-t^2)x^2]^{-1/q}}$$
and apply Theorem 3.1 to get
\begin{eqnarray}
   &   & \int_0^1\tau\left[\frac{d}{dt}(F_t)(1-F_t^2)^{-r} e^{-(1-F_t^2)^{-1/q}}
\right] dt \nonumber\\
   & = & \int_0^1 f_t(1)ind(QP)dt \nonumber\\
   & = & ind(QP)\int_0^1 2[4(t-t^2)]^{-r}e^{-[4(t-t^2)]^{-1/q}}dt \nonumber\\
   & = & C_{r,q}\,ind(QP) \nonumber\\
   & = & C_{r,q}\, sf(\{F_t\}) \nonumber
\end{eqnarray}
\noindent where the penultimate equality is obtained by
using the change of variable $u=2t-1$ and the last
by the definition of spectral flow \cite{Ph1}.
\end{proof}

We pause at this point to draw some
conclusions from the previous analysis which we believe to be of
independent interest. These results are $\theta$-summable
versions of Theorem 3.3 of \cite{Ph1}. These results
treat the case of (bounded) Breuer-Fredholm modules and do not need the full
machinery of one-forms nor of the appendices. The latter
are however necessary for the more general case of pre-Breuer-Fredholm modules
which arise naturally when we reduce the unbounded setting to the
bounded setting.

In order to emphasise the elegance of the following results we
consider only the case $r=0$ from Theorem 4.1.
The case of general $r$ is covered in the next section.

\begin{thm}
Let $\mathcal A$ be a unital Banach $*$-algebra and let
$({\mathcal N}, F_0)$ be an odd $\theta_q$-summable Breuer-Fredholm module
for $\mathcal A$ for some $q, 0<q\leq 1$.
Let $P=\frac{1}{2}(1+F_0).$ For each unitary $u\in \mathcal A$
with $[F_0,u]\in Li_0^{q/2}$, the path $F_t^u=F_0+t(uF_0u^*-F_0)$
lies in the self-adjoint Breuer-Fredholms and
$$ind(PuP)=sf(\{F_t^u\})=\frac{1}{C_{0,q}}\int_0^1
\tau \left(\frac{d}{dt}(F^u_t)e^{-|1-(F^u_t)^2|^{-1/q}}\right) dt$$
where $C_{0,q}= \int_{-1}^1 e^{-(1-u^2)^{-1/q}}du$.
\end{thm}

\begin{proof}[\bf Proof]
The first equality follows from the discussion at the
beginning of section 3 of \cite{Ph1} and the second
equality from the previous theorem with $Q=\frac{1}{2}(uF_0u^*+1)$
since
$$P-Q = \frac{1}{2}(F_0-uF_0u^*)= \frac{1}{2}[F_0,u]u^*\in Li_0^{q/2}.$$
\end{proof}



\begin{cor}
Let $\mathcal A$ be a unital Banach $*$-algebra and let
$({\mathcal N}, F_0)$ be an odd weakly $\theta$-summable Breuer-Fredholm module
for $\mathcal A$.
Let $P=\frac{1}{2}(1+F_0)$. For each unitary $u\in \mathcal A$
with $[F_0,u]\in Li^{1/2}$, the path $F_t^u=F_0+t(uF_0u^*-F_0)$
lies in the self-adjoint Breuer-Fredholms and
$$ind(PuP)=sf(\{F_t^u\})=\frac{1}{C^{\epsilon}}\int_0^1
\tau \left(\frac{d}{dt}(F^u_t)e^{-|1-(F^u_t)^2|^{-1-\epsilon}}\right) dt
$$
where $C^{\epsilon}= \int_{-1}^1 e^{-(1-u^2)^{-1-\epsilon}}du$.\\
If either the module is $\theta$-summable or if $||\;|\;[F_0,u]\;|^2||_{Li}
<\frac{2}{3}$, then we can set $\epsilon = 0$ in the formula.
\end{cor}

\begin{proof}[\bf Proof]
We let $q=1/(1+\epsilon)$, then
$$\mu_x(1-F_0^2)\leq \frac{K}{\log x}=\frac{K}{(\log x)^{1-q}}
(\frac{1}{\log x})^q
=o((\frac{1}{\log x})^q)$$ so that $({\mathcal N}, F_0)$
is  $\theta_q$-summable for  ${\mathcal A}$.
\end{proof}



\section{\bf SPECTRAL FLOW FORMULAE, BOUNDED CASE}

We suppose we have a $\theta_q$-summable
pre-Breuer-Fredholm module $({\mathcal N}, F_0)$
for the unital Banach $*$-algebra, $\mathcal A$, where $0<q\leq 1$.

We recall:
$$(Li_0^q)_{F_0}= \{X\in (Li_0^{q/2})_{sa}\: |\:
F_0X+XF_0\in Li_0^q\}.$$
Now we set ${\mathcal M}_q=F_0+(Li_0^q)_{F_0}.$

See appendices A and B for more details on these spaces.
There we will also show that
if $F_1\in{\mathcal M}_q$ then  $F_0X+XF_0\in Li_0^q$
if an only if  $F_1X+XF_1\in Li_0^q$ so that the definition of
${\mathcal M}_q$ is independent of base point.
Moreover ${\mathcal M}_q$ is contained in the self-adjoint
Breuer-Fredholms so that if $F_u$ with $u\in [0,1]$ is a
norm continuous path in this
space then $sf\{F_u\}$ is well-defined as in Section 2.

\begin{prop}
Let $({\cn},F_0)$ be an odd $\theta_q$-summable pre-Breuer-Freholm module
for the Banach $*$-algebra $\ca$
We define a one-form $\alpha_r$
on ${\mathcal M}_q$ via:
$$(\alpha_r)_F(X)=\frac{1}{C_{r,q}}\tau\left(X|1-F^2|^{-r}
e^{-|1-F^2|^{-1/q}}\right),$$
for
$F\in{\mathcal M}_q$, $X\in T_F({\mathcal M}_q)=(Li_0^q)_{F_0}.$
Then, the one-form $\alpha_r$ is closed (recall that $C_{r,q}$ is
defined in Theorem 2 of Section 4).
\end{prop}

\begin{proof}[\bf Proof]
By Theorem C.5 of Appendix C,
the following derivative exists
$$\frac{d}{ds}|_{s=0}\left\{\tau(Y|1-(F+sX)^2|^{-r}e^{-|1-(F+sX)^2|^{-1/q}})
\right\}$$
and equals
$$\frac{i}{2\pi}\int_\sigma\tau\left\{Y\left[g_r(T) \;,\; R_{\lambda}(T) \left[
[F,X]_{+} \;,\; T \right]_+R_{\lambda}(T)\right]_+\right\} m(\lambda)d\lambda$$
where $T=1-F^2$ ; where $[\cdot,\cdot]_+$ denotes the anticommutator;
where $g_r(T)=|T|^{-r/2}e^{-1/2|T|^{-1/q}}$; where
$R_\lambda(T)=(\lambda T^2-1)^{-1}$; and
$m(\lambda)=\lambda^{r/4}e^{-(\lambda^{1/2q})/2}.$

Within the integral there are four terms, one
of which is:
\begin{eqnarray}
& &\tau\{YR_\lambda (T)[F,X]_+TR_\lambda (T)g_r(T)\}\nonumber\\
&=&\tau\{XTR_\lambda (T) g_r(T)YR_\lambda (T)F\}
+\tau\{XFTR_\lambda (T) g_r(T)YR_\lambda (T) \}\nonumber\\
&=&\tau\{Xg_r(T)R_\lambda (T) TYFR_\lambda (T) \}
+\tau\{Xg_r(T)R_\lambda (T) TFYR_\lambda (T) \}\nonumber\\
&=&\tau\{Xg_r(T)R_\lambda (T) T[F,Y]_+R_\lambda (T) \},
\nonumber
\end{eqnarray}
which is precisely one of the other terms in the integral
with the roles of $X$ and $Y$ interchanged.

The other two terms are handled in a similar fashion.
Thus
\begin{eqnarray}
                &   & \frac{d}{ds}|_{s=0}[\tau(Y|
 1-(F+sX)^2|^{-r}e^{-|1-(F+sX)^2|^{-1/q}})] \nonumber \\
                & = & \frac{d}{ds}|_{s=0}[\tau(X|
1-(F+sY)^2|^{-r}e^{-|1-(F+sY)^2|^{-1/q}})] \nonumber
\end{eqnarray}
proving the result.
\end{proof}



\begin{cor}
If $({\cn},F_0)$ is a weakly $\theta$-summable pre-Breuer-Fredholm module for
$\ca$ and we let ${\mathcal M}=F_0+(Li)_{F_0},$ then:
for each $\epsilon>0$ we define a one-form on ${\mathcal M}$
by
$$(\alpha_{r,\epsilon})_F(X)=\frac{1}{C_{r,1/(1+\epsilon)}}
\tau\left(X|1-F^2|^{-r}e^{-|1-F^2|^{-1-\epsilon}}\right),$$
where
$F\in{\mathcal M}, X\in T_F({\mathcal M})=(Li)_{F_0}$.
Then $\alpha_{r,\epsilon}$ is closed. Here
$$C_{r,1/(1+\epsilon)}=\int_{-1}^1 (1-u^2)^{-r}e^{-(1-u^2)^{-1-\epsilon}}du.$$
If the module is $\theta$-summable then the conclusion holds on
${\mathcal M_1}=F_0+(Li_0)_{F_0}$  for the one-form $\alpha_{r,0}$.
\end{cor}

\begin{proof}[\bf Proof]
 With $q=1/(1+\epsilon)$ this follows from the fact that
the inclusion of $Li$ into $Li_0^q$ is bounded by Lemma A.5.
\end{proof}

\begin{defn} Fix $({\mathcal N}, F_0)$ with
$1-F_0^2\in {Li}_0^q$.
Let $F\in {\mathcal M}_q$ and let $F_t=F_0+t(F-F_0)$ for $t\in[0,1]$
be the linear path in ${\mathcal M}_q$ from $F_0$ to $F$. Define
$\theta_r:{\mathcal M}_q\to \real$ for $r\geq 0$ via:
$$\theta_r(F)=\frac{1}{C_{r,q}}
\int_0^1\tau\left[(F-F_0)|1-F_t^2|^{-r}e^{-|1-F_t^2|^{-1/q}}\right]dt.$$
\end{defn}

Recall that by definition
$$d\theta_F(X)=\frac{d}{ds}|_{s=0}(\theta(F+sX))$$
for $X\in T_F({\mathcal M}_q)$.

\begin{prop}(Poincar\'{e} Lemma)
With the above definitions, $d\theta_r=\alpha_r$.
\end{prop}

\begin{proof}[\bf Proof]
Fix $F\in {\mathcal M}_q$ and let $Y=F-F_0.$ For $T\in {Li}_0^q,$ let
$f_r(T)=|T|^{-r}e^{-|T|^{-1/q}}$ so that
$$\theta_r(F)=\frac{1}{C_{r,q}}\int_0^1\tau\left[Yf_r(1-F_t^2)\right]dt$$
where $F_t=F_0+tY$ for $t\in [0,1]$. Let $X\in (Li_0^q)_{F_0}$
so that
$$\theta_r(F+sX)=\frac{1}{C_{r,q}}\int_0^1\tau\left[(Y+sX)f_r(1-({F^s_t})^2)
\right]dt$$
where $F_t^s =F_0 + t(F+sX-F_0)=F_0+t(Y+sX)=F_t+s(tX)$
for real $s$. Now by the product rule
$$\frac{d}{ds}\tau[(Y+sX)f_r(1-({F^s_t})^2)]=\tau[(X)f_r(1-({F^s_t})^2)]+
\tau[(Y+sX)\frac{d}{ds}f_r(1-({F^s_t})^2)].$$
By Theorem C.5 of Appendix C with $F_t$ in place of $F_0$, $tX$ in place of
$X$, and an additive change of the variable $s$ we get:
$$\frac{d}{ds}(f_r(1-({F^s_t})^2))$$
$$=\frac{i}{2\pi}\left\{\int_\sigma\left
[f_r(T_{t,s}), (\lambda T_{t,s}^2-1)^{-1}\left[
(F_t^stX+tXF_t^s), T_{t,s}\right]_+
(\lambda T_{t,s}^2-1)^{-1}\right]_+
\lambda^{r/4}e^{-\lambda^{1/2q}/2}d\lambda\right\}$$
where $T_{t,s} = 1-(F_t^s)^2$
and $[\cdot,\cdot]_+$ denotes the anti-commutator.
Similar estimates to those of the proof
of Theorem C.5 show that
$$(t,s)\mapsto \frac{d}{ds}\tau\left[(Y+sX)f_r(1-({F^s_t})^2)\right]$$
is continuous on $[0,1]\times \real$. Since
$$(t,s)\mapsto \tau\left[(Y+sX)f_r(1-({F^s_t})^2)\right]$$
is also continuous, Theorem 11 of Chapter X of \cite{L},
allows us to differentiate under the integral sign to compute
$$\frac{d}{ds}|_{s=0}\int_0^1\tau\left[(Y+sX)f_r(1-({F^s_t})^2)\right]dt$$
$$
=\int_0^1\left\{\tau[Xf_r(1-{F_t}^2)]+\tau[Y\frac{d}{ds}|_{s=0}
f_r(1-({F^s_t})^2)]\right\}dt
\eqno (**).$$
Now
$$\tau\left[Y\frac{d}{ds}|_{s=0}f_r(1-({F^s_t})^2)\right]=
\tau\left[tX\frac{d}{ds}|_{s=0}f_r(1-(F_t+sY)^2)\right]
$$
$$=t\tau\left[X\frac{d}{ds}|_{s=0}
f_r(1-(F_{t+s})^2)\right]=t\tau\left[X\frac{d}{dt}f_r(1-{F_t}^2)\right].$$
Substitution in (**) and integration by parts gives, for the RHS of (**)
$$\int_0^1\frac{d}{dt}\left[t\tau(Xf_r(1-{F_t}^2))\right]dt=\tau(Xf_r(1-F^2))$$
as $F_1=F$. Dividing by the normalisation constant gives the result.
\end{proof}

\begin{cor}
The integral of the one-form $\alpha_r$ along a piecewise
$C^1$ path $\Gamma$ in ${\mathcal M}_q$
depends only on the endpoints of the path $\Gamma$.
\end{cor}

\begin{proof}[\bf Proof]
This follows as in Proposition 1.5 of \cite{CP1}.
\end{proof}

\begin{defn} Fix $F_0$ a self-adjoint Breuer-Fredholm with
$1-F_0^2\in Li_0^q$.
Let $F\in {\mathcal M}_q$ and
recall from Section 2 that $\tilde F=sign(F)\in {\mathcal M}_q$.
 Let $\{F_t\}$ for $t\in[0,1]$ be a $C^1$ path in ${\mathcal M}_q$
beginning at $F$ and ending at $\tilde F$.
For example
$F_t=F+t(\tilde F-F)$ will do. Define
$$\gamma_{r,q}(F)=\frac{1}{C_{r,q}}
\int_0^1\tau\left[\frac{d}{dt}(F_t) f_r(1-F_t^2)\right]dt.$$
\end{defn}

It follows by the previous corollary that $\gamma_{r,q}$ is well-defined.
Moreover, it is clear by considering the linear
path that if $F_1$ and $F_2$ are unitarily equivalent in
${\mathcal M}_q$
then $\gamma_{r,q}(F_1)=\gamma_{r,q}(F_2)$.

\begin{thm}
Let $\mathcal A$ be a unital Banach $*$-algebra and let
$({\mathcal N}, F_0)$ be an odd $\theta_q$-summable
pre-Breuer-Fredholm module
for $\mathcal A$ for some $q$, $0<q\leq 1$.
Let $r\geq 0$ then
if $F_j\in {\mathcal M}_q$ for $j=1,2$, the spectral flow along
any piecewise $C^1$ path $\{F_t\}$ in ${\mathcal M}_q$, $t\in [1,2]$
joining $F_1$ and $F_2$  is given by
$$sf(F_1,F_2)=\frac{1}{C_{r,q}}\int_1^2\tau\left[\frac{d}{dt}(F_t)f_r(1-F_t^2)
\right]dt
+\gamma_{r,q}(F_2)-\gamma_{r,q}(F_1)$$
where $f_r(T)=|T|^{-r}e^{-|T|^{-1/q}}$.
\end{thm}

\begin{proof}[\bf Proof]
The formula on the right is just the integral
of $\alpha_r$ along the piecewise $C^1$ path in ${\mathcal M}_q$
from $\tilde F_1$ to  $\tilde F_2$ made up of the
three parts, $\tilde F_1$ to $F_1$, then  $\{F_t\}$ for $t\in [1,2]$ and
finally $F_2$ to $\tilde F_2$. But since the integral of $\alpha_r$ is
independent
of the path in ${\mathcal M}_q,$ we may equally use the straight line path
joining $\tilde F_1$ and $\tilde F_2$.
Then we have  by Theorem 4.1:
$$ sf(\tilde F_1, \tilde F_2)=\frac{1}{C_{r,q}}
\int_1^2\tau\left[\frac{d}{dt}(F_t)f_r(1-F_t^2)^{-1}\right]dt
+\gamma_{r,q}(F_2)-\gamma_{r,q}(F_1).$$
Finally,
$$sf(\tilde F_1, \tilde F_2) =sf(\tilde F_1, F_1) +sf( F_1, F_2)+sf( F_2,
\tilde F_2)$$
$$
= sf(F_1, F_2)$$
as there is no spectral flow along the paths joining $F_j$ and $\tilde F_j$
as noted in the proof of Theorem 1.7, p.683 of \cite{CP1}.
\end{proof}


\begin{cor}
If we assume that
$({\mathcal N}, F_0)$ is an odd weakly $\theta$-summable
pre-Breuer-Fredholm module
for $\mathcal A$.
Then for
$$F_j\in {\mathcal M}=F_0+(Li)_{F_0},\quad j=1,2,$$
and the remaining hypotheses intact we get:
$$sf(F_1,F_2)=\frac{1}{C^\epsilon}\int_1^2\tau\left[\frac{d}{dt}(F_t)
e^{-|1-F_t^2|^{-1-\epsilon}}\right]dt
+\gamma^\epsilon(F_2)-\gamma^\epsilon(F_1)$$
(recall $C^\epsilon=\int_{-1}^1 e^{-(1-u^2)^{-1-\epsilon}}du$)
and $$\gamma^\epsilon(F)=\frac{1}{C^\epsilon}\int_0^1
\tau\left[(\tilde F- F)e^{-|1-F_t^2|^{-1-\epsilon}}\right]dt$$
where this last integral is along the linear path from $F$ to $\tilde F=
sign(F)$.\\
Again, if the module is $\theta$-summable we can take $\epsilon = 0$ on
${\mathcal M_1}=F_0+(Li_0)_{F_0}.$
\end{cor}

\begin{proof}[\bf Proof]
This follows from the fact that the inclusion
of $Li$ in $Li^q_0$ for $q=1/(1+\epsilon)$ is bounded by Lemma A.5.
\end{proof}

\begin{thm}
Let
$({\mathcal N}, F_0)$ be an odd $\theta_q$-summable
pre-Breuer-Fredholm module
for the unital Banach $*$-algebra $\mathcal A$ and for some $q$, $0<q\leq 1$.
Let $P=\frac{1}{2}(sign(F_0)+1)$.
For each unitary $u$ in $\mathcal A$ with $[F_0,u]\in Li_0^{q/2}$
the path $F_t^u=F_0+t(uF_0u^*-F_0)$
lies in ${\mathcal M}_q=F_0+(Li_0^q)_{F_0}$
and
$$ind(PuP)=sf(\{F_t^u\})
=\frac{1}{C_{r,q}}\int_0^1\tau\left[\frac{d}{dt}(F_t^u)f_r(1-(F_t^u)^2)
\right]dt,$$
where $f_r(T)=|T|^{-r}e^{-|T|^{-1/q}}$ for $r\geq0.$
\end{thm}

\begin{proof}[\bf Proof]
 The second equality follows from the previous theorem
since
$$\gamma_{r,q}(uF_0u^*)=\gamma_{r,q}(F_0).$$
The first equality follows from the discussion at the beginning of
section 3 of \cite{Ph1}.
\end{proof}

An analogue of Corollary 5.8 follows from Theorem 5.9 in the
obvious fashion.

\section{\bf PRELIMINARIES FOR THE UNBOUNDED CASE}

\begin{lemma} [{\it cf} \cite{CP1} Lemma 6, Appendix B]
If $D_0$ is an unbounded self-adjoint operator, $A$
a bounded self-adjoint operator and $D=D_0+A$
then

\noindent(1) $(1+D^2)^{-1}\leq f(||A||)(1+D_0^2)^{-1}$ and

\noindent(2) $-(f(||A||)-1)(1+D_0^2)^{-1} \leq (1+D^2)^{-1}-(1+D_0^2)^{-1}
\leq (f(||A||)-1)(1+D_0^2)^{-1}$ where

\noindent$f(a)=1+\frac{1}{2}(a^2 +a\sqrt{a^2+4}).$
\end{lemma}

\begin{proof}[\bf Proof]
The first result is the one cited. For (2)
the right hand inequality follows from (1) by subtracting $(1+D_0^2)^{-1}$.
The left hand inequality of (2) follows by noting first that
 $(1+D_0^2)^{-1}\leq f(||A||)(1+D^2)^{-1}$ by interchanging the roles of
 $D$ and $D_0$. Then
$$(\frac{1}{f(||A||)}-1)(1+D_0^2)^{-1} \leq (1+D^2)^{-1}-(1+D_0^2)^{-1}$$
or
$$-\frac{(f(||A||)-1)}{f(||A||)}(1+D_0^2)^{-1}
 \leq (1+D^2)^{-1}-(1+D_0^2)^{-1}.$$

\noindent Since $-(f(||A||)-1)\leq-\frac{(f(||A||)-1)}{f(||A||)}$ we are done.
\end{proof}

\begin{lemma}[\cite{CP1} Lemma 2.7]  If $D=D_0+A$ with $A\in {\mathcal N}_{sa}$
and with $D_0$ affiliated to $\mathcal N$ we define
$F_A=D(1+D^2)^{-1/2}$ and $F_0=D_0(1+D_0^2)^{-1/2}$
then  for $0<\sigma<1$,
$(F_A-F_0)^2=B_\sigma(1+D_0^2)^{\sigma/2}$
where $B_\sigma\in \mathcal N$ and
$||B_\sigma||\leq C(\sigma)||A||$
where $C(\sigma)$ is a constant depending only on $D_0$ and $\sigma$.
\end{lemma}

\begin{cor}
Let $({\mathcal N}, D_0)$ be weakly $\theta$-summable and let
$D_t=D_0+A_t\in D_0+{\mathcal N}_{sa}$ be an operator-norm
continuous path. Then for fixed $q$, $0<q<1$
$$t\mapsto F_t=D_t(1+D_t^2)^{-1/2}\in F_{D_0}+(Li_0^{q})_{F_{D_0}}$$
is continuous, where $F_{D_0}=D_0(1+D_0^2)^{-1/2}.$
\end{cor}

\begin{proof}[\bf Proof]
Let $q<\sigma<1$. Since $(1+D_0^2)^{-1}\in Li$ by definition
we have $X_t=F_t-F_{D_0}\in Li^{\sigma/2}\subseteq Li_0^{q/2}$
by Lemma 6.2. Using Lemma 6.2 again with
$q$ in place of $\sigma$ we get
$$||X_t-X_{t_0}||_{Li_0^{q/2}}
\leq C(q)||A_t-A_{t_0}||.||(1+D_{t_0}^2)^{-q/2}||_{Li_0^{q/2}}
\rightarrow 0$$
as $t\rightarrow t_0$ (noting that $(1+D_{t_0}^2)^{-q/2}
\in Li_0^{q/2}$ by part (1) of Lemma 6.1).
 Thus
$ F_t \in  F_{D_0}+Li_0^{q/2}$ and is continuous in that space.
Now $t\rightarrow 1- ( F_{D_0}+X_t)^2=1- F_t^2=(1+D_{t}^2)^{-1}$
is in $Li\subseteq Li_0^{q/2}$ by part (1) of Lemma 6.1.
Moreover it is continuous in
$Li_0^{q}$
by part (2) and the fact that $Li_0^{q}$ is an invariant operator ideal
(see Appendix A). The result now follows from
part (3), Lemma B.12.
\end{proof}

\begin{prop}
Let $D_0$ be an unbounded self-adjoint operator
affiliated with $\mathcal N$.
If $\mathcal I$ is an invariant operator ideal in $\mathcal N$,
 $(1+D_0^2)^{-1}\in\mathcal I$ and $t\mapsto D_t=D_0+A_t
\in D_0+{\mathcal N}_{sa}$ is a $C^1$ path in the operator
norm and  $F_t=D_t(1+D_t^2)^{-1/2}$ then
$t\mapsto 1-F^2_t$ is  $C^1$ in  $\mathcal I$.
\end{prop}

\begin{proof}[\bf Proof]
 By part (1) of Lemma 6.1 we have
 $1-F^2_t=(1+D_t^2)^{-1}\in{\mathcal I}$
and depends continously on $t$
by part (2). Using Lemma 2.9 of [CP1] (in the notation
used there we set $x=1$) we have
$$\frac{1}{t-t_0}[(1+D_t^2)^{-1}-(1+D_{t_0}^2)^{-1}]$$
$$
=-D_{t_0}(1+D_{t_0}^2)^{-1}\frac{A_t-A_{t_0}}{t-t_0}(1+D_t^2)^{-1}
-(1+D_{t_0}^2)^{-1}\frac{A_t-A_{t_0}}{t-t_0}D_t(1+D_t^2)^{-1}.$$
Now, by assumption $\frac{A_t-A_{t_0}}{t-t_0}\rightarrow A'_{t_0}$
in operator norm and $(1+D_t^2)^{-1}\rightarrow(1+D_{t_0}^2)^{-1}$
in $\mathcal I $-norm by continuity and so the first term
converges to
$$D_{t_0}(1+D_{t_0}^2)^{-1} A'_{t_0}(1+D_{t_0}^2)^{-1}$$
in $\mathcal I $-norm. For the second term we note that
$(1+D_{t_0}^2)^{-1}$ is fixed in  $\mathcal I $ and the rest converges
in operator norm to
$ A'_{t_0}D_{t_0}(1+D_{t_0}^2)^{-1}$ by Corollary 2 of Appendix A of \cite{CP1}.
Thus the derivative of $1-F^2_t$ exists in $\mathcal I $-norm
and equals
$$D_t(1+D_t^2)^{-1} A'_{t}(1+D_{t}^2)^{-1}-
(1+D_t^2)^{-1} A'_{t}D_t(1+D_{t}^2)^{-1}.$$
As a function of $t$ this is continuous in  $\mathcal I $-norm
as $t\mapsto(1+D_t^2)^{-1}$ is  $\mathcal I $-norm continuous
and the rest is operator norm continuous.
\end{proof}

\begin{prop}[{\it cf} Proposition 2.10 of \cite {CP1}]
Let $({\mathcal N}, D_0)$ be weakly $\theta$-summable and let
$t\mapsto A_t$ be a $C^1$ path in ${\mathcal N}_{sa}$
then with $D_t=D_0+A_t$ we have that
$t\mapsto F_t=D_t(1+D_t^2)^{-1/2}$
is a path of Breuer-Fredholm operators in
$F_{D_0}+(Li_0^q)_{F_{D_0}}$ for $0<q<1$
which is $C^1$ in the norm on that space.
Moreover
$$\frac{d}{dt}F_t= \frac{1}{\pi}\int_0^\infty
\lambda^{-1/2}[(1+D_t^2+\lambda)^{-1}(1+\lambda)A_t'
(1+D_t^2+\lambda)^{-1}$$
$$-D_t(1+D_t^2+\lambda)^{-1}A_t'D_t(1+D_t^2+\lambda)^{-1}]
d\lambda$$
where the integral converges in the $Li_0^{q/2}$-norm
and hence also in operator norm.
\end{prop}

\begin{proof}[\bf Proof]
The convergence of the integral follows because each
$(1+D_t^2)^{-1}$ is in $Li$ by Lemma 6.1
and hence in $Li_0^q$ for any $q<1$.
Thus $(1+D_t^2)^{-q/2}$ is in  $Li_0^{q/2}$.
Now we use the proof of Proposition 2.10
of \cite{CP1} by setting, in the notation
used there, $\frac{q}{2}=\frac{1}{2}-\epsilon$
and replacing the ${\mathcal L}^q$-norm by the $Li_0^{q/2}$-norm.
Then the proofs of convergence and continuity go through verbatim
as they depend only on operator norm estimates and the fact that
the ideals in question are invariant operator ideals.

The same remarks apply to the proof of the existence of $\frac{d}{dt}F_t$
in the $Li_0^{q/2}$-norm and the integral formula.
To see that the derivative is continuous in $Li_0^{q/2}$-norm
we need only note that
 $(1+D_t^2)^{-q/2}\leq f(||A_t||)^{q/2}(1+D_0^2)^{-q/2}$
by Lemma 6.1 and operator monotonicity so that the result follows
by our previous remarks (and
$||(1+D_t^2)^{-q/2}||_{Li_0^{q/2}}\leq C
||(1+D_0^2)^{-q/2}||_{Li_0^{q/2}}$ for all $t$).

Finally, for the proof that $t\mapsto F_t$ is
$C^1$ in the norm of  $F_{D_0}+(Li_0^q)_{F_{D_0}}$
we first use Corollary 6.1 to see that it is at least continuous
there. Then we apply Proposition 6.1 with ${\mathcal I}=Li_0^q$
to see that $t\mapsto 1-F_t^2$ is $C^1$ in $Li_0^{q}$
(notice that this does not follow from the product rule as
the $F_t$ themselves are not in $Li_0^{q/2}$). The result
now follows by Lemma B.15.
\end{proof}

\begin{rems*}
Using Lemma B.15 in exactly the same fashion,
we can immediately improve the conclusion of Proposition 2.10 of
\cite {CP1} to read (in the notation of \cite {CP1}) that
$t\mapsto F_t=D_t(1+D_t^2)^{-1/2}$ is $C^1$ in $F_0+
{\mathcal L}_{sa}^{q,\frac{q}{2}}.$ In the notation of this paper,
the latter space is denoted $F_0 + ({\mathcal L}^q)_{F_0}.$

Remark 1.8 of \cite {CP1} is now unnecessary.
\end{rems*}

\begin{prop}
Let $({\mathcal N}, D_0)$ be weakly $\theta$-summable and let
$t\mapsto A_t$ be a $C^1$ path in ${\mathcal N}_{sa}$
then with $D_t=D_0+A_t$ and
$F_t=D_t(1+D_t^2)^{-1/2}$ we have for $0<q<1$ that
the maps
$$t\mapsto F_t'|1-F_t^2|^{-3/2}e^{-|1-F_t^2|^{-1/q}}$$
and $$t\mapsto D_t'e^{-(1+D_t^2)^{1/q}}$$
are both continuous ${\mathcal L}^1$-valued functions of
$t$ and
$$\tau\left( F_t'|1-F_t^2|^{-3/2}e^{-|1-F_t^2|^{-1/q}}\right)
=\tau\left(D_t'e^{-(1+D_t^2)^{1/q}}\right)$$
for all $t$.
\end{prop}

\begin{proof}[\bf Proof]
By Corollary 6.3, $t\mapsto 1-F_t^2 \in Li_0^{q}$
is continuous. Then, by Corollary B.10
$$t\mapsto |1-F_t^2|^{-3/2}e^{-|1-F_t^2|^{-1/q}}\in {\mathcal L}^1$$
is continuous. Apply Proposition 6.5 to see that
$$t\mapsto F_t'|1-F_t^2|^{-3/2}e^{-|1-F_t^2|^{-1/q}}\in {\mathcal L}^1$$
is continuous. Now since
$$t\mapsto (1+D_t^2)^{-1}= 1-F_t^2 \in Li_0^{q}$$
is continuous we have  $t\mapsto e^{-(1+D_t^2)^{1/q}}\in {\mathcal L}^1$
is continuous  by Corollary B.9.
Thus $t\mapsto D_t'e^{-(1+D_t^2)^{1/q}}\in {\mathcal L}^1$
is continuous. Finally
the last claim of the Proposition follows by using the
integral formula for $F_t'$ of Proposition 6.5
and then multiplying through by
$e^{-(1+D_t^2)^{1/q}}=e^{-|1-F_t^2|^{-1/q}}$
to get an integral converging in trace-norm. One then passes the trace through the
integral and uses the cyclicity of the trace so as to allow
an application of Lemma 2.11 of \cite{CP1}.
\end{proof}

\section {\bf SPECTRAL FLOW FORMULAE, UNBOUNDED CASE}

We begin by noting the following result.

\begin{prop}
Assume that $({\mathcal N}, D_0)$ is  an odd unbounded weakly $\theta$-summable
Breuer-Fredholm module for the Banach $*$-algebra $\mathcal A$
and let $F_0=D_0(1+D_0^2)^{-1/2}$. Then  $({\mathcal N}, F_{0})$
is an
odd  $\theta_q$-summable
pre-Breuer-Fredholm module for $\mathcal A$ if $0<q<1$.
\end{prop}

\begin{proof}[\bf Proof]
This follows by an argument similar to that in Proposition 2.4
of \cite{CP1} or alternatively as in \cite{CPS}.
\end{proof}

\begin{defn}
Let $({\mathcal N}, D_0)$ be  an odd unbounded weakly $\theta$-summable
(respectively,\\ $\theta$-summable)
Breuer-Fredholm module for the Banach $*$-algebra $\mathcal A$.
Let ${\mathcal M}_0=D_0+{\mathcal N}_{sa}$
and for $D\in {\mathcal M}_0$ and $X\in T_D({\mathcal M}_0)={\mathcal N}_{sa}$,
the tangent space to  ${\mathcal M}_0$ at $D$, the map
$$\alpha_q(X)=\frac{1}{C_q}\tau\left(Xe^{-(1+D^2)^{1/q}}\right)$$
defines the one-form $\alpha_q$ on  ${\mathcal M}_0$ for $0<q<1$
(respectively, $0<q\leq 1$) where\\
$C_q =C_{0,q} =\int_{-\infty}^\infty e^{-(1+x^2)^{1/q}} dx$.
(We note that $C_1=C_{0,1}=\int_{-\infty}^\infty e^{-(1+x^2)}dx =
 \frac{\sqrt\pi}{e}$)
\end{defn}

We recall from the  definitions of Section 5
that if $F\in F_{D_0}+ (Li_0^q)_{F_{D_0}}$
then
$$\gamma_{\frac{3}{2},q}(F) =
\frac{1}{C_{\frac{3}{2},q}}\int_0^1\tau\left
( F_t'|1-F_t^2|^{-3/2}e^{-|1-F_t^2|^{-1/q}}\right)dt$$
where $\{F_t\}$ is the linear path from $F$ to sign$F$.
We note that
$$C_{\frac{3}{2},q}=\int_{-1}^1(1-u^2)^{-3/2}e^{-(1-u^2)^{-1/q}}du
=\int_{-\infty}^\infty e^{-(1+x^2)^{1/q}} dx=C_q.$$

\begin{lemma}
Let $\beta(X)=\tau(Xg(D)),\  X\in{\mathcal N}_{sa}=T_D({\mathcal M}_0)$
for $D\in{\mathcal M}_0$ be a one-form where
$g:{\mathcal M}_0\rightarrow {\mathcal L}^1$
is continuous and the integral of
$\beta$ is independent of the piecewise $C^1$ path
in ${\mathcal M}_0$. Then $\beta = df$
where $f(D)=\int_0^1\tau(D_t'g(D_t))dt$
and $\{D_t\}$ is any such path in
${\mathcal M}_0$ from $D_0$ to $D$. That is, $\beta$ is exact.
\end{lemma}

\begin{proof}[\bf Proof]
Recall $df_D(X)=\frac{d}{ds}|_{s=0}(f(D+sX))$. For each
$s$ choose our path from $D_0$ to $D+sX$ to pass through
$D$ and be linear from $D$ to  $D+sX$
and be indexed by $r\in [0,s]$ (or $[s,0]$ if $s<0$).
Then $f(D+sX)=f(D)+\int_0^s\tau(Xg(D+rX))dr$ and therefore
$$\frac{d}{ds}|_{s=0}(f(D+sX))=\tau(Xg(D))=\beta(X)$$
as claimed.
\end{proof}

\begin{thm}
Let $({\mathcal N}, D_0)$ be an odd unbounded weakly $\theta$-summable
Breuer-Fredholm module for the Banach $*$-algebra $\mathcal A$.
Let ${\mathcal M}_0=D_0+{\mathcal N}_{sa}$ then for
$0<q<1$ the integral of the one-form
$$\alpha_q(X)=\frac{1}{C_q}\tau\left(Xe^{-(1+D^2)^{1/q}}\right)$$
is independent of the path in  ${\mathcal M}_0$
and hence $\alpha_q$ is exact. Moreover, if
$\{D_t\}_{t\in[a,b]}$ is any piecewise $C^1$ path
in  ${\mathcal M}_0$ then
$$sf(D_a,D_b)=\frac{1}{C_q}\int_a^b \tau\left(D_t'e^{-(1+D_t^2)^{1/q}}\right)dt
+\gamma_{\frac{3}{2},q}(D_b(1+D_b^2)^{-1/2})-
\gamma_{\frac{3}{2},q}(D_a(1+D_a^2)^{-1/2}).$$
\end{thm}

\begin{proof}[\bf Proof]
We know that with $F_0=D_0(1+D_0^2)^{-1/2}$, $({\mathcal N}, F_0)$
is $\theta_q$-summable
by Proposition 7.1. Moreover
$t\mapsto F_t=D_t(1+D_t^2)^{-1/2}$
is a piecewise $C^1$ path in $F_0+ (Li_0^q)_{F_0}$
by Proposition 6.5. Now we recall from definition 2.15 of \cite {CP1}
that
$sf(D_a,D_b)= sf(F_a,F_b)$
and so by Theorem 5.1
(with $r=3/2$) together with Proposition 6.6
we obtain our formula. It follows from this
formula that the integral of the one-form $\alpha_q$
is independent of the path in ${\mathcal M}_0$.
\end{proof}

\begin{cor}
Let $({\mathcal N}, D_0)$ be an odd unbounded weakly
$\theta$-summable
Breuer-Fredholm module for the Banach $*$-algebra $\mathcal A$,
and let $0<q<1$. Let $P=\chi_{[0,\infty)}(D_0)$, then
for each unitary $u\in \mathcal A$ with
$u(dom D_0)\subseteq dom(D_0)$
and $[D_0,u]$ bounded we have that $PuP$ is a Breuer-Fredholm
operator in $P{\mathcal N}P$ and if $\{D_t^u\}$
is any piecewise $C^1$
 path in  ${\mathcal M}_0=D_0+{\mathcal N}_{sa}$
from $D_0$ to $uD_0u^*$ (for example the linear path)
then
$$\text{ind}(PuP)=sf(\{D_t^u\})=
\frac{1}{C_q}\int_0^1 \tau\left((D_t^u)'e^{-(1+(D^u_t)^2)^{1/q}}\right)dt.
$$
\end{cor}

\begin{proof}[\bf Proof]
Since
$$uDu^*=D_0-[D_0,u]u^*\in D_0+{\mathcal N}_{sa}$$
the right hand equality follows from the previous theorem
since
$$\gamma_{\frac{3}{2},q}(uD_0(1+D_0^2)^{-1/2}u^*)=
\gamma_{\frac{3}{2},q}(D_0(1+D_0^2)^{-1/2}).$$
For $F_0=D_0(1+D_0^2)^{-1/2}$ we have that  $({\mathcal N}, F_0)$ is
$\theta_q$-summable
by Proposition 7.1 and so by Theorem 5.9
$ind(PuP)$ is given by the spectral flow along the linear
path from $F_0$ to $uF_0u^*$ and hence is the
spectral flow of any piecewise $C^1$ path from  $F_0$ to $uF_0u^*$.
In particular,
$$ind(PuP)=sf(\{D_t^u(1+(D^u_t)^2)^{-1/2}\})= sf(\{D_t^u\})$$
\end{proof}

\begin{rems*}
Let $({\cn},D_0)$ be $\theta$-summable (for $\comp$ say) and for each
$D\in {\cm}_0=D_0+{\cn}_{sa}$ we let $F_D=D(1+D^2)^{-1/2}$. Then
$\gamma_{\frac{3}{2},1}(F_D)$ is well-defined even though we do not know
whether $F_D \in F_{D_0} + (Li_0)_{F_{D_0}}.$ This follows from the fact
that $1-F_D^2=(1+D^2)^{-1}\in Li_0$ by Lemma 6.1, and
that the definition of $\gamma_{\frac{3}{2},1}$ in Section 5 only involves
$F_D$ and $\tilde F = sign(F_D).$
\end{rems*}

$\spadesuit$ The following lemma finally allows us to get rid of the 
annoying $q$ in our formula when our module is actually $\theta$-summable.
$\spadesuit$

\begin{lemma}
Let $({\mathcal N}, D_0)$ be  an odd unbounded $\theta$-summable
Breuer-Fredholm module for the Banach $*$-algebra $\mathcal A$,
and let $\{D_t\}_{t\in[a,b]}$ be a piecewise
$C^1$ path in ${\mathcal M}_0$.
Then, for each $D\in {\mathcal M}_0$
where $F_D=D(1+D^2)^{-1/2}$,
$$\lim_{q\to 1^-}\int_0^1 \tau\left(D_t'e^{-(1+D_t^2)^{1/q}}\right)dt
=\int_0^1 \tau\left(D_t'e^{-(1+D_t^2)}\right)dt,  \eqno (7.1)$$
$$\lim_{q\to 1^-}\gamma_{\frac{3}{2},q}(F_D)
=\gamma_{\frac{3}{2},1}(F_D), \eqno(7.2) $$
$$\lim_{q\to 1^-}C_q=C_1.   \eqno(7.3)$$
\end{lemma}

\begin{proof}[\bf Proof]
The proof rests on a subsidiary result. Suppose that
$X\in Li_0$ and $|X|\leq 1$. Then, for $0<q<1$,\;\;
$0\leq e^{-|X|^{-1/q}}\leq e^{-|X|^{-1}}$
and so
$$0\leq|X|^{-r} e^{-|X|^{-1/q}}\leq |X|^{-r}e^{-|X|^{-1}}$$
for any $r\geq 0$.
By Corollary B.11 these operators are
in ${\mathcal L}^1$ and so
\begin{eqnarray}
&   &||\ |X|^{-r}e^{-|X|^{-1}}-|X|^{-r} e^{-|X|^{-1/q}}||_1\nonumber\\
& = & \tau(|X|^{-r}e^{-|X|^{-1}}-|X|^{-r} e^{-|X|^{-1/q}}) \nonumber\\
& = & \tau(|X|^{-r}e^{-|X|^{-1}})-\tau(|X|^{-r} e^{-|X|^{-1/q}})\nonumber\\
& = & \int_0^\infty \mu_s(X)^{-r}\left[e^{- \mu_s(X)^{-1}}
     -e^{- \mu_s(X)^{-1/q}}\right]ds\nonumber
\end{eqnarray}
by \cite{FK} Remark 3.3. Now the integrand converges pointwise to
0 as $q\to 1^-$ and is dominated by the integrable function
$s\mapsto  \mu_s(X)^{-r}e^{- \mu_s(X)^{-1}}$. Thus the integral
converges to 0 as  $q\to 1^-$. That is
 $$\lim_{q\to 1^-}||\ |X|^{-r}e^{-|X|^{-1}}-|X|^{-r} e^{-|X|^{-1/q}}||_1
=0 \eqno(7.4).$$
Now, to see (7.1):
$$\left|\int_0^1 \tau\left(D_t'e^{-(1+D_t^2)}\right)dt
-\int_0^1 \tau\left(D_t'e^{-(1+D_t^2)^{1/q}}\right)dt
\right|$$
$$\leq \int_0^1||D_t'||.||e^{-(1+D_t^2)}-e^{-(1+D_t^2)^{1/q}}||_1dt
\leq C \int_0^1||e^{-(1+D_t^2)}-e^{-(1+D_t^2)^{1/q}}||_1dt.$$
By (7.4) with $X=(1+D_t^2)^{-1}$ and $r=0$ we see that the integrand goes
to 0 pointwise in $t$. However
$$||e^{-(1+D_t^2)}-e^{-(1+D_t^2)^{1/q}}||_1
=\tau\left(e^{-(1+D_t^2)}-e^{-(1+D_t^2)^{1/q}}\right)\leq \tau\left
(e^{-(1+D_t^2)}\right).$$
But
$$t\mapsto (1+D_t^2)^{-1}\mapsto e^{-(1+D_t^2)}\in {\mathcal L}^1$$
is continuous by Lemma 6.1 and Corollary B.11.
Hence,
$t\mapsto \tau\left(e^{-(1+D_t^2)}\right)$ is integrable. Thus,
$$\int_0^1||e^{-(1+D_t^2)}-e^{-(1+D_t^2)^{1/q}}||_1dt\to 0$$
as
$q\to 1^-$ and (7.1) follows.

Next we observe that (7.3) is an easy application of the
dominated convergence theorem as
$C_q =\int_{-\infty}^\infty e^{-(1+x^2)^{1/q}} dx$.
Since $C_q =C_{\frac{3}{2},q} $
we also get $\lim_{q\to 1^-}C_{\frac{3}{2},q}=C_{\frac{3}{2},1}$.

Finally, we obtain (7.2) by a method very similar to that used for
 (7.1) using (7.3) and recalling that
$$\gamma_{\frac{3}{2},q}(F_D) =
\frac{1}{C_{\frac{3}{2},q}}\int_0^1\tau\left
( F_t'|1-F_t^2|^{-3/2}e^{-|1-F_t^2|^{-1/q}}\right)dt$$
where $F_t=(1-t)F_D+t\tilde F_D$ is the linear path
and $1-F_t^2\in Li_0$ with $|1-F_t^2|\leq 1$.
\end{proof}

\begin{defn}
Let $({\mathcal N}, D_0)$ be  an odd unbounded $\theta$-summable
Breuer-Fredholm module for the Banach $*$-algebra $\mathcal A$
and let $D\in {\mathcal M}_0$. Also let
$$\gamma_0(D)=\gamma_{\frac{3}{2},1}(D(1+D^2)^{-1/2})$$
and for $X\in T_D( {\mathcal M}_0)={\cn}_{sa}$ define a one-form $\alpha$ on
${\mathcal M}_0$ by
$$\alpha(X)=\frac{1}{C_1}\tau\left(Xe^{-(1+D^2)}\right)=\frac{1}{\sqrt\pi}
\tau\left(Xe^{-D^2}\right).$$
\end{defn}

\begin{thm}
Let $({\mathcal N}, D_0)$ be  an odd unbounded $\theta$-summable
Breuer-Fredholm module for the Banach $*$-algebra $\mathcal A$.
Then the integral of the one-form $\alpha$ is independent
of the path in  ${\mathcal M}_0$ so that $\alpha$ is exact
and moreover if $\{D_t\}_{t\in[a,b]}$ is any
piecewise $C^1$ path in ${\mathcal M}_0$ then
$$sf(D_a,D_b)=\frac{1}{\sqrt\pi}\int_a^b \tau\left(D_t'e^{-D_t^2}\right)dt
+\gamma_0(D_b)-\gamma_0(D_a).$$
\end{thm}

\begin{proof}[\bf Proof]
As $\theta$-summable implies weakly
 $\theta$-summable the last formula follows from Theorem 7.4,
Lemma 7.6 and the definition preceding the theorem. That the integral
of $\alpha$ is independent of the path in
 ${\mathcal M}_0$
now follows from this formula.
\end{proof}

\noindent{\bf Remark}: Note that the methods of Appendix C may also be used 
to prove directly that $\alpha$ is closed. 
Then because our space ${\mathcal M}_0$
is affine we can deduce by a Poincar\'e lemma style argument that $\alpha$
is exact. This direct proof is not appreciably shorter and we will
discuss it elsewhere \cite{CPRS2}.

\begin{cor}
Let $({\mathcal N}, D_0)$ be  an odd unbounded $\theta$-summable
Breuer-Fredholm module for the Banach $*$-algebra $\mathcal A.$
Let $P=\chi_{[0,\infty)}(D_0)$, then
for each unitary $u\in \mathcal A$ with
$u(dom D_0)\subseteq dom D_0$
and $[D_0,u]$ bounded we have that $PuP$ is a Breuer-Fredholm
operator in $P{\mathcal N}P$ and if $\{D_t^u\}$
is any piecewise $C^1$ path in  ${\mathcal M}_0=D_0+{\mathcal N}_{sa}$
from $D_0$ to $uD_0u^*$ (for example the linear path)
then
$$ind(PuP)=sf(\{D_t^u\})=
\frac{1}{\sqrt\pi}\int_0^1 \tau\left((D_t^u)'e^{-(D_t^u)^2}\right)dt.
$$
\end{cor}

\begin{proof}[\bf Proof]
See the proof of Corollary 7.5.
\end{proof}

\begin{cor}
Let $({\cn},D_0)$ be an odd unbounded $\theta$-summable Breuer-Fredholm
module for the Banach $*$-algebra $\ca.$ For any $\epsilon >0$ we define
a one-form, $\alpha^{\epsilon}$ on ${\cm}_0=D_0+{\cn}_{sa}$ by
$$\alpha^{\epsilon}(X)=\sqrt{\frac{\epsilon}{\pi}}\tau\left(Xe^{-\epsilon D^{2}}
\right)$$
for $D\in {\cm}_0$ and $X\in T_D({\cm}_0)={\cn}_{sa}.$ Then the integral
of $\alpha^{\epsilon}$ is independent  of the piecewise $C^1$ path in
${\cm}_0$
and if $\{D_t\}_{t\in [a,b]}$ is any piecewise $C^1$ path in ${\cm}_0$ then
$$sf(D_a,D_b)=
\sqrt{\frac{\epsilon}{\pi}}\int_a^b \tau\left(D_t'e^{-\epsilon D_t^2}\right)dt
+ \gamma_0(\sqrt{\epsilon}D_b) - \gamma_0(\sqrt{\epsilon}D_a).$$
\end{cor}

\begin{proof}[\bf Proof]
 Clearly, $({\cn},\sqrt{\epsilon}D_0)$ is also a
$\theta$-summable module for ${\ca}$ and $sf(\sqrt{\epsilon}D_a,
\sqrt{\epsilon}D_b) = sf(D_a,D_b)$ so that this is immediate from
Theorem 7.8.
\end{proof}

\section{\bf ETA INVARIANTS}

Corollary 7.10 is very similar to Theorem 2.6 of \cite{G} with some important
differences. First, Getzler's theorem is for the type $I_\infty$ case only; second,
he assumes that the endpoints of the path, $D_a$ and $D_b$ are invertible;
and third, his correction terms are truncated $\eta$-invariants. It is the
purpose of this section to show that his correction terms are identical
to ours when the endpoints are invertible, even in the type $II_\infty$ setting.
Moreover, when the endpoints are not invertible, we show how to modify
the truncated $\eta$-invariants to get the right correction terms:
namely, $\gamma_0(\sqrt{\epsilon}D)$.

First we define the truncated $\eta$-invariants and show that they make sense
in our general setting.

\begin{defn}
If $D$ is an unbounded self-adjoint operator affiliated with ${\cn}$, and
$e^{-tD^{2}}$ is trace-class for all $t>0$ (briefly,
$D$ is $\theta$-summable relative to $\cn$) then we define
$$\eta_{\epsilon}(D) = \frac{1}{\sqrt{\pi}}\int_{\epsilon}^{\infty}
\tau\left(De^{-tD^{2}}\right)t^{-1/2}dt.$$
\end{defn}

We first observe that $De^{-tD^{2}} = De^{-(t/2)D^{2}}e^{-(t/2)D^{2}}$
where the first factor is a bounded operator (in ${\cn}$) and the second
factor is trace-class by hypothesis. Thus the integrand is finite-valued
for each $t>0$. Moreover by the functional calculus, the map
$t\mapsto De^{-(t/2)D^{2}}$ is operator norm continuous. As the second
term equals $e^{t/2}e^{-(t/2)(1+D^{2})}$ and $(1+D^2)^{-1}\in Li_0$ by
Corollary B.6, we have
that the second term is trace-class continuous by Corollary B.9.
Thus the integrand is a continuous real-valued function. To see that the
integral converges, we first prove the following.

\begin{lemma}
Let $D$ be an unbounded self-adjoint operator affiliated with $\cn$
such that \\$(1+D^2)^{-1} \in {\mathcal K_{\cn}}.$ Let $\{E_{\lambda}\}$
denote the spectral resolution of $|D|$ and suppose $f:\real^{+} \to
\real^{+}$ is continuous and $f(|D|)$ is trace-class. Then,
$$\tau(f(|D|))=\int_0^{\infty}f(\lambda)d\phi_{\lambda}$$
where $\phi_{\lambda}=\tau(E_{\lambda}).$
\end{lemma}

\begin{proof}[\bf Proof]
By spectral theory, $(1+D^2)^{-1}\geq\frac{1}{1+\lambda^2}E_{\lambda}$
so that $\phi_{\lambda}$ is a finite-valued increasing function. If we fix
$\lambda_0 > 0$, then $f(|D|)E_{\lambda_0}$ lies in the $II_1$ algebra
$E_{\lambda_0}{\cn}E_{\lambda_0}$, on which $\tau$ is operator-norm
continuous. Then it is easy to see that
$$\tau(f(|D|)E_{\lambda_0})=
\tau\left(\int_0^{\lambda_0}f(\lambda)dE_{\lambda}\right)
=\int_0^{\lambda_0}f(\lambda)d\phi_{\lambda}.$$
Now, as $\lambda_0 \to \infty$ the RHS approaches $\int_0^\infty f(\lambda)
d\phi_{\lambda}$. By the lower semicontinuity of $\tau$ we get
$$\tau(f(|D|)\leq \liminf_{\lambda_0 \to \infty}
 \tau(f(|D|)E_{\lambda_0}) \leq \limsup_{\lambda_0 \to \infty}
\tau(f(|D|)E_{\lambda_0}) \leq \tau(f(|D|),$$
and we are done.
\end{proof}

A version of the following result is implicit in \cite{M}.

\begin{lemma}
If $D$ is $\theta$-summable relative to ${\cn}$, then
the integral $$\int_1^\infty \tau\left(|D|e^{-tD^2}\right)t^{-1/2}dt$$
converges.
\end{lemma}

\begin{proof}[\bf Proof]
We denote the spectral resolution of $|D|$
by $\{E_\lambda\}$.
We let $\phi_\lambda= \tau(E_\lambda)$ then by the previous lemma,
$$\int_1^\infty \tau(|D|e^{-tD^2})t^{-1/2}dt
=\int_1^\infty \int_0^\infty \lambda e^{-t\lambda^2}
 d\phi_\lambda t^{-1/2} dt.$$
To see that this double integral converges we first use Tonelli's theorem
to interchange the order of integration
$$\int_0^\infty\int_1^\infty \lambda e^{-t\lambda^2}t^{-1/2} dt d\phi_\lambda
=\int_0^\infty e^{-\lambda^2}
\int_1^\infty \lambda e^{-(t-1)\lambda^2}t^{-1/2} dt d\phi_\lambda.$$
Make the substitution $v=(t-1)\lambda^2$ then we see that our double integral
\begin{eqnarray}
&=&\int_0^\infty\int_0^\infty e^{-v}(v+\lambda^2)^{-1/2} e^{-\lambda^2}
dv d\phi_\lambda\nonumber\\
&\leq& \int_0^\infty e^{-v}v^{-1/2}dv \int_0^\infty  e^{-\lambda^2}
d\phi_\lambda\nonumber
\end{eqnarray}
which is finite
as required since $\int_0^{\infty}e^{-\lambda^{2}}d\phi_{\lambda} =
\tau(e^{-D^2})$.
\end{proof}

\begin{cor}
If $D$ is $\theta$-summable relative to ${\cn}$ then for each $\epsilon>0$,
the integral
$$\eta_\epsilon (D) = \frac{1}{\sqrt{\pi}}\int_{\epsilon}^{\infty}
\tau\left(De^{-tD^2}\right)t^{-1/2}dt$$ converges.
\end{cor}

\begin{proof}[\bf Proof]
By replacing $D$ with $\sqrt{\epsilon}D$, we can take $\epsilon = 1.$
Convergence now follows from the previous lemma since $|\tau\left(
De^{-tD^2}\right)| \leq \tau\left(|D|e^{-tD^2}\right).$
\end{proof}

In order to reconcile $\eta_\epsilon(D)$ with
the correction terms, $\gamma_0(\sqrt{\epsilon}D)$ of Corollary 7.10,
we recall that $$\gamma_0(\sqrt{\epsilon}D) = \gamma_{\frac{3}{2},1}
(\sqrt{\epsilon}D(1+{\epsilon}D^2)^{-1/2}).$$
Since $\sqrt{\epsilon}D(1+{\epsilon}D^2)^{-1/2})=D(\frac{1}{\epsilon}
+D^2)^{-1/2}$,
we are led to consider the transformations
$$F_s=D(s+D^2)^{-1/2}\; for\; s>0.$$

\begin{lemma}
Let $D$ be an unbounded self-adjoint operator and let $r>0$,
$s>0$. Then,\\
(1) $||\left[(s+D^2)^{1/2}+(r+D^2)^{1/2}\right]^{-1}||\leq\frac{1}
{\sqrt{s} + \sqrt{r}},$ and\\
(2) $||\left[(s+D^2)^{1/2}+(r+D^2)^{1/2}\right]^{-1}-
\left[2(s+D^2)^{1/2}\right]^{-1}||\leq\frac{1}{\sqrt{s}+\sqrt{r}}
|\sqrt{r}-\sqrt{s}|\frac{1}{2\sqrt{s}}.$
\end{lemma}

\begin{proof}[\bf Proof]
Item (1) follows from the functional calculus and the numerical inequality:
$$\frac{1}{\sqrt{s+x^2}+\sqrt{r+x^2}}\leq\frac{1}{\sqrt{s}+\sqrt{r}}$$
for all real $x$.\\
Item (2) follows from the numerical identity:
$$\frac{1}{\sqrt{s+x^2}+\sqrt{r+x^2}}-\frac{1}{2\sqrt{s+x^2}}=
\left(\frac{\sqrt{s+x^2}-\sqrt{r+x^2}}{\sqrt{s+x^2}+\sqrt{r+x^2}}\right)
\left(\frac{1}{2\sqrt{s+x^2}}\right)$$
and the easily proved estimate:
$$|\sqrt{s+x^2}-\sqrt{r+x^2}|\leq|\sqrt{s}-\sqrt{r}|.$$
\end{proof}

\begin{prop}
If $D$ is $\theta$-summable relative to $\cn$ and $F_s=D(s+D^2)^{-1/2}$
for all
$s>0$ then $F_s\in F_D+(Li_0)_{F_D}$ and the mapping $s\mapsto F_s$
is $C^1$ in this space. Moreover, $\frac{d}{ds}(F_s)=-\frac{1}{2}
F_s(s+D^2)^{-1}.$
\end{prop}

\begin{proof}[\bf Proof]
We first observe that for $s>0,\;(s+D^2)^{-1}=(1/s)(1+(1/s)D^2)^{-1}$ is
in $Li_0$ by Corollary B.6. Now we fix $s>0$. Then,
\begin{eqnarray}
F_r-F_s & = & D(r+D^2)^{-1/2}\left[(s+D^2)^{1/2}-(r+D^2)^{1/2}\right]
               (s+D^2)^{-1/2}\nonumber\\
        & = & F_r\left[\{(s+D^2)-(r+D^2)\}\{(s+D^2)^{1/2}+(r+D^2)^{1/2}\}^{-1}
               \right](s+D^2)^{-1/2}\nonumber\\
        & = & F_r\left[(s-r)\left\{(s+D^2)^{1/2}+(r+D^2)^{1/2}\right\}^{-1}
               \right](s+D^2)^{-1/2}.\nonumber
\end{eqnarray}
So, $F_s-F_r\in Li_0^{1/2}$ and in particular, $F_s-F_D\in Li_0^{1/2}.$
Since $1-F_s^2=s(s+D^2)^{-1}\in Li_0$ we see that
$F_s\in F_D+(Li_0)_{F_D}$ as claimed.

Now, $$\frac{1}{r-s}(F_r-F_s)=-F_r\left[(s+D^2)^{1/2}+(r+D^2)^{1/2}\right]^{-1}
(s+D^2)^{-1/2}$$
and so by the estimates of the previous lemma and the fact that $r\mapsto F_r$
is (at least!) operator norm continuous by the previous equations, we can take
the limit as $r\to s$ in the norm of $Li_0^{1/2}$ to get:
$$\frac{d}{ds}(F_s)=-\frac{1}{2}F_s(s+D^2)^{-1}\;in\;the\;Li_0^{1/2}\;sense.$$
It is easily seen that $s\mapsto -(1/2)F_s(s+D^2)^{-1}$ is continuous in the
norm of $Li_0$ and {\it a fortiori} in the norm of $(Li_0)_{F_D}.$
To see that this derivative exists in the sense of the norm
on $(Li_0)_{F_D},$ it suffices (by Lemma B.15) to see that
$$s\mapsto (1-F_s^2)=s(s+D^2)^{-1}\in Li_0 \; is\; C^1.$$
Clearly it suffices to see that $s\mapsto (s+D^2)^{-1}$ is $C^1$ in $Li_0.$
This is an easy resolvent equation calculation. In fact:
$$\frac{d}{ds}(s+D^2)^{-1}=-(s+D^2)^{-2}.$$
\end{proof}

\begin{lemma}
If $D$ is $\theta$-summable relative to $\cn$ and $\epsilon >0$, then
$$-\frac{1}{2}\eta_{\epsilon}(D)=\frac{1}{C_{\frac{3}{2},1}}\int_0^1\tau
\left(\frac{d}{ds}(F_s)(1-F_s^2)^{-\frac{3}{2}}e^{-(1-F_s^2)^{-1}}\right)ds$$
where for $s\geq0$, $F_s=\sqrt{\epsilon}D(s+{\epsilon}D^{2})^{-1/2}.$
\end{lemma}

\begin{proof}[\bf Proof]
Since $\eta_{\epsilon}(D)=\eta_{1}(\sqrt{\epsilon}D)$, we can let
$\epsilon=1$ by replacing $\sqrt{\epsilon}D$ with $D.$ Now,
\begin{eqnarray}
-\frac{1}{2}\eta_{1}(D)&=&-\frac{1}{2\sqrt{\pi}}\int_1^\infty \tau \left(
t^{-1/2}De^{-tD^2}\right)dt\nonumber\\
&=&-\frac{e}{\sqrt{\pi}}\int_1^\infty\tau\left(\frac{1}{2}t^{-1/2}D
e^{-(1+tD^2)}\right)dt.\nonumber
\end{eqnarray}
Since this integral converges absolutely, we can make the (scalar)
change of variable $t=1/s$ to obtain the absolutely convergent integral:
\begin{eqnarray}
-\frac{1}{2}\eta_{1}(D) & = & -\frac{e}{\sqrt{\pi}}\int_0^1\tau\left(\frac{1}{2}
                               s^{-3/2}De^{-(1+(1/s)D^2)}\right)ds\nonumber\\
                        & = & -\frac{e}{\sqrt{\pi}}\int_0^1\tau\left(
                               \frac{1}{2}s^{-3/2}D(1+(1/s)D^2)^{-3/2}
                               (1+(1/s)D^2)^{3/2}e^{-(1+(1/s)D^2)}\right)ds
                               \nonumber\\
                        & = & \frac{e}{\sqrt{\pi}}\int_0^1\tau\left(
                               \frac{d}{ds}(F_s)(1-F_s^2)^{-3/2}
                               e^{-(1-F_s^2)^{-1}}\right)ds\nonumber
\end{eqnarray}
by the previous proposition and the fact that $1-F_s^2=s(s+D^2)^{-1}$.
We have previously noted that $C_{\frac{3}{2},1}=\frac{\sqrt{\pi}}{e}$
so we are done.
\end{proof}

\begin{rems*}
Because the path $F_s = D(s+D^2)^{-1/2}$ is not continuous at zero
in the $(Li_0)_{F_D}$-norm, in order
to prove that this latter integral equals the integral of
the one-form $\alpha_{\frac{3}{2}}$ along some (any) $C^1$ path from
$F_0$ to $F_1=F_D$ we cannot just appeal to the exactness of our
one-form since we are integrating along a discontinuous path.
To overcome this we argue as follows.

First we truncate our path at $\delta>0$ where $\delta$ is small.
Then we have $\{F_s\}$ for $\delta\leq s\leq 1$
is a $C^1$ path joining $F_{\delta}$ and $F_1=F_D$. If we extend
this path at its
beginning with the straight line path from $F_0$ to $F_{\delta}$,
we obtain a piecewise $C^1$ path from $F_0$ to $F_D$.
Thus the integral of our one-form along this new path is the same
as the integral of our one-form along any $C^1$ path from
$F_0$ to $F_D$. For small $\delta$, the piece we have thrown away
is small by the absolute convergence of the integral. To complete the
argument we must show that the piece we have added, namely the integral of
our one-form along the straight line from $F_0$ to $F_{\delta}$ is also
small. This is not obvious, since in the generic type $II_\infty$ case,
we would have
$||F_0-F_{\delta}||=1$ for all $\delta>0$, so that $F_0$ and $F_{\delta}$
would be even farther apart in $(Li_0)_{F_D}$-norm.
\end{rems*}

\begin{not*}
For the purposes of the rest of this section we will use the notation
$\theta(F_1,F_2)$ to denote the integral of the one-form
$\alpha_{\frac{3}{2}}$ from $F_1$ to $F_2$ along a piecewise $C^1$
path in $F_D+(Li_0)_{F_D}.$
\end{not*}

\begin{lemma}
If $D$ is $\theta$-summable relative to $\cn$ and $0<\delta<1$,
then $$\lim_{\delta\to 0}\theta(F_0,F_{\delta})=0,$$
where $F_{\delta}=D(\delta+D^2)^{-1/2}.$
\end{lemma}

\begin{proof}[\bf Proof]
Let $F_{\delta,t}=F_0 + t(F_{\delta}-F_0)$
for $t\in [0,1]$ be the
straight line path from $F_0$ to $F_{\delta}$.
Then $$\theta(F_0,F_{\delta})=\frac{1}{C_{\frac{3}{2},1}}
\int_0^1\tau\left((F_{\delta}-F_0)(1-F_{\delta,t}^2)^{-3/2}
e^{-(1-F_{\delta,t}^2)^{-1}}\right)dt.$$
We observe that the operator in the integrand is $0$ on $ker(D)$
because of the term $(F_{\delta}-F_0)$ and so all of the functions
of $D$ can be regarded as being restricted to $ker(D)^{\perp}$.
That is, for the purposes of this calculation, we can (and do) assume that
$ker(D)=\{0\}.$
With this in mind, we factor the operator in the integrand into three pieces:

$(F_{\delta}-F_0)$,

$(1-F_{\delta,t}^2)^{-3/2}e^{-\frac{1}{2}(1-F_{\delta,t}^2)^{-1}},\; and$

$e^{-\frac{1}{2}(1-F_{\delta,t}^2)^{-1}}.$

The first factor is operator-norm bounded by 1. The second factor is
operator-norm bounded (independent of $t$ and $\delta$) by:
$$\sup_{x\in [0,1]}\left[x^{-\frac{3}{2}}e^{-\frac{1}{2x}}\right]=
\left(\frac{3}{e}\right)^{3/2}<1.6.$$
The third factor is bounded as a positive operator by:
$$e^{-\frac{1}{2}(1-F_{\delta}^2)^{-1}}.$$
Thus, it suffices to see that
$$||e^{-\frac{1}{2}(1-F_{\delta}^2)^{-1}}||_{1}=\tau\left(e^{-\frac{1}{2}
(1-F_{\delta}^2)^{-1}}\right)\to 0\;as\;\delta\to 0.$$

Now, for $0<\delta\leq 1,$
$$1-F_{\delta}^2=\delta(\delta+D^2)^{-1}=\delta(1+D^2)^{-1}\left[1-(1-\delta)
(1+D^2)^{-1}\right]^{-1},$$
and since $f(x)=\delta x[1-(1-\delta)x]^{-1}$ is an increasing function
of $x$ for $x\in [0,1]$, we have by part (iv), Lemma 2.5 of \cite{FK} that
$$\mu_t(1-F_{\delta}^2)=\delta\mu_t((1+D^2)^{-1})\left[1-(1-\delta)
\mu_t((1+D^2)^{-1})\right]^{-1}.$$
Since we are assuming $kerD=\{0\}$, we have for each fixed
$t>0$ that $0\leq \mu_t((1+D^2)^{-1})$ is strictly less than $1$ and so:
$$\lim_{\delta \to 0}\mu_t(1-F_{\delta}^2)=0.$$
Another application of part (iv), Lemma 2.5 of \cite{FK} gives us:
$$\lim_{\delta \to 0}\mu_t \left(e^{-\frac{1}{2}(1-F_{\delta}^2)^{-1}}\right)
=\lim_{\delta \to 0}\left(e^{-\frac{1}{2}[\mu_t
(1-F_{\delta}^2)]^{-1}}\right)=0.$$
Therefore, by Corollary 2.8 of \cite{FK} and the Lebesgue Dominated
Convergence Theorem:
$$\tau\left(e^{-\frac{1}{2}(1-F_{\delta}^2)^{-1}}\right)=
\int_0^\infty e^{-\frac{1}{2}[\mu_t(1-F_{\delta}^2)]^{-1}}dt
\to 0 \;as\;\delta \to 0.$$
This completes the proof.
\end{proof}

\begin{thm}
If $D$ is $\theta$-summable relative to $\cn$ and $\epsilon>0$, then
$$\frac{1}{2}\eta_{\epsilon}(D)=\gamma_0(\sqrt{\epsilon}D)-
\frac{1}{2}\tau([ker(D)]),$$
where $[ker(D)]$ is the projection on $ker(D).$
\end{thm}

\begin{proof}[\bf Proof]
Since $\eta_{\epsilon}(D)=\eta_1(\sqrt{\epsilon}D)$, we can assume that
$\epsilon=1.$ Combining Lemma 8.7, the Remarks, and Lemma 8.8 we now have:
\begin{eqnarray}
\frac{1}{2}\eta_1(D) & = & -\theta(F_0,F_D) \nonumber \\
                     & = & \theta(F_D,F_0) \nonumber \\
                     & = & \theta(F_D,\tilde F) - \theta(F_0,\tilde F)
                     \nonumber \\
                     & = & \gamma_{\frac{3}{2},1}(F_D)-\gamma_{\frac{3}{2},q}
                     (F_0) \nonumber \\
                     & = & \gamma_0(D)-\gamma_{\frac{3}{2},1}(F_0). \nonumber
\end{eqnarray}
Since $\tilde F-F_0=[ker(D)]$, it is an easy calculation that:
$$\gamma_{\frac{3}{2},1}(F_0)=\frac{1}{2}\tau([ker(D)]),$$
and we're done. Another explanation of this last equality which does
not directly involve calculating an integral is the following.
Let $E=[ker(D)]$ and let $\hat{F}=\tilde{F}-2E$. Then,
$\hat{F}$ and $\tilde{F}$ are unitarily equivalent and clearly,
$sf(\hat{F},\tilde{F})=\tau(E)$ is the integral of the one-form
$\alpha_{\frac{3}{2}}$ from $\hat{F}$ to $\tilde{F}$ (by Theorem 5.8).
Now, $F_0$ lies in a position of symmetry exactly half way
between $\hat{F}$ and $\tilde{F}$, and so the integrals from $\hat{F}$ to $F_0$
and $F_0$ to $\tilde{F}$ are identical with sum $\tau(E).$ Since
$\gamma_{\frac{3}{2},1}(F_0)$ is by definition the integral from
$F_0$ to $\tilde{F}$, we see that it is exactly $\frac{1}{2}\tau(E)$ as claimed.
\end{proof}

$\spadesuit$ The $\eta$-invariant is focussed on the spectral asymmetry of 
the operator
$D$, and so it treats $0$ in a symmetric manner: $\eta_1(D)$ is an integral
along a path connecting $F_D$ to $sgn(D)$ where the {\it signum} function,
$sgn$ has the value $0$ at $0.$ On the other hand, $\gamma$ is concerned
directly with spectral flow: it is the integral along a path from
$F_D$ to a ``universal'' symmetry (up to
unitary equivalence) associated with $D$ and so a natural choice is
the {\it signum} function $sign$ which takes the value $1$ at $0$. $\spadesuit$

Combining Corollary 7.10 and the previous theorem, we can now
deduce the the promised generalization of Theorem 2.6 of \cite {G}.
We again observe that it applies equally to the type $II_\infty$ situation and
that we do not have to assume invertibility of the endpoints.

\begin{cor}
Let $({\cn},D_0)$ be an odd unbounded $\theta$-summable Breuer-Fredholm
module for the Banach $*$-algebra $\ca.$ For any $\epsilon >0$ we define
a one-form, $\alpha^{\epsilon}$ on ${\cm}_0=D_0+{\cn}_{sa}$ by
$$\alpha^{\epsilon}(X)=\sqrt{\frac{\epsilon}{\pi}}\tau\left(Xe^{-\epsilon D^{2}}
\right)$$
for $D\in {\cm}_0$ and $X\in T_D({\cm}_0)={\cn}_{sa}.$ Then the integral
of $\alpha^{\epsilon}$ is independent of the piecewise $C^1$ path in ${\cm}_0$
and if $\{D_t\}_{t\in [a,b]}$ is any piecewise $C^1$ path in ${\cm}_0$ then
$$sf(D_a,D_b)=
\sqrt{\frac{\epsilon}{\pi}}\int_a^b \tau\left(D_t'e^{-\epsilon D_t^2}\right)dt
+\frac{1}{2}\eta_{\epsilon}(D_b) - \frac{1}{2}\eta_{\epsilon}(D_a)
+\frac{1}{2}\tau([ker(D_b)]-[ker(D_a)]).$$
\end{cor}

\section{\bf FINITELY SUMMABLE MODULES REVISITED}

We show in this section how an application of Corollary 7.10 and the
Laplace Transform combine to give a ``best possible'' version
of Theorem 2.17 of \cite{CP1}. The importance of this result is that we
can use it in the case of $(1,\infty)$-Breuer-Fredholm modules to obtain
Connes' Dixmier-trace formula for the index in a much wider setting
\cite{CPSu}. We do this by computing the limit as $p\to 1$ in the best 
possible formula.

\begin{lemma}
If $n>0$ (not necessarily an integer) and $D$ is an unbounded
self-adjoint operator affiliated with $\cn$,
such that $(1+D^2)^{-n}$ is trace-class,
then the integral $\int_0^\infty e^{-t(1+D^2)}t^{(n-1)}dt$
converges in both trace-norm and operator-norm to $\Gamma(n)(1+D^2)^{-n}$.
\end{lemma}

\begin{proof}[\bf Proof]
Since $(1+D^2)^{-1}$ is in ${\mathcal L}^n$,
$\mu_s((1+D^2)^{-1})$ is $O(\frac{1}{s^{1/n}})$ by Lemma B.2
and so  $(1+D^2)^{-1}\in Li_0$.
Now for $t>0$, the map
$t\mapsto \frac{1}{t}(1+D^2)^{-1}$ is clearly continuous
in $Li_0$ and hence $t\mapsto e^{-t(1+D^2)}\in{\mathcal L}^1 $
is continuous for $t>0$ by Corollary B.9.
Thus the integrand is a continuous  ${\mathcal L}^1$-
valued function of $t$.
Now
\begin{eqnarray}
\int_0^\infty||e^{-t(1+D^2)}t^{(n-1)}||_1dt
&=& \int_0^\infty\tau(e^{-t(1+D^2)}t^{(n-1)})dt\nonumber\\
&=&\int_0^\infty t^{(n-1)}\int_0^\infty e^{-\frac{t}{\mu_s((1+D^2)^{-1})}}ds dt
\nonumber
\end{eqnarray}
by \cite {FK} Corollary 2.8.
As the integrand is positive we may interchange the order of integration
by Tonelli's theorem to obtain
$$\int_0^\infty(\int_0^\infty t^{(n-1)}e^{-\frac{t}{\mu_s((1+D^2)^{-1})}}dt) ds
=\int_0^\infty\Gamma(n)\left(\frac{1}{\mu_s((1+D^2)^{-1}))}\right)^{-n}ds$$
by the Laplace Transform.
Now, this equals
$$\int_0^\infty\Gamma(n)\mu_s((1+D^2)^{-1})^nds=\Gamma(n)\tau((1+D^2)^{-n})$$
by \cite{FK} Corollary 2.8. Thus, the integral converges in $\mathcal L^1$-
norm.

Similarly the integrand of the statement of
the Lemma is operator-norm continuous
and
\begin{eqnarray}
\int_0^\infty||e^{-t(1+D^2)}t^{(n-1)}||dt
&\leq&\int_0^\infty t^{(n-1)}e^{-\frac{t}{||(1+D^2)^{-1}||}}dt\nonumber\\
&=&\Gamma(n)||(1+D^2)^{-1}||^n\nonumber
\end{eqnarray}
so the integral also converges in operator-norm. Clearly this
limit operator is non-negative.

Let $\{E_\lambda\}_{\lambda\in[0,1]}$
be the spectral resolution of $(1+D^2)^{-1}$.
Then for $\xi\in\mathcal H$ we have by the Spectral Theorem,
\begin{eqnarray}
<\Gamma(n)(1+D^2)^{-n}\xi,\xi>
&=&\int_0^1\Gamma(n)\lambda^n d<E_\lambda\xi,\xi>\nonumber\\
&=&\int_0^1(\int_0^\infty t^{n-1} e^{-t/\lambda} dt) d<E_\lambda\xi,\xi>
\nonumber\\
&=&\int_0^\infty t^{n-1}(\int_0^1 e^{-t/\lambda} d<E_\lambda\xi,\xi>)dt
\text{\ \ \    by  Tonelli's Theorem} \nonumber\\
&=&\int_0^\infty t^{n-1} <e^{-t(1+D^2)}\xi,\xi>dt
=<(\int_0^\infty t^{n-1}e^{-t(1+D^2)}dt)\xi,\xi>.\nonumber
\end{eqnarray}
Hence, $$\Gamma(n)(1+D^2)^{-n}=\int_0^\infty t^{n-1}e^{-t(1+D^2)}dt$$
where the integral converges in both norms as claimed.
\end{proof}

\begin{lemma}
If $(1+D_0^2)^{-1}$ is in ${\mathcal L}^n$
and $\{D_t\}$ is a piecewise
$C^1$ path in ${\mathcal M}_0=D_0+{\cn}_{sa}$
then
$$\int_0^\infty
\epsilon^{n-1}e^{-\epsilon}\int_0^1\tau(D_t'
e^{-\epsilon D_t^2)})dtd\epsilon$$
converges absolutely.
\end{lemma}

\begin{proof}[\bf Proof]
\begin{eqnarray}
\int_0^\infty
\epsilon^{n-1}e^{-\epsilon}\int_0^1|\tau(D_t'e^{-\epsilon D_t^2})|dtd\epsilon
&\leq& \sup_{t}||D_t'||\int_0^\infty
\epsilon^{n-1}e^{-\epsilon}\int_0^1
||e^{-\epsilon D_t^2}||_1dtd\epsilon\nonumber\\
&=&C\int_0^\infty
\epsilon^{n-1}\int_0^1 \tau(e^{-\epsilon(1+ D_t^2)}) dtd\epsilon\nonumber\\
&=&C\int_0^1\int_0^\infty\tau(\epsilon^{n-1}e^{-\epsilon(1+ D_t^2)})
d\epsilon dt\text{\;\;\;(Tonelli)}\nonumber\\
&=&C\int_0^1\tau(\Gamma(n)(1+ D_t^2)^{-n})dt \text{\;\;\;(previous Lemma)}
\nonumber\\
&\leq& C\Gamma(n) \sup_tf(||D_t-D_0||)^n\int_0^1 \tau((1+ D_0^2)^{-n})dt
\nonumber
\end{eqnarray}

\noindent by Lemma 6 and Corollary 4 of appendix B of \cite{CP1}.
\end{proof}

\begin{rems*}
For unbounded $p$-summable modules we can now prove a ``best possible''
result, at least when the endpoints of the path are unitarily
equivalent. The theorem below is optimal in two ways: first, the exponent
$p/2$ is the minimum for which the formula makes sense, and second,
we need no assumptions about the integrality of $p$ or $p/2$. This
result was conjectured in Appendix C of \cite{CP1}.

By similar methods, we can also derive an improved version of
Theorem 2.16 of \cite{CP1} when the
endpoints are not unitarily equivalent. However, the exponent we need in this
case, $(p+1)/2$, is not optimal. The reason for the extra $1/2$ in the
exponent is that when we apply the Laplace Transform trick to the (truncated
eta) correction terms, we need an exponent $n$ which makes
$D(1+D^2)^{-n}$ trace-class. When $D$ is $p$-summable, the minimum such $n$
is $(p+1)/2$. We omit this result, leaving the details to the interested reader.
\end{rems*}

\begin{thm}[{\it cf} Theorem 2.16 of \cite{CP1}]
Let $({\mathcal N}, D_0)$ be  an odd unbounded $p$-summable\\
Breuer-Fredholm module (for {\bf C})
and let ${\cm}_0=D_0+{\cn}_{sa}.$ Then for $D\in {\cm}_0$, $X\in T_D({\cm}_0)=
{\cn}_{sa},$
$$X\mapsto \frac{1}{\tilde{C}_{p/2}}\tau(X(1+D^2)^{-p/2})$$
is an exact one-form on ${\cm}_0.$ Moreover, if
$\{D_t\}_{t\in[a,b]}$ is a piecewise
$C^1$  path in ${\mathcal M}_0$ with $D_a$ and $D_b$ unitarily equivalent
then
$$sf(D_a,D_b)=
\frac{1}{\tilde C_{p/2}}
\int_a^b\tau(\frac{d}{dt}(D_t)(1+D^2_t)^{-p/2})dt,$$

\noindent where $\tilde{C}_{p/2}=\int_{-\infty}^\infty(1+x^2)^{-p/2}dx.$
\end{thm}

\begin{proof}[\bf Proof]
By Corollary 7.10, we have for each $\epsilon > 0$:
$$sf(D_a,D_b)=
\sqrt{\frac{\epsilon}{\pi}}\int_a^b \tau(D_t'e^{-\epsilon D_t^2})dt.$$
Letting $n=p/2$, the Laplace Transform
gives
$$1=\frac{1}{\Gamma(n-\frac{1}{2})}
\int_0^\infty\epsilon^{n-3/2}e^{-\epsilon}d\epsilon.$$
Thus combining these expressions yields
\begin{eqnarray}
sf(D_a,D_b)&=&\frac{1}{\Gamma(n-\frac{1}{2})}
\int_0^\infty\epsilon^{n-3/2}e^{-\epsilon}\sqrt{\frac{\epsilon}{\pi}}
\int_a^b \tau(D_t'e^{-\epsilon D_t^2})dtd\epsilon\nonumber\\
&=&\frac{1}{\Gamma(n-\frac{1}{2})\sqrt{\pi}}\int_a^b \int_0^\infty
\epsilon^{n-1}\tau(D_t'e^{-\epsilon(1+ D_t^2)})d\epsilon dt\nonumber\\
&=&\frac{1}{\Gamma(n-\frac{1}{2})\sqrt{\pi}}\int_a^b \int_0^\infty
\tau[D_t'\epsilon^{n-1}e^{-\epsilon(1+ D_t^2)}]d\epsilon dt\nonumber
\end{eqnarray}
using Lemma 9.2 and Fubini. Using Lemma 9.1
this double integral becomes
\begin{eqnarray}
&=&\frac{1}{\Gamma(n-\frac{1}{2})\sqrt{\pi}}\int_a^b\tau[D_t'
\int_0^\infty\epsilon^{n-1}e^{-\epsilon(1+D_t^2)}d\epsilon]dt\nonumber\\
&=&\frac{\Gamma(n)}{\Gamma(n-\frac{1}{2})\sqrt{\pi}}\int_a^b
\tau(D_t'(1+D_t^2)^{-n})dt.\nonumber
\end{eqnarray}
Finally, the normalization constant is the beta function:
$$B(n-1/2,1/2)=\frac{\Gamma(n-\frac{1}{2})\sqrt{\pi}}{\Gamma(n)}=
\frac{\Gamma(n-\frac{1}{2})\Gamma(1/2)}{\Gamma(n)}=\int_0^1t^{(n-3/2)}
(1-t)^{-1/2}dt$$
by \cite{Ru} Theorem 8.20. By the change of variables $t=1/(1+x^2)$ this
is $\int_{-\infty}^\infty(1+x^2)^{-n}dx$ which is the
constant $\tilde C_n.$

It follows that $sf(D_a,D_b)$ is given by the integral of our one-form
when the endpoints are unitarily equivalent. Thus, the integral of our
one-form around any closed path is $0$, and so the integral of the one-form
is independent of path, in general. Thus, the one-form is exact by Lemma 7.3.
\end{proof}

\begin{cor}[{\it cf} Theorem 2.17 of \cite{CP1}]
Let $({\cn},D_0)$ be an odd p-summable Breuer-\\Fredholm module
for the unital Banach $*$-algebra ${\ca}$, and let $P=\chi_{[0,\infty)}(D_0).$
Then for each $u\in U({\ca})$ with $u(domD_0)\subseteq dom(D_0)$ and $[D_0,u]$
bounded,
$PuP$ is a Breuer-Fredholm operator in $P{\cn}P$ and if $\{D_t^u\}$ is any
piecewise $C^1$ path in ${\cm}_0=D_0+{\cn}_{sa}$ from $D_0$ to $uD_0u^*$
(e.g., the linear path lies in ${\cm}_0$), then:
$$ind(PuP)=sf(\{D_t^u\})=\frac{1}{\tilde{C}_{p/2}}\int_0^1\tau\left(
\frac{d}{dt}(D_t^u)(1+(D_t^u)^2)^{-p/2}\right)dt,$$

\noindent the integral of the exact one-form, $\frac{1}{\tilde{C}_{p/2}}
\tau\left(X(1+D^2)^{-p/2}\right)$ along the path $\{D_t^u\}.$
\end{cor}

\section{\bf SPECTRAL FLOW AND THE JLO COCYCLE}

In this section we generalise the main theorem of \cite{G} 
relating the spectral flow formula and the JLO cocycle to 
the case of spectral flow in semifinite von Neumann algebras.
We adopt a more concrete functional analytic method 
than in \cite{G} to avoid having to introduce more background 
material and also because there are additional subtleties in the type
$II$ setting. Throughout this section, $(\mathcal N,D)$ is an odd
unbounded $\theta-$summable Breuer-Fredholm module for a
Banach $*-$algebra $\mathcal A$ (contained in $\mathcal N$) and $u\in\mathcal A$
is a unitary operator leaving $dom(D)$ invariant and satisfying $[D,u]$ is 
bounded. There are three steps which we divide into three subsections.

\subsection{The graded space}
We form a new graded Hilbert space
${\mathcal K}=\IC^2\otimes\IC^2\otimes\mathcal H$.
Introduce the Clifford algebra on $\IC^2$ with generators:
$$\sigma_1= \left(\begin{array}{cc}
                   0 & 1 \\
                   1 & 0
                   \end{array} \right),\ \  \sigma_2= \left(\begin{array}{cc}
                   0 & -i \\
                   i & 0
                   \end{array} \right), \ \ \sigma_3=
 \left(\begin{array}{cc} 1 & 0 \\
                         0 & -1
\end{array} \right).$$

Let $\sigma_0$ denote the $2\times 2$ identity matrix.
Then the grading on ${\mathcal K}$ is given by
$$\Gamma = \sigma_2\otimes\sigma_3\otimes I = \sigma_2\otimes
\left(\begin{array}{cc}
                   I & 0 \\
                   0 & -I
                   \end{array} \right).$$
where $I$ is the identity operator in $\mathcal N$.
Let $u\in \mathcal A$ be a unitary and
introduce the following operators on $\mathcal K$, all of which commute
with $\Gamma$, by
$$ D_0 = \sigma_2\otimes\sigma_0\otimes D = \sigma_2\otimes
\left(\begin{array}{cc}
                   D & 0 \\
                   0 & D
                   \end{array} \right),\ \ 
q=\sigma_3\otimes\left(\begin{array}{cc}
                   0 & -iu^{-1} \\
                   iu & 0
                   \end{array} \right),\ \
$$
$$
D_r= (1-r)D_0 - rqD_0q, \ \   
D_{r,s} = D_r + sq, \ \  r\in [0,1], s\in [0,\infty).$$
Notice that if we define $D_r\equiv D_{r,0}$ then,
$$ D_{r} = \sigma_2\otimes \left(\begin{array}{cc}
                   D+ru^{-1}[D,u] & 0 \\
                   0 & D+ru[D,u^{-1}]
                   \end{array} \right)
=D_0 +r\sigma_2\otimes\left(\begin{array}{cc}
                   u^{-1}[D,u] & 0 \\
                   0 & u[D,u^{-1}]
                   \end{array} \right).$$
So that
$$\dot D_r = \sigma_2\otimes \left(\begin{array}{cc}
                   u^{-1}[D,u] & 0 \\
                   0 & u[D,u^{-1}]
                   \end{array} \right).$$

Introduce a graded trace:
$Str(a) =\frac{1}{2\sqrt\pi}
 \tr(\Gamma a)$ for $a$ trace class and so for example
$$Str(\dot D_r e^{-D_r^2}) 
= \frac{1}{\sqrt\pi}\tr\{u^{-1}[D,u]e^{-(D+ru^{-1}[D,u])^2}- u[D,u^{-1}] 
e^{-(D+ru[D,u^{-1}])^2}\}.$$
Next we calculate
$$D_{r,s}^2= D_r^2+ s(1-2r)\sigma_1\otimes\left(\begin{array}{cc}
                   0 & [D,u^{-1}] \\
                   -[D,u] & 0
                   \end{array} \right) +s^2$$
which depends on the relation
$$D_0 q + qD_0 =\sigma_1\otimes\left(\begin{array}{cc}
                   0 & [D,u^{-1}] \\
                   -[D,u] & 0
                   \end{array} \right) $$
The preceding relation explains in part the
reason for introducing $\mathcal K$ and the grading: it converts
commutators to anticommutators, for example:
$$\int_0^1 Str(q e^{-tD_0^2}(D_0q+ qD_0) 
e^{-(1-t)D_0^2} )dt$$
$$=\frac{1}{2\sqrt\pi}
\int_0^1\tr\{\sigma_2\sigma_3\sigma_1\otimes \left(\begin{array}{cc}
                   iu^{-1}e^{-tD^2}[D,u]e^{-(1-t)D^2} & 0 \\
                   0 & -iue^{-tD^2}[D,u^{-1}]e^{-(1-t)D^2}
                   \end{array} \right)\}dt $$
$$= -\frac{1}{\sqrt\pi}\int_0^1\tr( u^{-1}e^{-tD^2}[D,u]e^{-(1-t)D^2} - 
ue^{-tD^2}[D,u^{-1}]e^{-(1-t)D^2}).$$

\subsection{Changing the path of integration}

We now want to see what our spectral flow formula looks like on $\mathcal K$.
\begin{lemma}
$$\int_0^1Str(\dot D_r e^{-D_r^2}) dr= 2  sf\{D,u^{-1}Du\}.$$
\end{lemma}

\begin{proof}

$$D_r^2= \sigma_0\otimes \left(\begin{array}{cc}
                   (D+ru^{-1}[D,u])^2 & 0 \\
                   0 & (D+ru[D,u^{-1}])^2)
                   \end{array} \right).$$
so that 
$$\int_0^1Str(\dot D_r e^{-D_r^2}) dr$$
\begin{eqnarray*}
&=& \frac{1}{2\sqrt\pi}
\int_0^1\tr\{\Gamma \sigma_2\otimes \left(\begin{array}{cc}
                   u^{-1}[D,u]e^{- (D+ru^{-1}[D,u])^2} & 0 \\
                   0 & u[D,u^{-1}]e^{-(D+ru[D,u^{-1}])^2}
                   \end{array} \right)\}dr\\
&=& \frac{1}{\sqrt\pi}\int_0^1\tr\{\left(\begin{array}{cc}
                   u^{-1}[D,u]e^{- (D+ru^{-1}[D,u])^2} & 0 \\
                   0 & -u[D,u^{-1}]e^{-(D+ru[D,u^{-1}])^2}
                   \end{array} \right)\}dr\\
&=&\frac{1}{\sqrt\pi}\int_0^1\tr\{ u^{-1}[D,u]e^{- (D+ru^{-1}[D,u])^2}\}dr 
                   -\frac{1}{\sqrt\pi}
\int_0^1\tr\{u[D,u^{-1}]e^{-(D+ru[D,u^{-1}])^2}\}dr
\end{eqnarray*}
Now the first term is $sf\{D,u^{-1}Du\}$
and the second is  $sf\{D,uDu^{-1}\}$.
As 
$$sf\{D,u^{-1}Du\}= sf\{uDu^{-1}, D\}=-sf\{D,u^{-1}Du\}$$
so that, as required,
$$\int_0^1Str(\dot D_r e^{-D_r^2}) dr= 2 sf\{D,u^{-1}Du\}.$$
\end{proof}

The main idea of \cite{G} is to change the path of integration
used to compute the spectral flow. This is achieved with the help
of the next result.

\begin{lemma} Consider the affine space $\Phi$
of perturbations of $D_0$
given by 
$$\{D_0+X \ \|\  X\in M_2\otimes M_2\otimes\mathcal N\;\;
\mbox{is self-adjoint and even} \;( i.e.,
\Gamma X=X\Gamma)\}$$
Then 
the map
$$X\to Str(Xe^{-(D_0+Y)^2})$$ is an exact one-form.
\end{lemma}

\begin{proof}
Closedness of this one-form is proved by using the methods of Appendix C 
and then exactness follows from a Poincar\'e lemma for 
the affine space $\Phi$ (see the remark after Theorem 7.8).
\end{proof}

Now consider the rectangle $R$ in $\IR^2$ given by
$0\leq r\leq 1,0\leq s\leq s_0$ for some $s_0$. 
Using the previous exactness result we conclude that the integral:
$$ \int_{\partial R} Str(\frac {dD_{r,s_0}}{dr}e^{-D_{r,s_0}^2})$$
around the boundary $\partial R$ is zero. So we can replace 
our original integral by a sum of three integrals and 
calculate the contribution of each to the spectral flow.

\begin{lemma}
$$\lim_{s_0\to \infty} \int_0^1 Str(\frac {dD_{r,s_0}}{dr}e^{-D_{r,s_0}^2})dr
=0$$
\end{lemma}
\begin{proof}
First we observe that
$$ Str(\frac {dD_{r,s_0}}{dr}e^{-D_{r,s_0}^2})
=Str (\dot D_re^{-D_{r,s_0}^2}).$$
Next notice that $D_{1/2}=\frac{1}{2}(D_0-qD_0q)$
so that $D_{1/2}$ anticommutes with $q$. Then
$(D_{1/2} +sq)^2=D^2_{1/2} +s^2$ so that
$$Str (\dot D_{1/2} e^{-(D_{1/2} +sq)^2})
=e^{-s^2}Str (\dot D_{1/2}e^{-D^2_{1/2}})$$
which decays exponentially to zero as $s\to\infty$.
Let $A_r= D_r-D_{1/2}=(\frac{1}{2} -r)(D_0+qD_0q).$ Since
$$D_0+qD_0q = -\dot D_r = \sigma_2\otimes \left(\begin{array}{cc}
                   -u^{-1}[D,u] & 0 \\
                   0 & -u[D,u^{-1}]
                   \end{array} \right),$$
we see that $A_r$ is bounded
by a constant independent of $r\in [0,1]$.
Using \cite{CP1} Corollary 8 Appendix B we know there are constants
$C,C'$ depending only on $||D_0+qD_0q||$ such that
$$\tr(e^{-(D_r+sq)^2})\leq C\tr(e^{-C'(D_{1/2}+sq)^2}.$$
Thus
$$|Str (\dot D_re^{-D_{r,s_0}^2})|\leq ||D_0+qD_0q||
\tr(e^{-(D_r+s_0q)^2})\leq C||D_0+qD_0q||\tr(e^{-C'(D_{1/2}+sq)^2})
$$
so that
$$ Str(\frac {dD_{r,s_0}}{dr}e^{-D_{r,s_0}^2})$$ decays
exponentially to zero uniformly in $r$ as $s_0\to \infty$
proving the result.
\end{proof}

\begin{lemma}
$$\int_0^\infty Str (\frac{dD_{1,s}}{ds}e^{-D_{1,s}^2})ds=
-\int_0^\infty Str (\frac {dD_{0,s}}{ds}e^{-D_{0,s}^2})ds$$
\end{lemma}
\begin{proof}
Note first that 
$$D_{1,s} = -qD_0 q +sq = -q(D_0 -sq)q$$
so that
$$Str (\frac {dD_{1,s}}{ds}e^{-D_{1,s}^2})=Str(qe^{-(D_0-sq)^2})$$
using invariance under conjugation by $q$.
Now let $\rho =\sigma_2\otimes \sigma_0\otimes I$.
Then using $\rho q\rho = -q$ and $\rho \Gamma = \Gamma \rho$
we have
$$2\sqrt\pi Str (\frac {dD_{1,s}}{ds}e^{-D_{1,s}^2})
=\tr(\Gamma qe^{-(D_0-sq)^2})
=\tr(\rho^2\Gamma  qe^{-(D_0-sq)^2})$$
$$
=\tr(\rho\Gamma  qe^{-(D_0-sq)^2})\rho)
=\tr(\rho\Gamma  q\rho e^{-(D_0+sq)^2}))
$$
$$
=-\tr(\rho\Gamma  \rho q e^{-(D_0+sq)^2})) 
=2\sqrt{\pi}Str (\frac {dD_{0,s}}{ds}e^{-D_{0,s}^2}) 
\eqno(*)$$
\end{proof}
Combining the above results yields the key observation of this subsection:
\begin{cor}
$$sf\{D,u^{-1}Du\}=
\int_0^\infty Str (\frac {dD_{0,s}}{ds}e^{-D_{0,s}^2})ds \eqno (**)$$
\end{cor}

\subsection{The Duhamel argument}

Given the last corollary
the essential observation in \cite{G} is to use the Duhamel Principle
to evaluate 
\begin{eqnarray*}
\hspace{.8in} 
\int_0^\infty Str (\frac {dD_{0,s}}{ds}e^{-D_{0,s}^2})ds
&=&\frac{1}{2\sqrt\pi}
\int_0^\infty\tr(\Gamma q e^{-({D_0}^2+s(D_0 q+qD_0)+s^2)})ds\\
&=&\frac{1}{2\sqrt\pi}
\int_0^\infty e^{-s^2}\tr(\Gamma q e^{-({D_0}^2+s(D_0 q+qD_0))})ds.\hspace{.3in}
(***)
\end{eqnarray*}
To use Duhamel (see \cite{Y}, p. 438) we write $[D_0, q]_+= D_0 q+qD_0$ and
then, checking (by use of the Spectral Theorem) that the formal derivatives
actually exist in operator norm,
\begin{eqnarray*}
e^{-({D_0}^2+s[D_0, q]_+)}-e^{-{D_0}^2} &=& 
- \int_0^1\frac{d}{dt}\left(e^{-tD_0^2}e^{-(1-t)(D_0^2+s[D_0,q]_+)}\right)dt\\
&=& -s\int_0^1e^{-tD_0^2}[D_0,q]_+
e^{-(1-t)(D_0^2+s[D_0,q]_+)}dt
\end{eqnarray*}
\begin{eqnarray*}
&=& -s\int_0^1e^{-tD_0^2}[D_0,q]_+\left(e^{-(1-t)D_0^2}
-\int_0^{1-t}\frac{d}{dx}\left(e^{-xD_0^2}e^{-(1-t-x)(D_0^2+s[D_0,q]_+)}\right)
dx\right)dt\\
&=&\ldots\\
&=& \sum_{k=1}^\infty(-s)^k\int_{\Delta_k}
e^{-t_k{D_0}^2}[D_0, q]_+e^{-t_{k-1}{D_0}^2}[D_0, q]_+
\ldots[D_0, q]_+e^{-t_0{D_0}^2} dt_k\ldots dt_1
\end{eqnarray*}

Where $\Delta_n$ is the standard $n$-simplex in $\IR^{n+1}$ given by
$$\{(t_0,t_1,\ldots, t_n)\vert t_j\in [0,1],\  \sum_0^n t_j =1\}.$$
Since $vol(\Delta_k) = 1/k!$, and the integrand over each $\Delta_k$
can be estimated in trace-norm (using the general H\"{o}lder
inequality) by $||[D_0,q]_+||^k tr(e^{-D_0^2})$, it is not hard to see that 
the series converges in trace-norm. Thus, we can substitute the resulting
formula for $e^{-({D_0}^2+s[D_0, q]_+)}$ into the integral formula $(***)$ for
$sf\{D,u^{-1}Du\}$. After a little manipulation we must evaluate the following.

\begin{lemma}
With the above hypotheses and notation, we have two cases.\\
If $n=2k+1$ is odd, then:
$$\int_{\Delta_n}
Str\{qe^{-t_0{D_0}^2}[D_0,q]_+
e^{-t_1{D_0}^2}[D_0,q]_+
\ldots [D_0,q]_+e^{-t_{n}D_0^2}\}
dt_1\ldots dt_{n}$$
$$=-\frac{1}{\sqrt\pi}\int_{\Delta_n}(-1)^k\tr\{u^{-1} e^{-t_0D^2}[D,u]
e^{-t_1D^2}[D,u^{-1}]\ldots [D,u]e^{-t_{n}D^2}\}dt_1\ldots dt_{n}
$$
$$+\frac{1}{\sqrt\pi}\int_{\Delta_n}
(-1)^k\tr\{u e^{-t_0D^2}[D,u^{-1}]e^{-t_1D^2}[D,u]\ldots
[D,u^{-1}]e^{-t_{n}D^2}\}dt_1\ldots dt_{n}.$$
If $n$ is even, then
$$\int_{\Delta_n}
Str\{qe^{-t_0{D_0}^2}[D_0,q]_+
e^{-t_1{D_0}^2}[D_0,q]_+
\ldots [D_0,q]_+e^{-t_{n}D_0^2}\}
dt_1\ldots dt_{n}=0.$$

\end{lemma}
\begin{proof}
For $n=2k+1$ odd, this is a straightforward calculation. For $n$ even,
we observe that the $\Gamma q$ contributes a first tensor factor of $i\sigma_1$;
each of the $n$ copies of $[D_0,q]_+$ contributes a first tensor factor of
$\sigma_1$; and each of $n+1$ exponential terms contributes a first tensor 
factor of $\sigma_0$. This yields a first tensor factor of $i\sigma_1$ and
so the trace is zero.
\end{proof}

\begin{lemma} With the above hypotheses and notation,
$$\sum_{k=0}^\infty k!\int_{\Delta_n}(-1)^k\tr\{u^{-1}
e^{-t_0D^2}[D,u]e^{-t_1D^2}[D,u^{-1}]
e^{-t_2D^2}\ldots [D,u]e^{-t_{2k+1}D^2}\}dt_1\ldots dt_{2k+1}$$
$$=-\sum_{k=0}^\infty k!\int_{\Delta_n}
(-1)^k\tr\{u e^{-t_0D^2}[D,u^{-1}]e^{-t_1D^2}[D,u]e^{-t_2D^2}\ldots
[D,u^{-1}]e^{-t_{2k+1}D^2}\}dt_1\ldots dt_{2k+1}.$$
\end{lemma}

\begin{proof}
This result is a simple consequence of the fact that the JLO formula
defines a cocycle in the $(b,B)$-bicomplex (Theorem IV.21 in \cite{Co4}:
see also \cite{JLO,GZ}).
It is well known to experts but we could not find a good exposition in the 
literature and so for completeness we indicate the proof.
The JLO cocycle is the sequence of multilinear functionals 
$(\phi_n)=\phi_{JLO}$
on $({\mathcal A}^{n+1})$, $n=0,1, 2,\ldots$
where for $n=2k+1$ (this is the only case we need):
$$\phi_{2k+1}(a^0,a^1,\ldots,a^{2k+1})
$$
$$=\sqrt{2i}\int_{\Delta_n}
(-1)^k\tr\{a^0e^{-t_0D^2}[D,a^1]e^{-t_1D^2}[D,a^2]e^{-t_2D^2}\ldots
[D,a^{2k+1}]e^{-t_{2k+1}D^2}\}dt_1\ldots dt_{2k+1}.$$
Introduce the sequence for $k=1,2,\ldots$
$$(1,u^{-1},u,\ldots,u^{-1},u)_{2k+2}\in{\mathcal
A}^{2k+3}$$
where $1$ is the identity of ${\mathcal A}$
and the subscript indicates the number of terms to the right of the first 
element
in each term of the sequence.
We calculate, using the standard formula for the Hochschild coboundary
operator $b$ \cite{Co4}
$$b\phi_{2k+1} (1,u^{-1},u,\ldots,u^{-1},u)=
\phi_{2k+1} (u^{-1},u,\ldots,u^{-1},u)+\phi_{2k+1}
(u,u^{-1},u,\ldots,u^{-1})
.$$
For the operator $B$ of the $(b,B)$-bicomplex we use the formula
$B=AB_0$ \cite{Co4} where $A$ is the antisymmetrisation operator and
$B_0$ acts on an $n+2$ linear functional $\psi$ by
$$B_0\psi(a^0,a^1,\ldots,a^{n})=\psi(1,a^0,a^1,\ldots,a^{n})-
(-1)^{n+1}\psi(a^0,a^1,\ldots,a^{n},1).$$

Thus, 
\begin{eqnarray*}
B_0\phi_{2k+1} (1,u^{-1},u,\ldots,u^{-1},u)&=&\\
\phi_{2k+1} (1,1,u^{-1},u,\ldots,u^{-1},u)&-&
\phi_{2k+1} (1,u^{-1},u,\ldots,u^{-1},u,1)=0
\end{eqnarray*}
because the commutator with $D$ in the JLO formula kills
all terms with two copies of the identity operator. 
It follows that $$B\phi_{2k+1}(1,u^{-1},u,\ldots,u^{-1},u)=0,$$
and since in the  $(b,B)$-bicomplex we have that $(b+B)\phi_{2k+1}=0$, we get 
\begin{eqnarray*}
0&=&b\phi_{2k+1} (1,u^{-1},u,\ldots,u^{-1},u)\\
&=&\phi_{2k+1} (u^{-1},u,\ldots,u^{-1},u)+\phi_{2k+1}
(u,u^{-1},u,\ldots,u^{-1})
\end{eqnarray*}
which proves the statement of the lemma.
\end{proof}

Hence we can use the preceding two lemmas to obtain our final formula:
$$\int_0^\infty Str (\frac {dD_{0,s}}{ds}e^{-D_{0,s}^2})ds=
$$
$$-\int_0^\infty\int_{\Delta_{n}}\sum_{n\;odd}^\infty s^{n}e^{-s^2}
Str\{ q
e^{-t_{n}{D_0}^2}[D_0, q]_+
e^{-t_{n-1}{D_0}^2}[D_0, q]_+
\ldots[D_0, q]_+e^{-t_0{D_0}^2}\} dt_{n} \ldots dt_1ds
$$
$$=-\sum_{k=0}^\infty \frac{k!}{2} \int_{\Delta_{2k+1}}
Str\{ q
e^{-t_{2k+1}{D_0}^2}[D_0, q]_+
e^{-t_{2k}{D_0}^2}[D_0, q]_+
\ldots[D_0, q]_+e^{-t_0{D_0}^2}\} dt_{2k+1} \ldots dt_1
$$
$$=\frac{1}{\sqrt\pi}\sum_{k=0}^\infty(-1)^k {k!} \int_{\Delta_{2k+1}}
\tr\{u^{-1} e^{-t_0D^2}[D,u]e^{-t_1D^2}[D,u^{-1}]
\ldots [D,u]e^{-t_{2k+1}D^2}\} dt_{2k+1} \ldots dt_1.
$$

Now we have, by Corollary 10.5  the connection between cyclic cohomology
in the form of the JLO cocycle and spectral flow.

\begin{thm} Let $(\mathcal N,D)$ be an odd
unbounded $\theta-$summable Breuer-Fredholm module for the
Banach $*-$algebra $\mathcal A$ and let $u\in\mathcal A$
be a unitary operator leaving $dom(D)$ invariant and satisfying $[D,u]$ is 
bounded. Then,
$$sf\{D,u^{-1}Du\}$$
$$=\frac{1}{\sqrt\pi}\sum_{k=0}^\infty (-1)^k {k!} \int_{\Delta_{2k+1}}
\tr\{u^{-1} e^{-t_0D^2}[D,u]e^{-t_1D^2}[D,u^{-1}]
\ldots [D,u]e^{-t_{2k+1}D^2}\} dt_{2k+1} \ldots dt_1.
$$
\end{thm}

One may interpret this result as a Chern character 
formula pairing the $K_1({\mathcal A})$ class $[u]$
of $u\in\mathcal A$ with an explicit entire cyclic cocycle for $\mathcal A$
obtained via the semifinite $K$-homology class defined by $D$.
We will not pursue this point of view here mentioning only the fact
that the main interest in this result is that it allows
us to prove an analogue of the Connes-Moscovici local index theorem
\cite{CoMo}in the semi-finite case.

\appendix

\section{\bf OPERATOR IDEALS}

This appendix establishes some properties
of certain ideals of operators in
a von Neumann algebra $\mathcal N$  with a faithful, normal
semifinite trace $\tau$. A similar discussion is
contained in section 5 of \cite{Suk}. One can prove versions of the
results in this section for the Marcinkiewicz spaces defined in section
5 of \cite{Suk} (see also \cite{DD}).
Our point of view is slightly
different from and more na\"ive than his: we try to establish the
basic properties of these normed ideals with as little machinery as
possible in hopes of making the material more accessible to those as
na\"ive as ourselves. We claim no originality for the results themselves.
We first recall Definition 2.1.
\begin{defn}
 If $S\in\mathcal N$ the {\bf t-th generalized
singular value of S} for each real $t>0$ is given by
$$\mu_t(S)=\inf\{||SE||\ \vert \ E \text{ is a projection in }
{\mathcal N} \text { with } \tau(1-E)\leq t\}.$$
\end{defn}

For the basic properties of these singular values we refer to
\cite{FK}.

\begin{defn}
If $\mathcal I$ is a $*$-ideal in $\mathcal N$
which is complete in a norm $||\cdot||_{\mathcal I}$
then we will call $\mathcal I$ an {\bf invariant
operator ideal} if\\
(1) $||S||_{\mathcal I}\geq ||S||$ for all $S\in \mathcal I$,\\
(2) $||S^*||_{\mathcal I} = ||S||_{\mathcal I}$ for all $S\in \mathcal I$,\\
(3) $||ASB||_{\mathcal I}\leq ||A|| \:||S||_{\mathcal I}||B||$ for all
$S\in \mathcal I$, $A,B\in \mathcal N$.\\
Since $\mathcal I$ is an ideal in a von Neumann algebra, it follows
from I.1.6, Proposition 10 of \cite{Dix} that if $0\leq S \leq T$
and $T \in {\mathcal I}$, then $S \in {\mathcal I}$ and $||S||_{\mathcal I}
\leq ||T||_{\mathcal I}$. Much more is true, especially in the type $I$ case
but we shall not need it here, see \cite{GK}.
\end{defn}

\noindent {\bf Examples}.\\
(1) Let ${\mathcal I}=Li =\{T\in {\mathcal N}\ \vert\
\mu_s(T)=O(1/\log s)\}$. The norm on $Li$ is:
$$\|T\|_{Li}=\sup_{r>0}\left\{\frac{\int_0^r\mu_s(T) ds}{\int_0^r
(\log{(s+e)})^{-1}ds}\right\}.$$
We observe that $\|T\|_{Li}\geq ||T||$.\\
(2) We let ${\mathcal I}=Li_0 =\{T\in {\mathcal N}\ \vert\
\mu_s(T)=o(1/\log s)\}$ with the norm inherited from $Li$.
Using the estimate
$$r(\log(r+e))^{-1}\leq \int_0^r(\log(s+e))^{-1} ds\leq 3r(2\log(r+e))^{-1}$$
of Lemma A.4 below it follows that $Li_0$
 is the closure of the ideal $\mathcal F_N$
of `finite rank' operators in the $Li$ norm, see \cite{GK}.
Here
$$\mathcal F_N=\{T\in{\mathcal N}\ |\ T=ET \text{ for some projection }
E\in{\mathcal N} \text { with } \tau(E)<\infty\}.$$
(3) For $0< q\leq 1$, let ${\mathcal I}=Li^q= \{T\in {\mathcal N}\ \vert\
\mu_s(T)=O((1/\log s)^q)\}$
$$= \{T\in {\mathcal N}\ \vert\ |T|^{1/q}\in  Li\}.$$
We use the following norm on $Li^q$:
$$\|T\|_{Li^q}=(\||T|^{1/q}\|_{Li})^q.$$
We prove below that this is in fact a norm on $Li^q$ and that
H\"older's inequality is satisfied for these spaces.
There is an equivalent norm on $Li^q$ in which it is complete \cite{Suk},
given by:
$$\sup_{r>0}\left\{\frac{\int_0^r\mu_s(T) ds}{\int_0^r
(\log(s+e))^{-q}ds}\right\}.$$
We will not use this norm explicitly, however.

\begin{lemma}
(1) For each $q$ with $0<q\leq 1$, we have $||\cdot||_{Li^q}$ is
an invariant norm on $Li^q$.\\
(2)If $0<q,q'\leq 1$ where $q+q'\leq 1$ and
if $S\in Li^q$ and $T\in Li^{q'}$
then $ST\in Li^{q+q'}$ and
$$ ||ST||_{Li^{q+q'}}
\leq ||S||_{Li^{q}}||T||_{Li^{q'}}.$$
\end{lemma}

\begin{proof}[\bf Proof]
(1) The only nontrivial part is subadditivity. So, suppose $T,S \in Li^q$.
Then,
\begin{eqnarray}
||T+S||_{Li^q} &=& \left(||\:|T+S|^{1/q}||_{Li}\right)^q\nonumber\\
               &=& \sup_{t>0}\left[\frac{\int_0^t\mu_s(|T+S|^{1/q})ds}
                    {\int_0^t\frac{1}{\log(s+e)}ds}\right]^q\nonumber\\
      &\leq& \sup_{t>0}\left[\frac{\int_0^t[\mu_s(T)+\mu_s(S)]^{1/q}ds}
                    {\int_0^t\frac{1}{\log(s+e)}ds}\right]^q\text{\;\;(Theorem
                    4.4 part (iii), \cite{FK})}\nonumber\\
     &\leq& \sup_{t>0}\frac{\left[\int_0^t(\mu_s(T))^{1/q}ds\right]^{q}+
                    \left[\int_0^t(\mu_s(S))^{1/q}ds\right]^{q}}
                    {\left[\int_0^t\frac{1}{\log(s+e)}ds\right]^q}\nonumber\\
     &\leq& ||T||_{Li^q} + ||S||_{Li^q} .\nonumber
\end{eqnarray}
(2) This proof is very similar to part (1) except we cite Theorem 4.2 part
(iii) of \cite{FK} and then apply the usual H\"older inequality for the
interval $[0,t].$
\end{proof}

\begin{lemma} For $r\geq e$,
$$\int_e^r \frac{dx}{\log x} \leq \frac{3(r-e)}{2\log r}
$$
\end{lemma}

\begin{proof}[\bf Proof]
A straightforward calculus exercise.
\end{proof}

Note that we may reformulate this inequality as
$$\int_0^r \frac{dx}{\log(x+e)} \leq \frac{3r}{2\log(r+e)}.
$$

\begin{lemma}(1) For $S\in Li^q$, $0<q\leq 1$ let
$f_q(t)=\mu_t(S)(\log(t+e))^q$
then $$(2/3)^q||f_q||_\infty\leq ||S||_{Li^q}\leq ||f_q||_\infty.$$
(2) For $0<q\leq 1$ the embedding $Li\hookrightarrow Li^q$ is bounded
(by $3/2$) and ,in fact, for $0<q<1$,we have $Li \subseteq Li_0^q$.
\end{lemma}

\begin{proof}[\bf Proof]
(1) On the one hand,
\begin{eqnarray}
||S||_{Li^q} & = & \sup_{r> 0}\left[\frac{\int_0^r(\mu_t(S))^{1/q}dt}
                    {\int_0^r\frac{1}{\log(t+e)}dt}\right]^{q}\nonumber\\
             & = & \sup_{r> 0}\left[\frac{\int_0^r\frac{(f_q(t))^{1/q}}
                   {\log(t+e))} dt}
                    {\int_0^r\frac{1}{\log(t+e)}dt}\right]^{q}
                     \leq  ||f_q||_\infty.\nonumber
\end{eqnarray}
On the other hand by the previous lemma,
\begin{eqnarray}
||S||_{Li^q} &\geq& \sup_{r>0}\left[\frac{r\mu_r(S)^{1/q}}{1.5r/\log(r+e)}
\right]^q \nonumber\\
&=& (2/3)^q \sup_{r>0}\;[\mu_r(S)(\log(r+e))^q]\nonumber\\
&=& (2/3)^q||f_q||_\infty.\nonumber
\end{eqnarray}
(2) By part (1),
\begin{eqnarray}
||S||_{Li^q}\leq ||f_q||_{\infty} &=& \sup_{t>0}\;[\mu_t(S)\log(t+e)^q]\nonumber
\\
&\leq& \sup_{t>0}\;[\mu_t(S)\log(t+e)]\nonumber\\
&=& ||f_1||_{\infty} \leq (3/2)||S||_{Li}.\nonumber
\end{eqnarray}
If $0<q<1$, and $S \in Li$ so that $\mu_t(S) = O(1/\log{t})$ then clearly
$\mu_t(S) = o(1/(\log{t})^q)$ so that $S \in Li_0^q.$
\end{proof}

\section{\bf TRACE-CLASS CONTINUITY OF CERTAIN MAPS}

The following result is well-known in the case of type $I_\infty$ von
Neumann algebras. Fyodor Sukochev pointed out to us that, in general, it
follows from a result of \cite{FK}. Moreover, the
result has been extended by him and his co-authors
to many symmetric operator spaces including the spaces $Li^q$ ,
\cite{CDS}.

\begin{prop}.
Let $\{T_n\}_{n=1}^\infty$ and $T$ be positive trace class operators
in $\mathcal N$. If both $||T_n-T||\to 0$ and $||T_n||_1\to||T||_1$ then
$||T_n-T||_1\to 0$.
\end{prop}

\begin{proof}[\bf Proof]
It follows from $||T_n-T||\to 0$ that $T_n$ converges to $T$ in the
measure topology, and so the result follows from \cite{FK} Theorem 3.7.
\end{proof}

The following two results were first pointed out to
us by Chris Bose.
They are well-known as the containment of
$L^p$ in $L^{p,\infty}$.
We include the simple proofs of these known results for completeness,
see \cite{Suk}.

\begin{lemma} Let $f$ be a non-negative decreasing function
on $\real^{+}$ and suppose $f\in L^p$ for some $p>0$
then $f(x)\leq \frac{||f||_p}{x^{1/p}}$ for all $x>0$.
In other words $f$ is $O(\frac{1}{x^{1/p}})$ as $x\rightarrow\infty$.
\end{lemma}

\begin{proof}[\bf Proof]
As $f$ is decreasing,
$f(x)\chi_{[0,x]}\leq f$ for each $x>0$ and so
$$f(x)^px= \int_0^\infty f(x)^p\chi_{[0,x]} du\leq \int_0^\infty f(u)^p du
=||f||_p^p.$$
Hence
$f(x)\leq \frac{||f||_p}{x^{1/p}}$.
\end{proof}

\begin{cor}
If $f\in L^p(\real^{+})$ satisfies the hypotheses of the lemma then $f$
is $o(\frac{1}{x^{1/p}})$ as $x\rightarrow\infty$.
\end{cor}

\begin{proof}[\bf Proof]
 Given $\epsilon>0$, there is an $x_0$ such that
$\int_{x_0}^\infty f^p \leq \epsilon^p/2.$ Then for $x>x_0$:
$$f(x)^p(x-x_0)\leq \int_{x_0}^\infty f^p \leq \frac{\epsilon^p}{2}$$
or
$$f(x)^p\leq \frac{\epsilon^p}{2(x-x_0)}$$
implying $f(x)^p\leq \frac{\epsilon^p}{x}$ if $x\geq 2x_0$.
That is, $f(x)\leq \frac{\epsilon}{x^{1/p}}$ if $x\geq 2x_0$.
\end{proof}

\begin{lemma}
An operator $X$ is in ${Li}$ if and only if there is a $t_0>0$
such that
$$\tau(e^{-t|X|^{-1}})<\infty$$ for all $t>t_0$.
\end{lemma}

\begin{proof}[\bf Proof]
By definition,
$$X\in {Li}\iff \mu_x(|X|)=O(\frac{1}{\log x}).$$
Thus  $X\in {Li}$ is equivalent to the existence
of $t_0$ such that for all $x\geq x_0$, $\mu_x(|X|)\leq \frac{t_0}{\log x}$
which in turn is equivalent to
$e^{\frac{-t_0}{\mu_x(|X|)}}\leq 1/x$ for all $x\geq x_0$.
This is equivalent to
$$e^{-\frac{t}{\mu_x(|X|)}}\leq \frac{1}{x^{t/t_0}}$$
for all $x\geq x_0$, $t>t_0$.
Applying lemma 2.5 of \cite{FK} now gives
$X\in\ {Li}\iff$  there is a $t_0$ such that for all
$t>t_0$, $\mu_x(e^{-t|X|^{-1}})\leq
\frac{1}{x^{t/t_0}}$ for all $x\geq x_0$.
Hence the forward implication of the lemma follows
by lemma 2.7 of \cite{FK} as there is a $t_0$ such that for all $t>t_0$
$$\tau(e^{-t|X|^{-1}})= \int_0^\infty\mu_x(e^{-t|X|^{-1}}) dx
\leq C+ \int_{x_0}^\infty \frac{1}{x^{t/t_0}}dx<\infty.$$
The reverse implication follows because
$$\tau(e^{-t_1|X|^{-1}})=\int_0^\infty\mu_x(e^{-t_1|X|^{-1}}) dx<\infty$$
for some $t_1>0$ implies, by Corollary B.3 for $p=1$
that $\mu_x(e^{-t_1|X|^{-1}})$ is $o(\frac{1}{x})$.
So there is an $x_0\geq 0$ so that $x\geq x_0$ implies
 $\mu_x(e^{-t_1|X|^{-1}})\leq\frac{1}{x}.$
But then $\mu_x(|X|)= O(\frac{1}{\log x})$ which implies
that  $X\in {Li}$.
\end{proof}

\begin{rems*} By the first half of the proof, we can choose for $t_0$
any positive real number for which eventually in $x$,
$ \mu_x(|X|)\leq \frac{t_0}{\log x}$. So, if $X$ is in $Li$ and
$||X||_{Li}<K$ then by Lemma A.5 we have:
$$\mu_t(X)\log (t)\leq\mu_t(X)\log (t+e)\leq(1.5)||X||_{Li}<(1.5)K,$$
and so $\tau(e^{-(1.5)K|X|^{-1}})<\infty.$
\end{rems*}

\begin{cor} We have $X\in {Li}_0$
if and only if $\tau(e^{-t|X|^{-1}})<\infty$ for all $t>0$.
\end{cor}

\begin{cor}Let $D$ be an unbounded self-adjoint operator
affiliated with $\mathcal N$. Then\\
 (i) $\tau(e^{-tD^2})<\infty$ for all $t>0$ iff
 $(1+D^2)^{-1}\in {Li}_0$ (i.e., $({\mathcal N}, D)$ is
 $\theta$-summable).\\
(ii)  $\tau(e^{-tD^2})<\infty\text{ for all }t>t_0>0$ if and only if
 $  (1+D^2)^{-1}\in {Li}$ (i.e., $({\mathcal N}, D)$ is weakly
 $\theta$-summable).
\end{cor}

\begin{lemma} If $S, T\in {Li}^q$ then
 $$||\ |S|-|T|\
||_{Li^q}
\leq \sqrt 2 ||S-T||^{1/2}_{Li^q}  ||S+T||^{1/2}_{Li^q}.$$
\end{lemma}

\begin{proof}[\bf Proof]
 This is a special case of Theorem 2.1 of \cite{DD} as $Li^q$
 is fully symmetric by \cite{Suk}.
\end{proof}

\begin{rems*} Fyodor Sukochev has pointed out to us that the following
theorem can be proved in much greater generality using results of
O. E. Tikhonov on the continuity of operator functions in the
measure topology combined with some results of V. Chilin and F. Sukochev
on weak convergence in operator spaces. Using these techniques one is
reduced to proving (in the case of the following theorem) that
$$\int_0^\infty f(\mu_t(T_n))dt\to \int_0^\infty f(\mu_t(T))dt.$$
Since this is the main difficult point in our proof, we prefer not
to confuse the issue with unnecessary generality. However, his comments
have led to a streamlining of our proof for which we are very grateful.
\end{rems*}

\begin{thm}
Let $\mathcal I$ be an invariant operator ideal in $\mathcal N$
and let $f$ be a continuous increasing function from $\real^+$ to itself
such that $f(T)$ is trace-class for each $T\in \mathcal I_+$. Then,
$T\mapsto f(T)$ mapping $\mathcal I_+ \to L^1$ is continuous.
\end{thm}

\begin{proof}[\bf Proof]
Suppose $||T_n-T||_{\mathcal I}\to 0$ in $\mathcal I_+$. Then $S_n=T_n-T$ is
self-adjoint and
$$0\leq T_n=T+S_n\leq T+|S_n|$$
where $|S_n|\geq 0$ is also in $\mathcal I$ with
$$||\ |S_n|\ ||_{\mathcal I}=||S_n||_{\mathcal I}\to 0.$$
Now, since $||.||_{\mathcal I}\geq||.||$ we have $||S_n||\to 0.$
Next, we note that by \cite{FK} (Lemma 2.5 (i), (v)) we have for all $t>0$
\begin{eqnarray}
\mu_t(T_n)\;=\;\mu_t(T+S_n)&\leq&\mu_t(T)+||S_n||\;\;and\nonumber\\
\mu_t(T)\;=\;\mu_t(T_n\!-\!S_n)&\leq& \mu_t(T_n)+||S_n||\nonumber
\end{eqnarray}
Thus, $\mu_t(T_n) \to \mu_t(T)$ uniformly in $t$.
Since the $T_n$ are uniformly bounded in operator norm by, say $C$,
and since $f$ is uniformly continuous on $[0,C]$, we see
that $f(\mu_t(T_n))\to f(\mu_t(T))$ uniformly in $t$.
We claim that
$$\tau(f(T_n))=\int_0^\infty f(\mu_t(T_n))dt\to
\int_0^\infty f(\mu_t(T))dt=\tau(f(T)).$$
The two equalities follow by \cite{FK} (Corollary 2.8).
The convergence
 $$\tau(f(T_n))\to\tau(f(T))$$
 is a little subtle because we have a non-finite measure.
 It suffices to show that every
subsequence has itself a subsequence which converges to
$\tau(f(T))$. So given a subsequence
 $\{\tau(f(T_{n_k}))\}$
choose a further subsequence $T_{n_{k_j}}= T+ S_{n_{k_j}}:= T+R_j$
with the property that $\sum_{j=1}^\infty||R_j||_{\mathcal I}<\infty$.
Then $\sum_{j=1}^\infty |R_j|$ converges in $\mathcal I_+$ to say $R\geq 0$.
So we have $T+R\geq T+R_j$ for all $j$ and also $T+R\geq T$ so that by
\cite{FK} (Lemma 2.5(iii)):
$$f(\mu_t(T+R))\geq f(\mu_t(T_{n_{k_j}}))\ {\text {and }}
f(\mu_t(T+R))\geq f(\mu_t(T)).$$
Now
$$\int_0^\infty f(\mu_t(T+R))dt = \tau(f(T+R))<\infty$$
since $T+R\in \mathcal I$. By dominated convergence
  $\tau(f(T_{n_{k_j}}))\to\tau(f(T))$
and hence
 $\tau(f(T_n))\to\tau(f(T))$.
By the functional calculus,
$$||f(T_n) - f(T)||\to 0$$
as $f$ is bounded and continuous on $[0,C]$ so that using
$$||f(T_n)||_1=\tau(f(T_n))\to
\tau(f(T))= ||f(T)||_1$$
we obtain $||f(T_n) - f(T)||_1\to 0$
by Proposition B.1 of Appendix B.
\end{proof}

\begin{cor}
For any $b>0$, the map $T\mapsto e^{-b|T|^{-1/q}}$ from $Li_0^q$
to the trace-class operators is continuous.
\end{cor}

\begin{proof}[\bf Proof]
 Since $T\mapsto |T|$ is continuous by Lemma B.7 we can
assume that $T\geq 0$. The result follows from Corollary B.5 and
Theorem B.8.
\end{proof}

\begin{cor}
For any $c\geq 0$, $b>0$
 the map $T\mapsto |T|^{-c}e^{-b|T|^{-1/q}}$ from $Li_0^q$ to
the trace-class operators is continuous.
\end{cor}

\begin{proof}[\bf Proof]
 As in the previous proof  we can
assume that $T\geq 0$. Then,
$$T^{-c}e^{-bT^{-1/q}}=T^{-c}e^{-(b/2)T^{-1/q}}e^{-(b/2)T^{-1/q}}.$$
For $t\geq 0$, $t\mapsto t^{-c}e^{-(b/2)t^{-1/q}}$ is bounded and continuous
with the understanding that it is zero at $t=0$.
So, $T\mapsto T^{-c}e^{-(b/2)T^{-1/q}}$
is operator-norm to operator-norm continuous $Li_0^q\to \mathcal N$. Since
the other factor is continuous $Li_0^q\to L^1$ by the previous corollary
and since the $Li_0^q$ norm dominates the operator norm, the product
is also continuous $Li_0^q\to L^1$.
\end{proof}

\begin{cor}
If $c\geq 0$ and $\epsilon>0$, $b>0$ then
 $T\mapsto |T|^{-c}e^{-b|T|^{-(1+\epsilon)}}$ from $Li$ to
the trace-class operators is continuous.
\end{cor}

\begin{proof}[\bf Proof]
Let $q=1/(1+\epsilon)$. It is clear that ${Li}\subset {Li}^q_0$
and the inclusion is continuous by Lemma A.5. The result
then follows from the preceeding corollary.
\end{proof}

\begin{rems*}
If $T\in Li$ and $||T||_{{Li}}<\frac{2}{3}$ then by the Remarks
after Lemma B.4, $\tau(e^{-|T|^{-1}})<\infty.$
The proof of Theorem B.8 now shows that the map
$T\mapsto e^{-|T|^{-1}}$ from\\
$\{T\in {Li}\ |\ ||T||_{{Li}}<\frac{2}{3}\}$ to $L^1$ is continuous.
More generally, for any $C > 0$, the map $T\mapsto e^{-C|T|^{-1}}$ is continuous
from $\{T\in {Li}\ |\ ||T||_{{Li}}<\frac{2}{3}C\}$ to $L^1$.
\end{rems*}

Now we choose a fixed self-adjoint
Breuer-Fredholm operator $F_0\in \mathcal N$ with $1-F_0^2\in {\mathcal I}$
for some invariant operator ideal ${\mathcal I}$,
and recall the definition of Section 2:
$${\mathcal I}_{F_0}=\{X\in {\mathcal I}^{1/2}_{sa}\
|\ 1-(F_0+X)^2\in {\mathcal I}\}$$

\begin{lemma} Suppose that $\mathcal I$ and $\mathcal I^{1/2}$
are invariant operator ideals in $\cn$ satisfying the ``Cauchy-Schwarz''
inequality: $||XY||_{\mathcal I} \leq ||X||_{\mathcal I^{1/2}}
||Y||_{\mathcal I^{1/2}},$ then\\
(1) ${\mathcal I}_{F_0}$
is a real vector space and if $F_1\in F_0+{\mathcal I}_{F_0}$ then
$${\mathcal I}_{F_0}=\{Y\in {\mathcal I}^{1/2}_{sa}\
|\ 1-(F_1+Y)^2\in {\mathcal I}\}$$
so that the definition is independent
of the base point.\\
(2) In the norm $|||X|||_{F_0}=||X||_{\mathcal I^{1/2}}+
||XF_0+F_0X||_{\mathcal I}$,
${\mathcal I}_{F_0}$ is a real Banach space and different base points
define equivalent norms.\\
(3) If $\{F_n\}$ and $F$ are in $F_0 + {\mathcal I}_{F_0}$ then
$|||F_n - F|||\to 0$ if and only if $||F_n -F||_{\mathcal I^{1/2}}\to 0$
and $||(1-F_n^2) - (1-F^2)||_{\mathcal I}\to 0.$
\end{lemma}

\begin{proof}[\bf Proof]
  The first part of the lemma is immediate from the definition
and note (1) on page 678 of \cite{CP1}.
For the second part choose a Cauchy sequence $\{X_n\}$ in
${\mathcal I}_{F_0}$
then there is a limit say $X$ in ${\mathcal I^{1/2}}$ norm.
As $||.||_ {\mathcal I^{1/2}}$ dominates the operator norm,
 $X_n\rightarrow X$
in $\mathcal N$. Similarly there is a limit $Z$ of $F_0X_n+X_nF_0$
in  ${\mathcal I}$ norm and therefore in $\mathcal N$ as well.
Hence $Z=F_0X+XF_0$ is in  ${\mathcal I}$ so that $X_n\rightarrow X$
in the norm on ${\mathcal I}_{F_0}$. To see that the norm is
independent of the base point let $F_1=F_0+Y$ and observe
that
\begin{eqnarray}
|||X|||_{F_1} & = & ||X||_{\mathcal I^{1/2}}+||XF_1+F_1X||_{\mathcal I}
                      \nonumber\\
              &\leq& ||X||_{\mathcal I^{1/2}}+||XF_0+F_0X||_{\mathcal I}
                     +||XY+YX||_{\mathcal I}\nonumber\\
              &\leq& |||X|||_{F_0}+2||XY||_{\mathcal I}\nonumber\\
              &\leq& |||X|||_{F_0}+2||X||_{\mathcal I^{1/2}}
                     ||Y||_{\mathcal I^{1/2}} \nonumber\\
              &\leq& |||X|||_{F_0}(1 + 2||Y||_{\mathcal I^{1/2}}).\nonumber
\end{eqnarray}

The reverse inequality is similar.

The third part of the lemma is immediate.
\end{proof}

\begin{lemma} With ${\mathcal I}$ as in the previous lemma,
the map from ${\mathcal I}_{F_0}\to {\mathcal I}$
given by $X\mapsto 1-(F_0+X)^2$ is continuous.
\end{lemma}

\begin{proof}[\bf Proof]
 We have
\begin{eqnarray}
&    & || 1-({F_0}+X)^2- (1-({F_0}+Y)^2)||_{\mathcal I}\nonumber\\
&\leq& ||Y^2-X^2||_{\mathcal I}+||(Y-X)F_0+F_0(Y-X)||_{\mathcal I}
        \nonumber\\
&\leq& \frac{1}{2}||(Y-X)(Y+X)+(Y+X)(Y-X)||_{\mathcal I}
          +|||Y-X|||_{F_0}\nonumber\\
&\leq& ||Y-X||_{\mathcal I^{1/2}}||Y+X||_{\mathcal I^{1/2}}
         +|||Y-X|||_{F_0}\nonumber\\
&\leq& |||Y-X|||_{F_0}|||Y+X|||_{F_0}+|||Y-X|||_{F_0}.\nonumber
\end{eqnarray}
\end{proof}

\begin{cor} For a self-adjoint $F_0\in \mathcal N$ with
$1-F_0^2\in {Li}^q_0$
and $r\geq 0$ the map\\
$F\mapsto |1-F^2|^{-r}e^{-|1-F^2|^{-1/q}}$ from
the affine space $F_0+(Li_0^q)_{F_0}$ to the trace-class
operators is continuous.
\end{cor}

\begin{proof}[\bf Proof]
 This follows from Lemma B.13 and Corollary B.10.
\end{proof}

\begin{lemma}
Suppose that $\mathcal I$ and $\mathcal I^{1/2}$ are invariant
operator ideals in $\cn$
satisfying the ``Cauchy-Schwarz'' inequality
and let $F_0\in {\cn}_{sa}$ satisfy
$1-F_0^2\in {\mathcal I}.$
If $t\mapsto F_t\in F_0+{\mathcal I}_{F_0}$ is a path, then
it is $C^1$ in that space if and only if:\\
(1) $t\mapsto F_t$ is $C^1$ in ${\mathcal I^{1/2}}$-norm, and\\
(2) $t\mapsto (1-F_t^2)$ is $C^1$ in ${\mathcal I}$-norm.
\end{lemma}

\begin{proof}[\bf Proof]
Suppose conditions (1) and (2) are satisfied, and suppose
$F_t=F_0+X_t$ so that $t\mapsto X_t$ is $C^1$ in ${\mathcal I^{1/2}}$-norm.
Then, $${\mathcal I^{1/2}}-\lim_{s\to t}\frac{1}{s-t}(X_s-X_t)=
X_t^{\prime} \; {exists}$$ and $t\mapsto X_t^{\prime}$ is
${\mathcal I^{1/2}}$-norm continuous
(thus, $X_t^{\prime}$ also exists in operator-norm and is operator-norm
continuous).

Now by the product rule and the ``Cauchy-Schwarz'' inequality, we have that
$t\mapsto X_t^2\in {\mathcal I}$ is $C^1$ in the norm of
$\mathcal I$. Hence,
$$t\mapsto (F_0X_t+X_tF_0)=(1-F_0^2-X_t^2)-(1-F_t^2)$$
is $C^1$ in the norm of $\mathcal I$ by condition (2). Then:
\begin{eqnarray}
Z_t:&=& ||\cdot||_{\mathcal I}\cdot \frac{d}{dt}\left(F_0X_t+
              X_tF_0\right)\nonumber\\
    &=& ||\cdot||\cdot \frac{d}{dt}\left(F_0X_t+X_tF_0\right)\nonumber\\
    &=& \cdots=F_0X_t^{\prime}+X_t^{\prime}F_0.\nonumber
\end{eqnarray}
That is, the difference quotients for $X_t^{\prime}$ also converge to
$X_t^{\prime}$ in the norm of ${\mathcal I}_{F_0}.$

The proof of the other implication is similar and a little easier.
Since we do not use this implication anywhere, we omit the proof.
\end{proof}

\section{\bf INTEGRAL FORMULAE}

In order to explicitly compute the derivatives used
in Section 5 we need to be able to express the map
$T\mapsto |T|^{-r}e^{-|T|^{-1}}$ in terms of the resolvents of $T^2,$
where $T$ is self-adjoint and bounded and $r\geq 0$ (the cases $r=0$
and $r=3/2$ are the ones of interest). In order to do this, we are
forced to consider Cauchy integrals along unbounded contours.
In what follows, we take $\lambda^r$ to be the principal
branch of the usual analytic function of $\lambda$ on
$\comp\backslash(-\infty,0]$.

\begin{lemma}
For $S$ bounded and self-adjoint
and with $\lambda = \pm (t\pm i)$,
$$||(\lambda S-1)^{-1}||\leq (1+t^2)^{1/2}\leq 1+|t| \eqno (C.1)$$
\end{lemma}

\begin{proof}[\bf Proof]
It suffices, by the functional calculus, to prove the numerical
inequality
$|(\lambda x-1)^{-1}|\leq  (1+t^2)^{1/2}$ for all $x\in \real$
which follows by elementary calculus.
\end{proof}

\begin{lemma} For $T$ bounded and self-adjoint
and for any real $a$ in $(0,||T||^{-2})$ we let
 $\sigma=\sigma_1+\sigma_2+\sigma_3$
 be the piecewise
smooth curve in $\comp$ with
$\sigma_1(t) =-t+i, \ t\in (-\infty,-a]$, $\sigma_2(t) =a-ti, \ t\in
 [-1,-1]$
and $\sigma_3(t) =t-i, \ t\in [a,\infty)$. Then
for $c\geq 0$, $k>0$, $b>0$ the integral
$\int_\sigma T^2(\lambda T^2-1)^{-1} \lambda^c e^{-b\lambda^k} d\lambda$
converges absolutely in operator norm.
\end{lemma}

\begin{proof}[\bf Proof]
The contour $\sigma$ is pictured below.

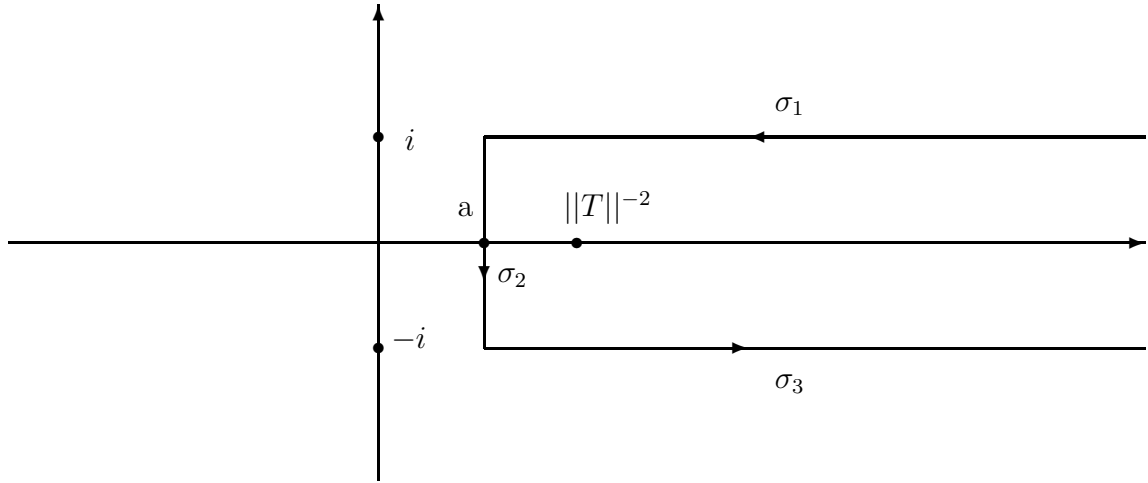
\begin{picture}(450,200)(-150,-100)
\thicklines
\put(0,0){\vector(0,1){90}}
\put(0,0){\line(0,-1){90}}
\put(0,0){\vector(1,0){290}}
\put(0,0){\line(-1,0){140}}
\put(290,40){\vector(-1,0){150}}
\put(140,40){\line(-1,0){100}}
\put(40,40){\vector(0,-1){55}}
\put(40,0){\line(0,-1){40}}
\put(40,-40){\vector(1,0){100}}
\put(140,-40){\line(1,0){150}}
\put(150,50){$\sigma_1$}
\put(150,-55){$\sigma_3$}
\put(45,-15){$\sigma_2$}
\put(40,0){\circle*{4}}
\put(30,10){a}
\put(75,0){\circle*{4}}
\put(70,10){$||T||^{-2}$}
\put(0,40){\circle*{4}}
\put(0,-40){\circle*{4}}
\put(10,35){$i$}
\put(5,-40){$-i$}
\end{picture}

 Now, the integrand is a well-defined continous function
of $\lambda$ as
$\frac{1}{a}>||T^2||$ so that $\frac{1}{a}\notin\text{sp}(T^2)$.
For $\lambda=\sigma_1(t) =-t+i$ (C.1) gives
\begin{eqnarray}
|| T^2(\lambda T^2-1)^{-1} \lambda^c e^{-b\lambda^k} ||&\leq&||T||^2
(1+t^2)^{1/2} |\lambda|^c e^{-\Re (b\lambda^k)}\nonumber\\
&\leq&||T||^2(1+t^2)^{(c+1)/2} e^{-b(1+t^2)^{k/2}/2}\nonumber
\end{eqnarray}
as soon as $|k\arg\lambda|\leq\pi/3$.
This is clearly integrable as $t\to-\infty$. The rest is similar.
\end{proof}

\begin{lemma}
Let $t_0$ in  $[0,||T||^2]$ be fixed
and let $\sigma$ be as above.
Then $$\frac{1}{2\pi i}\int_\sigma t_0(\lambda t_0-1)^{-1}\lambda^c
e^{-b\lambda^k} d\lambda$$
converges absolutely and equals  $t_0^{-c}e^{-b/t_0^k}$ for $t_0>0$
and is zero for $t_0=0$.
\end{lemma}

\begin{proof}[\bf Proof]
We take a cutoff version, $\sigma_N$, of $\sigma$ by truncating
$\sigma_1$ and $\sigma_3$ to
$\sigma_1'(t) =-t+i, \ t\in (-N,-a]$
and $\sigma_3'(t) =t-i, \ t\in [a,N)$ and joining the endpoints. Then
for $\frac{1}{t_0}<N$
we get the value of the integral around $\sigma_N$ as stated using
Cauchy's theorem.
It suffices to show that the integral along the complementary path,
$\sigma_N' = \sigma - \sigma_N,$ converges to $0$ as $N\to \infty$.
The integrals along the horizontal pieces go to zero as $N\to\infty$ using the
previous lemma. The vertical piece where $\lambda(t)=N-ti$ for $t\in[-1,1]$
can be estimated (when $N\geq2/t_0$) using
$\frac{t_0}{|\lambda t_0-1|}\leq t_0.$
This bounds the vertical piece by
$\frac{t_0}{\pi}(N^2+1)^{c/2}e^{-bN^k/2}$ provided
$|k\arg(N\pm i)|\leq \pi/3.$
Hence this piece also goes to zero as $N\to\infty$.
\end{proof}

\begin{lemma} Let $T$ be bounded and self-adjoint and let $\sigma$
 be as in the previous
lemmas. Let $0<a<||T||^{-2}$ then we have for $c\geq 0$ and $k>0$
that the integral
$$\frac{1}{2\pi i}\int_\sigma  T^2(\lambda T^2-1)^{-1}
\lambda^c e^{-b\lambda^k} d\lambda$$
converges absolutely in operator norm to
$ (T^2)^{-c}e^{-bT^{-2k}}$ (which
equals $|T|^{-2c}e^{-b|T|^{-2k}}$).
\end{lemma}

\begin{proof}[\bf Proof]
First let $\{E_x\}$ be the spectral resolution
of $T^2$ (assumed as usual to be strong operator
continuous from the right except at zero where
we take $E_0=0$ while $\lim_{x\to 0^+}E_x$ is the kernel
projection) and let $\xi\in H$. By the Spectral Theorem:
\begin{eqnarray}
\langle (T^2)^{-c}e^{-bT^{-2k}}\xi,\xi\rangle&=&
\int_0^{||T||^2}x^{-c}e^{-b/x^k} d\langle E_x \xi,\xi\rangle\nonumber\\
&=&\int_0^{||T||^2}\{\frac{1}{2\pi i}\int_\sigma x(\lambda
x-1)^{-1}\lambda^ce^{-b\lambda^k} d\lambda\}
 d\langle E_x \xi,\xi\rangle\nonumber\\
&=&\frac{1}{2\pi i}\int_\sigma \{\int_0^{||T||^2}x(\lambda
x-1)^{-1}d\langle E_x \xi,\xi\rangle\}\lambda^ce^{-b\lambda^k} d\lambda
\nonumber\\
&=&\frac{1}{2\pi i}\int_\sigma \langle T^2(\lambda T^2-1)^{-1}\xi,\xi\rangle\
\lambda^c e^{-b\lambda^k} d\lambda\nonumber\\
&=&\langle(\frac{1}{2\pi i}\int_\sigma  T^2(\lambda T^2-1)^{-1}
\lambda^c e^{-b\lambda^k} d\lambda) \xi,\xi\rangle\nonumber
\end{eqnarray}
as the last integral converges absolutely in operator norm.
Observe that the change
in the order of integration is justified for the following reasons.
The measures $|d\lambda|$ and $d\langle E_x\xi,\xi\rangle$ are
positive and $\sigma$-finite and the function
$|x(\lambda x-1)^{-1}\lambda^ce^{-b\lambda^k}|$ is continuous
and hence product measurable.
Moreover the iterated integral
$$\int_0^{||T||^2}\int_\sigma |x(\lambda x-1)^{-1}\lambda^{c}e^{-b\lambda^k}|
\cdot |d\lambda| d\langle E_x \xi,\xi\rangle$$
is clearly finite so an application of the theorems of Tonelli and Fubini
justifies the interchange.
\end{proof}

We now use the notation and results of Appendix B and the integral formula
of the previous lemma to compute the exterior derivatives of the one-forms
of Section 5.

\begin{rems*}
Note that if $g:[0,1]\to L^1$ is continuous in the trace-norm and
differentiable in the operator norm then $f=g^2$ is differentiable
in {\em trace-norm} with $f'(s)=g'(s)g(s)+g(s)g'(s)$
and so $\frac{d}{ds}\tau(f(s))= 2\tau(g(s)g'(s))$. The proof
just uses the usual product rule method.
\end{rems*}

 As usual $F_0$ is
self-adjoint Breuer-Fredholm operator with $1-F_0^2 \in Li_0^q$. Let $X\in
(Li_0^q)_{F_0}$ and let $F_s =F+sX$ so that $F_s^*=F_s$ and hence
$1-F_s^2$ is in $Li_0^q$. Let $T_s=1-F_s^2$ and $T_0=1-F_0^2$.
To differentiate
$s\mapsto |T_s|^{-r}e^{-|T_s|^{-1/q}}$
we factorise it as the square of $s\mapsto |T_s|^{-r/2}e^{-(1/2)|T_s|^{-1/q}}$
and apply the above remark noting that $s\mapsto T_s$ is $Li^q$-norm continuous
so that $s\mapsto |T_s|^{-r/2}e^{-(1/2)|T_s|^{-1/q}}$ is continuous in
trace-norm.

\begin{thm}
Assume the notation above and fix $r\geq 0$ with 
$$g_r(T)=|T|^{-r/2}e^{-(1/2)|T|^{-1/q}}$$ for $T\in Li_0^q$.
 Then in trace norm
$$\frac{d}{ds}|_{s=0}(|1-F_s^2|^{-r}e^{-|1-F_{s}^2|^{-1/q}})=
\frac{d}{ds}|_{s=0}[g_r(T_s)]^2$$
exists and equals
$$\frac{i}{2\pi}\int_\sigma
\left[g_r(T_0) \;,\; R_{\lambda}(T_0) \left[
[F_0,X]_{+} \;,\; T_0 \right]_+R_{\lambda}(T_0)\right]_+ m(\lambda)d\lambda$$
where $T_0=1-F_0^2$ ;  $[\cdot,\cdot]_+$ denotes the anticommutator;
$R_\lambda(T)=(\lambda T^2-1)^{-1}$; and \\
$m(\lambda)=\lambda^{r/4}e^{-(\lambda^{1/2q})/2}.$
\end{thm}

\begin{proof}[\bf Proof]
By the remark it suffices to show that in operator norm
$\frac{d}{ds}|_{s=0}(g_r(T_s))$ exists and
equals
$$\frac{i}{2\pi}\int_\sigma
R_{\lambda}(T_0)\{(F_0X+XF_0)T_0+
T_0(F_0X+XF_0)\}R_{\lambda}(T_0)m(\lambda)d\lambda$$
where $T_s=1-(F_0+sX)^2$. Now,
\begin{eqnarray}
& &\frac{1}{s}(g_r(T_s)-g_r(T_0))\nonumber\\
&=&\frac{1}{s2\pi i}\int_\sigma[T_s^2 R_\lambda(T_s)
-T_0^2 R_\lambda(T_0)]m(\lambda) d\lambda
\nonumber\\
&=&\frac{1}{s2\pi i}\int_\sigma(R_\lambda(T_s)
[T_s^2(\lambda T_0^2-1)-(\lambda T_s^2-1)T_0^2
]R_\lambda(T_0)m(\lambda) d\lambda
\nonumber\\
&=&\frac{1}{s2\pi i}\int_\sigma R_\lambda(T_s)(T_0^2-T_s^2)
R_\lambda(T_0) m(\lambda) d\lambda
\nonumber
\end{eqnarray}
\begin{eqnarray}
&=&\frac{1}{s2\pi i}\int_\sigma R_\lambda(T_s)\frac{1}{2}
[(T_0-T_s)(T_0+T_s)+(T_0+T_s)(T_0-T_s)]
R_\lambda(T_0) m(\lambda) d\lambda\nonumber\\
&=&\frac{-1}{4\pi i}\int_\sigma R_\lambda(T_s)\left[
(F_0X+XF_0+sX^2)\;,\; (T_0+T_s)\right]_+ R_\lambda(T_0)
m(\lambda) d\lambda\hspace{.5in}(C.2)\nonumber
\end{eqnarray}
where again $[\cdot,\cdot]_+$ denotes the anti-commutator.
We need to check that we get convergence
to the formal limit as $s\to 0$ of $(C.2)$:
$$\frac{i}{2\pi}\int_\sigma R_\lambda(T_0)
[(F_0X+XF_0)T_0+T_0(F_0X+XF_0)]R_\lambda(T_0)
m(\lambda) d\lambda. \eqno (C.3)$$
Now let $h_s(\lambda)$ denote the continuous operator-valued
function of $\lambda$ which is the integrand of $(C.2)$
and let $h(\lambda)$ denote the integrand of $(C.3)$.
It suffices to show that $h_s\to h$ uniformly on compact
subsets (of $\sigma$) and that $||h_s(\lambda)||\leq k(\lambda)$
where $\int_\sigma k(\lambda)|d\lambda|<\infty$.
For this to hold it suffices to show that
$$||(\lambda T_s^2-1)^{-1}-(\lambda T_0^2-1)^{-1}||\to 0$$
uniformly on compacta as the other factors converge uniformly on all of
$\sigma$.
Now, by the resolvent equation and Lemma C.1,
\begin{eqnarray}
||(\lambda T_s^2-1)^{-1}-(\lambda T_0^2-1)^{-1}||
&=&
||(\lambda T_s^2-1)^{-1}\lambda(T_0^2-T_s^2)(\lambda T_0^2-1)^{-1}||
\nonumber\\
&\leq& |\lambda|^3||T_0^2-T_s^2||\to 0\nonumber
\end{eqnarray}
uniformly for $\lambda$ in a bounded set.
Let
$$C=\sup_{s\in[-1,1]}\frac{1}{2}||[(F_0X+XF_0+sX^2)(T_0+T_s)+
(T_0+T_s)(F_0X+XF_0+sX^2)]||.$$
Then
$$||h_s(\lambda)||\leq C|\lambda|^{2+r/4}e^{-\Re(\lambda^{1/2q})/2}$$
which is integrable on $\sigma$.
\end{proof}


\begin{thebibliography}{APS3}

\bibitem[APS1]{APS1} M. F. Atiyah, V. Patodi, I. M. Singer, \emph{Spectral
Asymmetry and Riemannian Geometry.} \emph{I}, Proc. Camb. Phil. Soc.,
{\bfseries 77}(1975), 43--69.

\bibitem[APS3]{APS3} M. F. Atiyah, V. Patodi, I. M. Singer, \emph{Spectral
Asymmetry and Riemannian Geometry.} \emph{III}, Proc. Camb. Phil. Soc.,
{\bfseries 79}(1976), 71--99.

\bibitem[AS]{AS} M. F. Atiyah, I. M. Singer, \emph{Index Theory for
Skew-adjoint Fredholm Operators}, Publ. Math. Inst. Hautes Etudes Sci. (Paris),
{\bf series 37}(1969), 5--26.

\bibitem[ASS]{ASS} J. Avron, R. Seiler, B. Simon, \emph{ The
index of a pair of projections}, J. Funct. Analysis {\bf 120} (1994), 220--237.

\bibitem[BJ]{BJ} S. Baaj, P. Julg, \emph{Th\'{e}orie bivariante de Kasparov
et op\'{e}rateurs non born\'{e}s dans les $C^*$-modules hilbertiens},
C. R. Acad. Sci. Paris S\'{e}r. {\bf I} 296(1983),no.21,875--878.

\bibitem[BF]{BF} B. Boo$\beta$-Bavnbek, K. Furutani, \emph{The Maslov Index: A
Functional Analytic Definition}  \emph{and the Spectral Flow Formula},
Tokyo J. Math., {\bf 21},No. 1(1998), 1--34.

\bibitem[BW]{BW} B.Boo$\beta$-Bavnbek, K.P.Wojciechowski, Elliptic Boundary
Problems for Dirac Operators, Birkh\"{a}user, Boston, Basel, Berlin, 1993.

\bibitem[B1]{B1} M. Breuer, \emph{Fredholm Theories in von Neumann
algebras. I}, Math. Ann., {\bfseries 178}(1968), 243--254.

\bibitem[B2]{B2} M. Breuer, \emph{Fredholm Theories in von Neumann
algebras. II}, Math. Ann., {\bfseries 180}(1969), 313--325.

\bibitem[CDSS]{CDSS}  L.A. Coburn, R.G. Douglas, D.G. Schaeffer and I.M. Singer
\emph{C$^*$-algebras of operators on a half-space II: Index Theory}
Publ. IHES, {\bfseries 71} (1971) 69-79 

\bibitem[CDS]{CDS} V.I. Chilin, P.G.Dodds, F.A. Sukochev, \emph{The
Kadec-Klee Property in Symmetric Spaces of Measureable Operators},
Israel J. Math., {\bf vol. 97}(1997), 203--219.

\bibitem[CP]{CP} A. L. Carey, J. Phillips, \emph{Algebras Almost Commuting
with Clifford Algebras in a II$_\infty$ Factor}, $K$-Theory,
{\bfseries 4}(1991), 445--478.

\bibitem[CP1]{CP1} A. L. Carey, J. Phillips, \emph{Unbounded Fredholm
Modules and Spectral Flow}, Canadian J. Math., {\bf vol. 50}(4)(1998),
673--718.

\bibitem[CPRS2]{CPRS2} A.L. Carey, J. Phillips, A. Rennie, F.A. Sukochev,
paper in preparation.

\bibitem[CPS]{CPS} A. L. Carey, J. Phillips and F. A. Sukochev \emph{On
unbounded p-summable Fredholm modules}, Advances in Math., {\bf 151} (2000),
140-163.

\bibitem[CPSu]{CPSu} A. L. Carey, J. Phillips and F. A. Sukochev \emph
{Spectral flow and Dixmier traces}, Advances in Math., {\bf 173} (2003),
68-113.

\bibitem[Co1]{Co1} A. Connes, Noncommutative Differential Geometry, Publ.
Math. Inst. Hautes Etudes Sci. (Paris), {\bfseries 62}(1985), 41--44.

\bibitem[Co2]{Co2} A. Connes, \emph{Cyclic Cohomology of Banach Algebras and
Characters of $\theta$-summable Fredholm Modules}, $K$-Theory {\bf1}(1988),
519--548.

\bibitem[Co3]{Co3} A. Connes, \emph{Compact Metric Spaces, Fredholm Modules
and Hyperfiniteness}, Ergodic Theory and Dynamical Systems {\bf9}(1989),
207--220.

\bibitem[Co4]{Co4} A. Connes, Non-commutative geometry, Academic Press,
San Diego, 1994.

\bibitem[CoM]{CoM} A. Connes and H. Moscovici, \emph{ Cyclic cohomology,
the Novikov conjecture and hyperbolic groups},
Topology {\bfseries 29} (1990), 345--388.

\bibitem[CoMo]{CoMo} A. Connes and H. Moscovici, \emph 
{The Local Index Formula in Noncommutative Geometry} GAFA {\bf 5} 
(1995) 174-243

\bibitem[DD] {DD} P.G. Dodds and T.K. Dodds,
\emph{On a submajorization inequality of
T. Ando,} Operator Theory Advances and Applications {\bfseries 75}
 (1995), 113--133.

\bibitem [DHK]{DHK} R. G. Douglas, S. Hurder and J. Kaminker., \emph{Cyclic
cocycles, renormalisation and eta-invariants}, Invent. Math. {\bf 103}
(1991) 101--179.

\bibitem[Dix]{Dix} J. Dixmier, Les alg\`ebres d'op\'erateurs dans l'espace
Hilbertien (Alg\`ebres de von Neumann), Gauthier-Villars, Paris, 1969.

\bibitem[DS]{DS} N. Dunford, J. T. Schwartz,
Linear Operators, Part II, Wiley, New York, London, 1963.

\bibitem[DPSS]{DPSS} P. G. Dodds, B. de Pagter, E. M. Semenov, F. A. Sukochev
\emph{Symmetric functionals and singular traces},
Positivity {\bf 2} (1998), 47-75.

\bibitem[FK]{FK} T. Fack and H. Kosaki \emph{ Generalised $s$-numbers of
$\tau$-measurable operators} Pacific J. Math. {\bf 123} (1986), 269--300

\bibitem[G]{G} E. Getzler, \emph{The Odd Chern Character in Cyclic Homology
and Spectral Flow}, Topology, {\bfseries 32}(1993), 489--507.

\bibitem[GK]{GK} I. C. Gohberg, M. G. Krein, Introduction to the Theory
of Non-selfadjoint Operators, Translations of Mathematical Monographs,
{\bf vol. 18}, AMS, 1969.

\bibitem[GZ]{GZ} E. Getzler and A. Szenes, \emph{On the Chern character of a 
theta-summable Fredholm module}, J. Functional Anal. {\bf 84} (1989), 343-357.

\bibitem[H]{H} S. Hurder, \emph{Eta invariants and the odd index theorem for
coverings}, Contemp. Math., {\bf 105} (1990), 47--82.

\bibitem[JLO]{JLO} A. Jaffe, A. Lesniewski and K. Osterwalder, \emph{Quantum
K-theory: the Chern character}, Comm. Math. Phys. {\bf 112} (1988), 75-88.

\bibitem[Kam]{Kam} J. Kaminker, \emph{Operator algebraic invariants for
elliptic operators}, Proc. Symp. in Pure Math., {\bf 51} (1990), 307--314.

\bibitem[L]{L} S. Lang, Analysis I, Addison-Wesley, Reading, Menlo Park,
London, Don Mills, 1968.

\bibitem[Ma]{Ma} V. Mathai, \emph{$L^2$ invariants of covering spaces}
in Geometric Analysis and Lie Theory in Mathematics and Physics,
Cambridge University Press, Cambridge 1998.

\bibitem[M]{M} V. Mathai, \emph{Spectral flow, eta invariants
 and von Neumann algebras},
Journal of Functional Analysis, {\bfseries 109} (1992), 442--456.

\bibitem[M1]{M1} V. Mathai. preprint and private communication.

\bibitem[P1]{P1} V.S. Perera, \emph{Real Valued Spectral Flow in a Type
II$_\infty$ Factor}, Ph.D. Thesis, IUPUI, 1993.

\bibitem[P2]{P2} V.S. Perera, \emph{Real Valued Spectral Flow in a Type
II$_\infty$ Factor}, preprint, IUPUI, 1993.

\bibitem[Ph]{Ph} J. Phillips, \emph{Self-Adjoint Fredholm Operators and
Spectral Flow}, Canad. Math. Bull., {\bf 39}(1996), 460--467.

\bibitem[Ph1]{Ph1} J. Phillips,
\emph{Spectral Flow in Type I and Type II factors-a New Approach},
Fields Institute Communications, {\bf vol. 17}(1997), 137--153.

\bibitem[Ru]{Ru} W. Rudin, Principles of Mathematical Analysis, $3^{rd}$ ed.,
McGraw-Hill, New York, London, Toronto, 1976.

\bibitem[Si]{Si} I.M. Singer, \emph{Eigenvalues of the Laplacian and invariants
of manifolds}, Proceedings of the International Congress, Vancouver 1974,
{\bf vol. I}, 187-200.

\bibitem[Sp]{Sp} M. Spivak, A Comprehensive Introduction to Differential
Geometry, {\bf vol. 1}, $2^{nd}$ ed., Publish or Perish Inc., Berkeley, 1979.

\bibitem[Suk]{Suk} F. A. Sukochev, \emph{Operator esimates for Fredholm
modules}, preprint.

\bibitem[Suk2]{Suk2} F. A. Sukochev, \emph{Unbounded Fredholm Modules and
Submajorization}, preprint.

\bibitem[Y]{Y} K. Yosida, Functional Analysis, Springer, New York, 1971.

\end{thebibliography}
\end{document}